\documentclass[reqno,12pt]{amsart}
\usepackage{amsmath,amsthm,amssymb,amsfonts,amscd,graphicx}
\theoremstyle{plain} 
	\newtheorem{thm}{Theorem}[section]
	\newtheorem*{thm*}{Theorem}
	\newtheorem{cor}[thm]{Corollary}
	\newtheorem*{cor*}{Corollary}
	\newtheorem{lem}[thm]{Lemma}
	
	\newtheorem{prop}[thm]{Proposition}
	
	\newtheorem*{conj*}{Conjecture}
	\newtheorem{assum}[thm]{Assumption}
\theoremstyle{definition}
	\newtheorem{defn}[thm]{Definition}%[section]

\theoremstyle{remark}
	\newtheorem{rem}[thm]{Remark}
	
	\newtheorem*{pf}{Proof}
\numberwithin{equation}{section}
\def\AA{{\mathbb A}}
\def\CC{{\mathbb C}}

\def\PP{{\mathbb P}}
\def\QQ{{\mathbb Q}}
\def\RR{{\mathbb R}}
\def\ZZ{{\mathbb Z}}

\def\C{{\mathcal C}}

\def\F{{\mathcal F}}
\def\G{{\mathcal G}}
\def\H{{\mathcal H}}

\def\O{{\mathcal O}}

\def\T{{\mathcal T}}

\def\m{{\mathfrak m}}
\def\p{\partial }

\def\ns{{\nabla}\hspace{-1.4mm}\raisebox{0.3mm}{\text{\footnotesize{\bf /}}}}

\begin{document}
\title{Primitive Forms for Affine Cusp Polynomials}
\date{\today}
\author{Yoshihisa Ishibashi}
\address{Department of Mathematics, Graduate School of Science, Osaka University, 
Toyonaka Osaka, 560-0043, Japan}
\email{y-ishibashi@cr.math.sci.osaka-u.ac.jp}
\author{Yuuki Shiraishi}
\address{Department of Mathematics, Graduate School of Science, Osaka University, 
Toyonaka Osaka, 560-0043, Japan}
\email{sm5021sy@ecs.cmc.osaka-u.ac.jp}
\author{Atsushi Takahashi}
\address{Department of Mathematics, Graduate School of Science, Osaka University, 
Toyonaka Osaka, 560-0043, Japan}
\email{takahashi@math.sci.osaka-u.ac.jp}
\begin{abstract}
We determine a primitive form for a universal unfolding of an affine cusp polynomial.
Moreover, we prove that the resulting Frobenius manifold is isomorphic to 
the one constructed from the Gromov-Witten theory for an orbifold projective line with at most three orbifold points. 
\end{abstract}
\maketitle
%%%%%%%%%%%%%%%%%%%%%%%%%%%%%%%%%%%%%%%%%%%%%%%%%%%%%%%%%%%%%%%%
%%%%%%%%%%%%%%%%%%%%%%%%%%%%%%%%%%%%%%%%%%%%%%%%%%%%%%%%%%%%%%%%
\section{Introduction}
The notion of primitive forms is introduce by K. Saito \cite{sa:1, S1202-Saito} in his study of period mapping 
associated to a deformation theory of an isolated hypersurface singularity.
Roughly speaking,  a primitive form is a generalization of a differential of the first kind on an elliptic curve and 
the tools in order to define a primitive form are natural generalization of a polarized Hodge structure.
One of the most important ingredients is a Frobenius structure, which he called the {\it flat structure}, 
on the deformation space. 
The existence of primitive forms is proved for isolated hypersurface singularities \cite{M.Saito} (for comprehensive guides, see also \cite{her:1, sabb:1}).
However, explicit forms of primitive forms are rarely known except for a few cases; for ADE singularities \cite{sa:1},  
for simple elliptic singularities \cite{sa:1, sat:1} and for some Laurent polynomials \cite{douai-sabbah:1, douai-sabbah:2}.
In the present paper, we shall add a new class of polynomials to them.
Let $A$ be a triplet $(a_1,a_2,a_3)$ of positive integers such that $a_1\le a_2\le a_3$.
Set $\mu_A=a_1+a_2+a_3-1$ and $\chi_A:=1/a_1+1/a_2+1/a_3-1$.
We can associate to $A$ the polynomial $f_A({\bf x})\in\CC[x_1,x_2,x_3]$ given as
\begin{equation}
f_A({\bf x}):=x_1^{a_1}+x_2^{a_2}+x_3^{a_3}-q^{-1}\cdot x_1 x_2 x_3
\end{equation}
for some $q\in\CC\backslash\{0\}$, which we shall call a cusp polynomial of type $A$, 
in particular, an affine cusp polynomial of type $A$ if $\chi_A>0$.
The purpose of this paper is to determine a primitive form for a universal unfolding $F_A$ 
of an affine cusp polynomial $f_A$, a deformation $F_A$ of $f_A$ over a $\mu_A$-dimensional 
manifold $M:=\CC^{\mu_A-1}\times (\CC\backslash\{0\})$ given as
\begin{equation}
F_A({\bf x};{\bf s},s_{\mu_A}):=x_1^{a_1}+x_2^{a_2}+x_3^{a_3}-s_{\mu_A}^{-1}\cdot x_1 x_2 x_3
+s_1\cdot 1+\sum_{i=1}^{3} \sum_{j=1}^{a_i-1}s_{i,j}\cdot x_i^j
\end{equation}
which induces an $\O_M$-isomorphism
\begin{equation}\label{eqn:isom in intro}
\T_M\longrightarrow \O_M[x_1,x_2,x_3]
\left/\Big(\frac{\p F_A}{\p x_1},\frac{\p F_A}{\p x_2},\frac{\p F_A}{\p x_3}\Big)\right.
,\quad {\delta}\mapsto{\delta}F_A.
\end{equation}
Moreover, we shall prove that the resulting Frobenius manifold is isomorphic to 
the one constructed from the Gromov-Witten theory for an orbifold projective line with at most three orbifold points. 
First, we shall show the following theorem, which is our main theorem in this paper:
\begin{thm*}[Theorem~\ref{primitive form}]
The element $\zeta_A:=[s_{\mu_A}^{-1}dx_1 \wedge dx_2 \wedge dx_3]\in \H_{F_A}^{(0)}$ is a primitive form 
for the tuple $(\H_{F_A}^{(0)},\nabla, K_{F_A})$ with the minimal exponent $r=1$. 
\end{thm*}
There is an $\O_M$-algebra structure on $\T_M$ by the $\O_M$-isomorphism \eqref{eqn:isom in intro}, 
which we shall denote by $\circ$. Note also that we have special elements $[1]$ and $[F_A]$ on the right hand side
of \eqref{eqn:isom in intro} which give elements $e\in\Gamma(M,\T_M)$ which is the identity with respect to the product $\circ$
and $E\in\Gamma(M,\T_M)$. 
As a corollary of our main theorem, we obtain a non-degenerate symmetric $\O_M$-bilinear form on $\T_M$
such that the tuple $(\eta,\circ,e, E)$ defines on the deformation space $M$ 
a Frobenius structure of rank $\mu_A$ and dimension one (cf. \cite{st:1}).  
Next, we show that this Frobenius structure satisfies the conditions in Theorem 3.1 of \cite{ist:1}.
\begin{thm*}[Theorem~\ref{thm:satisfies IST1}]
For the Frobenius structure of rank $\mu_A$ and dimension one constructed from the pair $(f_A,\zeta_A)$, 
there exists flat coordinates $t_1$, $t_{1,1},\dots ,t_{3,a_3-1}$, $t_{\mu_A}$ satisfying the following conditions$:$
\begin{enumerate}
\item 
The unit vector field $e$ and the Euler vector field $E$ are given by
\[
e=\frac{\p}{\p t_1},\ E=t_1\frac{\p}{\p t_1}+\sum_{i=1}^3\sum_{j=1}^{a_i-1}\frac{a_i-j}{a_i}t_{i,j}\frac{\p}{\p t_{i,j}}
+\chi_A\frac{\p}{\p t_{\mu_A}}.
\]
\item 
The non-degenerate symmetric bilinear form $\eta$ on $\T_M$ defined by 
\begin{equation}
\eta(\delta,\delta'):=K_{F_A}^{(0)}(u\nabla_\delta\zeta_A,u\nabla_{\delta'}\zeta),\quad \delta,\delta'\in\T_M,
\end{equation}
satisfies
\begin{align*}
&\ \eta\left(\frac{\p}{\p t_1}, \frac{\p}{\p t_{\mu_A}}\right)=
\eta\left(\frac{\p}{\p t_{\mu_A}}, \frac{\p}{\p t_1}\right)=1,\\ 
&\ \eta\left(\frac{\p}{\p t_{i_1,j_1}}, \frac{\p}{\p t_{i_2,j_2}}\right)=
\begin{cases}
\frac{1}{a_{i_1}}\quad i_1=i_2\text{ and }j_2=a_{i_1}-j_1,\\
0 \quad \text{otherwise}.
\end{cases}
\end{align*}
\item 
The Frobenius potential $\F_{f_A,\zeta_A}$ satisfies $E\F_{f_A,\zeta_A}|_{t_{1}=0}=2\F_{f_A,\zeta_A}|_{t_{1}=0}$,
\[
\left.\F_{f_A,\zeta_A}\right|_{t_1=0}\in\CC\left[t_{1,1}, \dots, t_{1,a_1-1},
t_{2,1}, \dots, t_{2,a_2-1},t_{3,1}, \dots, t_{3,a_3-1},e^{t_{\mu_A}}\right].
\]
\item Assume the condition {\rm (iii)}. The restriction of the Frobenius potential $\F_{f_A,\zeta_A}$ to 
the submanifold $\{t_1=e^{t_{\mu_A}}=0\}$ is given as 
\[
\F_{f_A,\zeta_A}|_{t_1=e^{t_{\mu_A}}=0}=\G^{(1)}+\G^{(2)}+\G^{(3)}, 
\]
where $\G^{(i)}\in \CC[t_{i,1},\dots, t_{i,a_i-1}]$, $i=1,2,3$.
\item 
Assume the condition {\rm (iii)}. In the frame $\frac{\p}{\p t_1}$, $\frac{\p}{\p t_{1,1}},\dots, 
\frac{\p}{\p t_{3,a_3-1}}$, $\frac{\p}{\p t_{\mu_A}}$ of $\T_M$,
the product $\circ$ can be extended to the limit $t_1=t_{1,1}=\dots=t_{3,a_3-1}=e^{t_{\mu_A}}=0$.
The $\CC$-algebra obtained in this limit is isomorphic to 
\[
\CC[x_1,x_2,x_3]\left/\left(x_1x_2,x_2x_3,x_3x_1,a_1x_1^{a_1}-a_2x_2^{a_2},
a_2x_2^{a_2}-a_3x_3^{a_3}\right)\right.,
\]
where $\p/\p t_{1,1},\p/\p t_{2,1},\p/\p t_{3,1}$ are mapped to
$x_1,x_2,x_3$, respectively.
\item The term 
\[
\begin{cases}
e^{t_{\mu_A}}\quad \textit{if}\quad a_{1}=a_{2}=a_{3}=1,\\
t_{3,1}e^{t_{\mu_A}}\quad \textit{if}\quad 1=a_{1}=a_{2}<a_{3},\\
t_{2,1}t_{3,1}e^{t_{\mu_A}}\quad \textit{if}\quad 1=a_{1}<a_{2},\\
t_{1,1}t_{2,1}t_{3,1}e^{t_{\mu_A}}\quad \textit{if}\quad a_{1}=2,
\end{cases}
\]
occurs with the coeffcient $1$ in $\F_{f_A,\zeta_A}$. 
\end{enumerate}
\end{thm*}
As a consequence, 
we obtain the mirror isomorphism as Frobenius manifolds:
\begin{cor*}[Corollary~\ref{cor:classical mirror}]
There exists an isomorphism of Frobenius manifolds between the one constructed from the Gromov--Witten theory 
for $\PP^{1}_A$ and the one constructed from the pair $(f_A,\zeta_A)$. 
\end{cor*}
This result is already obtained by  
Milanov--Tseng \cite{mt:1} for the case $a_1=1$ and by Rossi \cite{r:1} for the case $\chi_A>0$. 
The present work is the second in a series of papers to simplify and generalise their works. 
In our subsequent papers, we shall show that the Frobenius structure constructed from
the pair $(f_A,\zeta_A)$ in this paper is isomorphic to the one constructed from the extended affine 
Weyl group $\widehat{W}_A$ by Dubrovin--Zhang in \cite{dz:1} 
and that the Frobenius structure constructed from cusp singularity (i.e., $f_A$ with $\chi_A<0$) 
with the ``canonical" primitive form satisfies the conditions in Theorem 3.1 of \cite{ist:1}.
\bigskip
\noindent
{\it Acknowledgement}\\
\indent
The second named author is deeply grateful to Claus Hertling and 
Christian Sevenheck for their valuable discussions and encouragement.  
He is supported by the JSPS International Training Program (ITP).
The third named author is supported by Grant-in Aid for Scientific Research 
grant numbers 24684005 from the Ministry of Education, 
Culture, Sports, Science and Technology, Japan.
%%%%%%%%%%%%%%%%%%%%%%%%%%%%%%%%%%%%%%%%%%%%%%%%%%%%%%%%%%%%%%%%%%%%%%%%%%%%%%
%%%%%%%%%%%%%%%%%%%%%%%%%%%%%%%%%%%%%%%%%%%%%%%%%%%%%%%%%%%%%%%%%%%%%%%%%%%%%%
\section{Notations and Terminologies}
Let $A$ be a triplet $(a_1,a_2,a_3)$ of positive integers such that $a_1\le a_2\le a_3$.
Set 
\begin{equation}
\mu_A:=a_1+a_2+a_3-1
\end{equation}
and
\begin{equation}
\chi_A:=\frac{1}{a_1}+\frac{1}{a_2}+\frac{1}{a_3}-1.
\end{equation}
%%%%%%%%%%%%%%%%%%%%%%%%%%%%%%%%%%%%%%%%%%%%%%%%%%%%%%%%%%%%%%%%%%%%%%%%%%%%%%
\subsection{Universal unfolding of affine cusp polynomials}
\begin{defn}
A polynomial $f_A({\bf x})\in\CC[x_1,x_2,x_3]$ given as
\begin{equation}
f_A({\bf x}):=x_1^{a_1}+x_2^{a_2}+x_3^{a_3}-q^{-1}\cdot x_1 x_2 x_3
\end{equation}
for some $q\in\CC\backslash\{0\}$ is called 
the {\it cusp polynomial of type $A$}. 
In particular, we shall call $f_A({\bf x})$ an {\it affine cusp polynomial} of type $A$ if $\chi_A>0$. 
\end{defn}
\begin{assum}
Throughout this paper, we shall assume that the triplet of integers $A$ satisfies $\chi_A>0$.
In particular, we shall only consider affine cusp polynomials.
\end{assum}
Let $f_A=f_A({\bf x})$ be an affine cusp polynomial.
One can easily show that 
\[
\CC[x_1,x_2,x_3]\left/\Big(\frac{\p f_A}{\p x_1},\frac{\p f_A}{\p x_2},\frac{\p f_A}{\p x_3}\Big)\right.
\]
is of rank $\mu_A$.
We can consider the {\it universal unfolding} of $f_A$,  
a deformation $F_A$ of $f_A$ defined on $\CC^{3}\times M$, 
$M:=\CC^{\mu_A-1}\times (\CC\backslash\{0\})$ 
over a $\mu_A$-dimensional parameters $({\bf s},s_{\mu_A})\in M$ given as follows:
\begin{defn}
Define a function $F_A({\bf x};{\bf s}, s_{\mu_A})$ defined on $\CC^{3}\times M$ as follows$;$ 
\begin{equation}
F_A({\bf x};{\bf s},s_{\mu_A}):=x_1^{a_1}+x_2^{a_2}+x_3^{a_3}-s_{\mu_A}^{-1}\cdot x_1 x_2 x_3
+s_1\cdot 1+\sum_{i=1}^{3} \sum_{j=1}^{a_i-1}s_{i,j}\cdot x_i^j.
\end{equation}
\end{defn}
Denote by
\[
p:\CC^{3}\times M\longrightarrow M,\quad
({\bf x};{\bf s},s_{\mu_A})\mapsto ({\bf s},s_{\mu_A})
\]
the projection map from the total space to the deformation space.
Set 
\begin{equation}
p_*\O_\C:=\O_{M}[x_1,x_2,x_3]
\left/\Big(\frac{\p F_A}{\p x_1},\frac{\p F_A}{\p x_2},\frac{\p F_A}{\p x_3}\Big)\right..
\end{equation}
$p_*\O_\C$ can be thought of as the direct image of 
the sheaf of relative {\it algebraic} functions on the relative critical set $\C$ of $F_A$ with respect to the projection 
$p:\CC^3\times M\longrightarrow M$.
\begin{prop}\label{unfolding}
The function $F_A({\bf x};{\bf s},s_{\mu_A})$ satisfies the following conditions$:$
\begin{enumerate}
\item $F_A({\bf x};{\bf 0},q)=f_A({\bf x})$.
\item The $\O_M$-homomorphism $\rho$ called the Kodaira--Spencer map defined as 
\begin{equation}\label{universal}
\rho: \T_M\longrightarrow p_*\O_\C,\quad \delta\mapsto \delta F_A,
\end{equation}
is an isomorphism.
\end{enumerate}
\end{prop}
\begin{pf}
Since $F_A$ is tame at any point on the deformation space $M$,
namely, there are no critical points coming from infinity.
Therefore, the statement follows from the same argument in
Section 2.a of \cite{douai-sabbah:1}.
\qed
\end{pf}
Note that the tangent bundle $\T_M$ naturally obtains an $\O_M$-algebra structure. 
\begin{defn}\label{defn:Jacobi ring}
We shall denote by $\circ$ the induced product structure on $\T_M$ by the $\O_M$-isomorphism~\eqref{universal}. 
Namely, for $\delta,\delta'\in \T_M$, we have
\begin{equation}
{(\delta\circ\delta')}F_A = 
{\delta}F_A\cdot{\delta'}F_A\ \text{in}\ p_*\O_\C.
\end{equation}
\end{defn}
Set $\overline{M}:=\CC^{\mu_A-1} \times \CC$. 
Let $\iota$ be the natural embedding
\begin{equation}
\iota:M\hookrightarrow \overline{M}, \quad
({\bf s}, s_{\mu_A}) \mapsto ({\bf s}, s_{\mu_A}),
\end{equation}
and let $M_\infty$ be the divisor $\overline{M}\backslash\iota(M):=\{ ({\bf s},s_{\mu_A})\in\overline{M}~|~s_{\mu_A}=0\}$.
Denote again by $p$ the projection map from the total space to the deformation space:
\[
p:\CC^{3}\times \overline{M}\longrightarrow \overline{M},\quad
({\bf x};{\bf s},s_{\mu_A})\mapsto ({\bf s},s_{\mu_A})
\]
Set 
\begin{equation}\label{eq:Jacobi-ring extension}
p_*\O_{\overline{\C}}:=\O_{\overline{M}}[x_1,x_2,x_3]
\left/\Big(s_{\mu_A}\frac{\p F_A}{\p x_1},s_{\mu_A}\frac{\p F_A}{\p x_2},s_{\mu_A}\frac{\p F_A}{\p x_3},
H_{1}({\bf x}, {\bf s}), H_{2}({\bf x}, {\bf s})
\Big)\right..
\end{equation}
where
\begin{equation}
H_{i}({\bf x}, {\bf s}):=a_{i}x_{i}^{a_i}-a_{i+1}x_{i+1}^{a_{i+1}}+
\sum_{j=1}^{a_i-1}j\cdot s_{i,j}\cdot x_i^j-\sum_{j=1}^{a_{i+1}-1}j\cdot s_{i+1,j}\cdot x_{i+1}^j,\ i=1,2.
\end{equation}
\begin{prop}
As an $\O_{\overline{M}}$-module, ${p}_*\O_{\overline{\C}}$ is free of rank $\mu_A$.
As an $\O_{\overline{M}}$-algebra, it is isomorphic to the $\O_{\overline{M}}$-subalgebra of $\iota_*\left(p_*\O_\C\right)$ 
whose $\O_{\overline{M}}$-basis is given by the set of residue classes of monomials 
\begin{equation}
\left\{1,\ x_i^j\  (i=1,2,3,\ j=1,\dots, a_i-1),\ s_{\mu_A}^{-1}x_1x_2x_3\right\}.
\end{equation}
In particular, we have the isomorphism of $\O_M$-modules
\begin{equation}
\iota^*\left({p}_*\O_{\overline{\C}}\right)\simeq p_*\O_\C.
\end{equation}
Namely, ${p}_*\O_{\overline{\C}}$ is an extension of $p_*\O_\C$ as a free $\O_{\overline{M}}$-module.
\end{prop}
\begin{pf}
It is almost obvious that $\iota^*\left({p}_*\O_{\overline{\C}}\right)\simeq p_*\O_\C$.
\begin{lem}\label{lem:limit of Jac}
We have
\begin{equation}
{p}_*\O_{\overline{\C}}\left/s_{\mu_A}{p}_*\O_{\overline{\C}}\right.
\simeq\O_{M_\infty}[x_1,x_2,x_3]\left/\Big(x_2x_3,x_3x_1,x_1x_2,H_{1}({\bf x}, {\bf s}), H_{2}({\bf x}, {\bf s})\Big)\right.,
\end{equation}
is a free $\O_{M_\infty}$-module of rank $\mu_A$. In particular, we have
\begin{equation}
{p}_*\O_{\overline{\C}}\left/\m_{({\bf 0},0)}{p}_*\O_{\overline{\C}}\right.
\simeq\CC[x_1,x_2,x_3]\left/\Big(x_2x_3,x_3x_1,x_1x_2,a_1x_1^{a_1}-a_2x_2^{a_2},a_2x_2^{a_2}-a_3x_3^{a_3}\Big)\right.,
\end{equation}
where $\m_{({\bf 0},0)}$ is the maximal ideal of $\O_{\overline{M}}$ corresponding to the point $({\bf 0},0)$.
\end{lem}
\begin{pf}
Some elementary calculations yield the statement.
\qed
\end{pf}
By this lemma, we see that ${p}_*\O_{\overline{\C}}$ is free of rank $\mu_A$ as an $\O_{\overline{M}}$-module.
The rest follows from the following equality in $\iota_*\left(p_*\O_\C\right):$
\begin{equation}
\left[s_{\mu_A}^{-1}x_1x_2x_3\right]=\left[a_ix_i^{a_{i}}+\sum_{j=1}^{a_i-1}j\cdot s_{i,j}\cdot x_i^j\right],\quad i=1,2,3.
\end{equation}
\qed
\end{pf}
Therefore, $p_*\O_{\overline{\C}}$ can be thought of as the direct image of 
the sheaf of relative algebraic functions on $\overline{C}$, 
a closure of the relative critical set $\C$ in $\CC^3\times \overline{M}$, which is flat and finite 
over $\overline{M}$.
\begin{prop}
Denote by $\T_{\overline{M}}\left(-\log M_\infty\right)$ the sheaf of holomorphic vector fields on $\overline{M}$ 
with logarithmic zeros along the divisor $M_\infty$ defined as 
\begin{equation}
\T_{\overline{M}}\left(-\log M_\infty\right):=\{\delta\in\T_{\overline{M}}~|~\delta s_{\mu_A}\in s_{\mu_A}\cdot \O_{\overline{M}}\}.
\end{equation}
Then, $\T_{\overline{M}}\left(-\log M_\infty\right)$ is a free $\O_{\overline{M}}$-module of rank $\mu_A$ 
and we have the isomorphism of $\O_M$-modules
\begin{equation}
\iota^*\left(\T_{\overline{M}}\left(-\log M_\infty\right)\right)\simeq \T_M.
\end{equation}
Namely, $\T_{\overline{M}}\left(-\log M_\infty\right)$ is an extension of $\T_M$ as a free $\O_{\overline{M}}$-module.
\end{prop}
\begin{pf}
The statement follows from the fact that 
\begin{equation}
\T_{\overline{M}}\left(-\log M_\infty\right)\simeq \O_{\overline{M}}\frac{\p}{\p s_1}
\bigoplus_{\substack{
1\le i\le 3, \\
1\le j\le a_{i}-1}} 
\O_{\overline{M}}\frac{\p}{\p s_{i,j}}
\bigoplus \O_{\overline{M}}s_{\mu_A}\frac{\p}{\p s_{\mu_A}}.
\end{equation}
\qed
\end{pf}
\begin{prop}\label{prop:KS-extension} 
The Kodaira--Spencer map $\rho$ induces the $\O_{\overline{M}}$-isomorphism $\overline{\rho}$
\begin{equation}
\overline{\rho}:\T_{\overline{M}}\left(-\log M_\infty\right)\longrightarrow {p}_*\O_{\overline{\C}}.
\end{equation}
\end{prop}
\begin{pf}
The statement easily follows since the Kodaira--Spencer map $\rho$ induces the $\O_{\overline{M}}$-isomorphism 
between $\T_{\overline{M}}\left(-\log M_\infty\right)$ and the free 
$\O_{\overline{M}}$-submodule of $\iota_*\left(p_*\O_\C\right)$ spanned by the residue classes of monomials 
\begin{equation}
1,\ x_i^j\  (i=1,2,3,\ j=1,\dots, a_i-1),\ s_{\mu_A}^{-1}x_1x_2x_3.
\end{equation}
\qed
\end{pf}
%%%%%%%%%%%%%%%%%%%%%%%%%%%%%%%%%%%%%%%%%%%%%%%%%%%%%%%%%%%%%%%%%%%%%%%%%%%%%%
\subsection{Primitive vector field and Euler vector field}
\begin{defn}\label{defn:vector fields}
The vector field $e$ and $E$ on $M$ corresponding to the unit 
$1$ and $F$ by the $\O_M$-isomorphism~\eqref{universal} 
is called the {\it primitive vector field} and the {\it Euler vector field}, respectively. 
That is, 
\begin{equation}
{e}F_A=1\ \text{and}\ {E}F_A=F_A\ 
\text{in}\ p_*\O_\C.
\end{equation}
\end{defn}
\begin{prop}
The primitive vector field $e$ and the Euler vector field $E$ on $M$ are given by 
\begin{equation}
e=\frac{\p}{\p s_1},\quad E=s_1\frac{\p }{\p s_1}+\sum_{i=1}^3\sum_{j=1}^{a_i-1}\frac{a_i-j}{a_i}s_{i,j}\frac{\p}{s_{i,j}}+
\chi_A s_{\mu_A}\frac{\p}{\p s_{\mu_A}}.
\end{equation}
In particular, they naturally define elements of $\Gamma(\overline{M},\T_{\overline{M}}\left(-\log M_\infty\right))$.
\end{prop}
\begin{pf}
This is clear.
\qed
\end{pf}
\begin{lem}
We have the ``Euler's identity"$:$
\begin{equation}
F_A=EF_A+\sum_{i=1}^3\frac{1}{a_i}x_i\frac{\p F_A}{\p x_i}.
\end{equation}
\end{lem}
\begin{pf}
This is also clear.
\qed
\end{pf}
\begin{defn}
An element $g\in\O_M$ is of {\it degree} $l$ for some $k\in \QQ$ if 
it satisfies the equation $Eg=lg$ where $E$ is the Euler vector field. 
If $Eg=lg$, then $l$ is denoted by $\deg(g)$.
\end{defn}
\begin{lem}
Suppose that $g\in\O_{\overline{M}}$ satisfies $Eg=\deg(g)\cdot g$.
Then, we have $\deg(g)\ge 0$ and $\deg(g)=0$ if and only if $g$ is a constant function.
\end{lem}
\begin{pf}
It follows from our assumption that $\chi_A>0$.
\qed
\end{pf}

%%%%%%%%%%%%%%%%%%%%%%%%%%%%%%%%%%%%%%%%%%%%%%%%%%%%%%%%%%%%%%%%%%%%%%%%%%%%%%
\subsection{Filtered de Rham cohomology}
For any non-negative integer $i$, we set
\begin{equation}
p_*\Omega^i_{\AA^3\times M/M}:=\Gamma(\AA^3,\Omega_{\AA^3}^i)\otimes_\CC\O_M,
\end{equation}
where $\Omega_{\AA^3}^i$ denotes the sheaf of {\it algebraic} differential $i$-forms on $\AA^3$, 
the algebraic variety underlying the complex manifold $\CC^3$.
$p_*\Omega^i_{\AA^3\times M/M}$ can be thought of as the direct image of 
the sheaf of relative algebraic differential $i$-forms with respect to the projection $p:\CC^3\times M\longrightarrow M$.
\begin{defn}
Define $\Omega_{F_A}$ as
\begin{equation}
\Omega_{F_A}:=p_{*}\Omega^{3}_{\AA^{3}\times M/M}/dF_{A}\wedge p_{*}\Omega^{2}_{\AA^{3}\times M/M}.
\end{equation}
\end{defn}
\begin{prop}\label{omega F}
$\Omega_{F_A}$ is a free $p_*\O_\C$-module of rank one and hence a free $\O_M$-module of rank $\mu_A$.
\end{prop}
\begin{pf}
A $p_*\O_\C$-free base of $\Omega_{F_A}$ can be chosen as $[\omega]$, $\omega:=dx_1\wedge dx_2\wedge dx_3$,
and hence an $\O_M$-free basis of $\Omega_{F_A}$ can be chosen as 
\begin{equation}
\left\{[\omega],\ [x_i^j\omega]\  (i=1,2,3,\ j=1,\dots, a_i-1),\ [x_1x_2x_3\omega]\right\}.
\end{equation}
\qed
\end{pf}
\begin{defn}
We set 
\begin{equation}
\H_{F_A}:=\RR^3 p_{*}(\Omega^{\bullet}_{\AA^{3}\times M/M}\otimes_{\O_{M}}\O_M((u)),
u d+dF_A\wedge)
\end{equation}
and call it the {\it filtered de Rham cohomology group} of the universal unfolding $F_A$.
\end{defn}
For any $k\in\ZZ$, put
\[
\H_{F_A}^{(-k)}:=\RR^3 p_{*}(\Omega^{\bullet}_{\AA^{3}\times M/M}\otimes_{\O_{M}}\O_M[[u]]u^k,
u d+dF_A\wedge).
\]
Obviously, one has an $\O_M[[u]]$-isomorphism for all $i\in\ZZ$
\begin{equation}
\H_{F_A}^{(0)}\simeq \H_{F_A}^{(-k)},\quad \omega\mapsto u^k\omega,
\end{equation}
and an $\O_M((u))$-isomorphism 
\begin{equation}
\H_{F_A}^{(-k)}\otimes_{\O_M[[u]]}\O_M((u))\simeq \H_{F_A}.
\end{equation}
Furthermore, $\H_{F_A}^{(-k)}$ is naturally a submodule of $\H_{F_A}$ 
so that $\{\H_{F_A}^{(-k)}\}_{i\in\ZZ}$ form an increasing and exhaustive filtration of $\H_{F_A}$:
\begin{equation}\label{filt}
\cdots \subset \H^{(-k-1)}_{F_A} \subset \H^{(-k)}_{F_A} \subset \cdots \subset 
\H^{(-1)}_{F_A}\subset \H^{(0)}_{F_A} \subset \cdots \subset \H_{F_A},  
\end{equation}
such that $\H_{F_A}$ is complete with respect to the filtration in the following sense:
\begin{equation}
\bigcup_{k\in \ZZ} \H^{(-k)}_{F_A}=\H_{F_A} \ \text{and} \ \bigcap_{k\in \ZZ} \H^{(-k)}_{F_A}=\{0\}.  
\end{equation}
\begin{prop}\label{b-lat-k}
For any $k\in\ZZ$, $\H_{F_A}^{(-k)}$ is an $\O_M[[u]]$-free module of rank $\mu_A$.
In particular, we have the following short exact sequence of $\O_M$-modules
\begin{equation}
0 \rightarrow \H^{(-1)}_{F_A}\hookrightarrow \H^{(0)}_{F_A} \xrightarrow{r^{(0)}} \Omega_{F_A}\rightarrow 0.
\end{equation}
\end{prop}
\begin{pf}
One can choose an $\O_M[[u]]$-free basis of $\H^{(0)}_{F_A}$ as 
\begin{equation}
\left\{[\omega],\ [x_i^j\omega]\  (i=1,2,3,\ j=1,\dots, a_i-1),\ [x_1x_2x_3\omega]\right\}
\end{equation}
where $\omega=dx_1\wedge dx_2\wedge dx_3$. The rest is clear.
\qed
\end{pf}
\begin{defn}
Define an element $\zeta_A$ of $\Gamma(M,\H_{F_A}^{(0)})$ as 
\begin{equation}
\zeta_A:=\left[s_{\mu_A}^{-1}dx_1\wedge dx_2\wedge dx_3\right].
\end{equation}
\end{defn}
\begin{defn}
For any $k\in \ZZ$, define an $\O_{\overline{M}}[[u]]$-free module $\overline{\H}_{F_A}^{(-k)}$ of rank $\mu_A$ 
as 
\begin{equation}
\overline{\H}_{F_A}^{(k)}:=
\O_{\overline{M}}[[u]]\cdot u^k\zeta_A
\bigoplus_{\substack{
1\le i\le 3, \\
1\le j\le a_{i}-1}} 
\O_{\overline{M}}[[u]]\cdot u^kx_i^j\zeta_A
\bigoplus \O_{\overline{M}}[[u]]\cdot u^ks_{\mu_A}^{-1}x_1x_2x_3\zeta_A,
\end{equation}
where we regard $\zeta_A, x_i^j\zeta_A, s_{\mu_A}^{-1}x_1x_2x_3\zeta_A$ as 
elements of $\Gamma(\overline{M},\iota_*\H_{F_A}^{(0)})$.
\end{defn}
\begin{prop}\label{extend de Rham}
We have the isomorphism of $\O_M[[u]]$-modules
\begin{equation}
\iota^*\left(\overline{\H}_{F_A}^{(-k)}\right)\simeq\H_{F_A}^{(-k)},\quad i\in\ZZ.
\end{equation}
Namely, $\overline{\H}_{F_A}^{(-k)}$ is an extension of $\H_{F_A}^{(-k)}$ as a free $\O_{\overline{M}}[[u]]$-module.
\end{prop}
\begin{pf}
It is clear by Proposition~\ref{b-lat-k}.
\qed
\end{pf}
\begin{defn}\label{defn-extend-omega F}
Define an $\O_{\overline{M}}$-free module $\overline{\Omega}_{F_A}$ of rank $\mu_A$ 
as 
\begin{equation}
\overline{\Omega}_{F_A}:=
\O_{\overline{M}}\cdot r^{(0)}(\zeta_A)
\bigoplus_{\substack{
1\le i\le 3, \\
1\le j\le a_{i}-1}} 
\O_{\overline{M}}\cdot r^{(0)}(x_i^j\zeta_A)
\bigoplus \O_{\overline{M}}\cdot r^{(0)}(s_{\mu_A}^{-1}x_1x_2x_3\zeta_A),
\end{equation}
where we regard $r^{(0)}(\zeta_A), r^{(0)}(x_i^j\zeta_A), r^{(0)}(s_{\mu_A}^{-1}x_1x_2x_3\zeta_A)$ as elements 
of $\Gamma(\overline{M},\iota_*\Omega_{F_A})$.
\end{defn}
\begin{prop}\label{extend omega F}
We have the isomorphism of $\O_M$-modules
\begin{equation}
\iota^*\left(\overline{\Omega}_{F_A}\right)\simeq\Omega_{F_A}.
\end{equation}
Namely, $\overline{\Omega}_{F_A}$ is an extension of $\Omega_{F_A}$ as a free $\O_{\overline{M}}$-module.
\end{prop}
\begin{pf}
It is almost clear.
\qed
\end{pf}
\begin{prop}\label{rank1-extend-omega F}
$\overline{\Omega}_{F_A}$ is a free ${p}_{*}\O_{\overline{\C}}$-module of rank one.
\end{prop}
\begin{pf}
A ${p}_*\O_{\overline{\C}}$-free base of $\overline{\Omega}_{F_A}$ can be chosen as $r^{(0)}(\zeta_A)$ in Definition \ref{defn-extend-omega F}.
\qed
\end{pf}
\begin{prop}\label{extend-b-lat}
We have the following short exact sequence of $\O_{\overline{M}}$-modules
\begin{equation}
0 \rightarrow \overline{\H}^{(-1)}_{F_A}\hookrightarrow \overline{\H}^{(0)}_{F_A} \xrightarrow{\overline{r}^{(0)}} \overline{\Omega}_{F_A}\rightarrow 0.
\end{equation}
\end{prop}
\begin{pf}
It is clear by their definitions.
\qed
\end{pf}
%%%%%%%%%%%%%%%%%%%%%%%%%%%%%%%%%%%%%%%%%%%%%%%%%%%%%%%%%%%%%%%%%%%%%%%%%%%%%%
\subsection{Gau\ss--Manin connection}
We define the free $\O_M[[u]]$-module $\T_{\CC_u\times M}$ of rank $\mu_A+1$ as follows:
\begin{equation}
\T_{\CC_u\times M}:=\O_M[[u]]\frac{d}{du}\oplus \O_M[[u]]\otimes_{\O_M}\T_M.
\end{equation}
\begin{defn}
We define a connection, called the {\it Gau\ss--Manin connection},
\begin{equation}
\nabla: \T_{\CC_u\times M} \otimes_{\O_M} \H_{F_A} \rightarrow \H_{F_A}
\end{equation}
by letting; for $\delta \in \T_{M}$ and $\zeta =[\phi dx_1\wedge \dots \wedge dx_n]\in \H_{F_A}$,
\begin{subequations}
\begin{equation}
\nabla_{\delta} \zeta :=[(\frac{1}{u}(\delta F_A)\phi +\delta (\phi ) )dx_1\wedge \dots \wedge dx_n],
\end{equation}
\begin{equation}
\nabla_{\frac{d}{du}} \zeta :=[(-\frac{1}{u^{2}} F_A\phi +\frac{d\phi}{du})dx_1\wedge \dots \wedge dx_n].
\end{equation}
\end{subequations}
\end{defn}
\begin{prop}\label{Gauss--Manin connection}
{\it Gau\ss--Manin connection}
$\nabla:\T_{\CC_u\times M}\otimes_{\O_M}\H_{F_A}\longrightarrow \H_{F_A}$ satisfies following;
\begin{enumerate}
\item $\nabla$ is integrable$:$
\[
\left[\nabla_{\frac{d}{du}},\nabla_{\frac{d}{du}}\right]=0,\ 
\left[\nabla_{\frac{d}{du}}, \nabla_\delta \right]=0,\ 
\left[\nabla_\delta, \nabla_{\delta'}\right]=
\nabla_{[\delta,\delta']},\ \delta,\delta'\in\T_M.
\]
\item $\nabla$ satisfies {\it Griffith transversality}$:$ that is, 
\begin{equation}
\nabla:\T_M\otimes_{\O_M}\H_{F_A}^{(-k)}\longrightarrow\H_{F_A}^{(-k+1)},\quad k\in\ZZ.
\end{equation}
\item The covariant differentiation $\nabla_{\frac{d}{du}}$ satisfies 
\begin{equation}\label{eq:GM-2}
\nabla_{u\frac{d}{du}}(\H_{F_A}^{(-k)})\subset\H_{F_A}^{(-k+1)},\quad k\in\ZZ.
\end{equation}
\end{enumerate}
\end{prop}
\begin{pf}
See Proposition 4.5 of \cite{st:1} and reference there in.
\qed
\end{pf}
We shall consider the extension the Gau\ss--Manin connection $\nabla$ on $\H_{F_A}$
to the one on $\overline{\H}_{F_A}$.
Define the free $\O_{\overline{M}}[[u]]$-module $\T_{\CC_u\times \overline{M}}(-\log M_\infty)$ of rank $\mu_A+1$ as follows:
\begin{equation}
\T_{\CC_u\times \overline{M}}(-\log M_\infty):=\O_{\overline{M}}[[u]]\frac{d}{du}\oplus 
\O_{\overline{M}}[[u]]\otimes_{\O_{\overline{M}}}\T_{\overline{M}}(-\log M_\infty).
\end{equation}
\begin{prop}\label{extension GM connection}
The Gau\ss--Manin connection $\nabla$ on $\H_{F_A}$ extends to $\overline{\H}_{F_A}$ 
with logarithmic poles along $M_\infty$, namely, $\nabla$ induces the connection
\begin{equation}
\overline{\nabla}:\T_{\CC_u\times \overline{M}}(-\log M_\infty)\otimes_{\O_{\overline{M}}}\overline{\H}_{F_A}\longrightarrow
\overline{\H}_{F_A},
\end{equation}
satisfying the following conditions$:$
\begin{enumerate}
\item $\overline{\nabla}$ is integrable$:$
\[
\left[\overline{\nabla}_{\frac{d}{du}},\overline{\nabla}_{\frac{d}{du}}\right]=0,\ 
\left[\overline{\nabla}_{\frac{d}{du}}, \overline{\nabla}_\delta \right]=0,\ 
\left[\overline{\nabla}_\delta, \overline{\nabla}_{\delta'}\right]=
\overline{\nabla}_{[\delta,\delta']},\ \delta,\delta'\in\T_{\overline{M}}(-\log M_\infty).
\]
\item $\overline{\nabla}$ satisfies {\it Griffith transversality}$:$ that is, 
\begin{equation}
\overline{\nabla}:\T_{\overline{M}}(-\log M_\infty)\otimes_{\O_{\overline{M}}}
\overline{\H}_{F_A}^{(-k)}\longrightarrow\overline{\H}_{F_A}^{(-k+1)},\quad k\in\ZZ.
\end{equation}
\item The covariant differentiation $\overline{\nabla}_{\frac{d}{du}}$ satisfies 
\begin{equation}
\overline{\nabla}_{u\frac{d}{du}}(\overline{\H}_{F_A}^{(-k)})\subset\overline{\H}_{F_A}^{(-k+1)},\quad k\in\ZZ.
\end{equation}
\end{enumerate}
\end{prop}
\begin{pf}
We shall check that 
\begin{equation}\label{extend Gauss manin mu}
s_{\mu_A}\nabla_{\frac{\p}{\p s_{\mu_A}}}\left(\overline{\H}_{F_A}\right)\subset \overline{\H}_{F_A}.
\end{equation}
\begin{equation}\label{extend Gauss manin s}
\nabla_{\frac{\p}{\p s_{1}}}\left(\overline{\H}_{F_A}\right)\subset \overline{\H}_{F_A}, \quad 
\nabla_{\frac{\p}{\p s_{i,j}}}\left(\overline{\H}_{F_A}\right)\subset \overline{\H}_{F_A}, \quad
i=1,2,3, \ j=1,\dots, a_{i}-1.
\end{equation}
\begin{equation}\label{extend Gauss manin u}
u\nabla_{\frac{d}{d u}}\left(\overline{\H}_{F_A}\right)\subset \overline{\H}_{F_A}.
\end{equation}
First we shall check the condition \eqref{extend Gauss manin mu}. One has
\begin{equation}\label{eGM mu1}
s_{\mu_A}\nabla_{\frac{\p}{\p s_{\mu_A}}}\zeta_A=\frac{1}{u}s^{-1}_{\mu_A}x_1x_2x_3\zeta_A-\zeta_A,
\end{equation}
\begin{equation}\label{eGM mu2}
s_{\mu_A}\nabla_{\frac{\p}{\p s_{\mu_A}}}x_{i}^{j}\zeta_A=\frac{1}{u}s^{-1}_{\mu_A}x_1x_2x_3\cdot x_{i}^{j}\zeta_A-x_{i}^{j}\zeta_A,
\end{equation}
\begin{equation}\label{eGM mu3}
s_{\mu_A}\nabla_{\frac{\p}{\p s_{\mu_A}}}s^{-1}_{\mu_A}x_1x_2x_3\zeta_A=\frac{1}{u}(s^{-1}_{\mu_A}x_1x_2x_3)^{2}\zeta_A-2s_{\mu_A}^{-1}x_1x_2x_3\zeta_A.
\end{equation}
The images of $s^{-1}_{\mu_A}x_1x_2x_3\cdot x_{i}^{j}\zeta_A$ and $(s^{-1}_{\mu_A}x_1x_2x_3)^{2}\zeta_A$ by $r^{(0)}$ can be extended to 
$\overline{\Omega}_{F_A}$.
By Proposition~\ref{rank1-extend-omega F} and 
Proposition~\ref{extend-b-lat}, we can show that 
$1/u\cdot s^{-1}_{\mu_A}x_1x_2x_3\zeta_A, \ 1/u\cdot (s^{-1}_{\mu_A}x_1x_2x_3)^{2}\zeta_A \in \overline{\H}_{F_A}$.    
Next we shall check that the condition \eqref{extend Gauss manin u}. 
The assertion for $s_1$ is obvious. One has
\begin{equation}
\nabla_{\frac{\p}{\p s_{i',j'}}} \phi \cdot \zeta_A=\frac{1}{u} \ x_{i'}^{j'}\cdot \phi \cdot \zeta_A,
\end{equation}
where $\phi$ is an element of the set $\left\{1, \p F_A/\p s_{i,j} \ (i=1,2,3. \ j=1,\dots, a_{i}-1), \ \p F_A/\p s_{\mu_A}\right\}$.
The image of $x_{i'}^{j'}\cdot \phi \cdot \zeta_A$ by $r^{(0)}$ can be extended to 
$\overline{\Omega}_{F_A}$. By Proposition~\ref{rank1-extend-omega F} and 
Proposition~\ref{extend-b-lat}, we can show that 
$1/u\cdot \ x_{i'}^{j'}\cdot \phi \cdot \zeta_A \in \overline{\H}_{F_A}$.
Finally we shall check the condition \eqref{extend Gauss manin u}. One has
\begin{equation}
u\nabla_{\frac{d}{d u}} \phi \cdot \zeta_A=-\frac{1}{u}F_A \cdot \phi \cdot \zeta_A,
\end{equation}
where $\phi$ is an element of the set $\left\{1, \p F_A/\p s_{i,j} \ (i=1,2,3. \ j=1,\dots, a_{i}-1), \ \p F_A/\p s_{\mu_A}\right\}$.
The image of $F_A \cdot \phi \cdot \zeta_A$ by $r^{(0)}$ can be extended to 
$\overline{\Omega}_{F_A}$ since we have 
\begin{equation}
F_A=EF_A+\sum_{i=1}^3\frac{1}{a_i}x_i\frac{\p F_A}{\p x_i}.
\end{equation}
By Proposition~\ref{rank1-extend-omega F} and 
Proposition~\ref{extend-b-lat}, we can show that $-1/u\cdot F_A \cdot \phi \cdot \zeta_A \in \overline{\H}_{F_A}$.
The conditions {\rm (i)}, {\rm (ii)} and {\rm (iii)} follow from Proposition~\ref{Gauss--Manin connection}.
\qed
\end{pf}
%%%%%%%%%%%%%%%%%%%%%%%%%%%%%%%%%%%%%%%%%%%%%%%%%%%%%%%%%%%%%%%%%%%%%%%%%%%%%%
\subsection{Higher residue pairing}

\begin{defn}
Define an $\O_M$-bilinear form $J_A$ on $\Omega_{F_A}$ by
\begin{align} 
J_{F_A}(\omega_1, \omega_2):=&
- {\rm Res}_{\CC^3\times M/M}\left[
\begin{gathered}
\phi_1\phi_2dx_1\wedge dx_2 \wedge dx_3\\
\frac{\p F_A}{\p x_1}\ \frac{\p F_A}{\p x_2}\ \frac{\p F_A}{\p x_3}
\end{gathered}
\right]\\
=&\displaystyle 
\frac{-1}{(2\pi \sqrt{-1})^3}\int_{C}\frac{\phi_1 \ \phi_2 \ dx_{1}\wedge dx_{2}\wedge
dx_{3}}{\frac{\p F_A}{\p x_1}\ \frac{\p F_A}{\p x_2}\ \frac{\p F_A}{\p x_3}}
\end{align}
where $\omega_1 = [\phi_1dx_1\wedge dx_2 \wedge dx_3]$ and $\omega_2 = [\phi_2dx_1\wedge dx_2 \wedge dx_3]$, $C$ is the distinguished 
boundary of a polydisc containing
all the critical points of $F_A$ (see (3.3.7) and (3.3.8) of \cite{sabb:1}).
\end{defn}
\begin{prop}\label{Gro}
The $\O_M$-bilinear form $J_{F_A}$ on $\Omega_{F_A}$ is non-degenerate.
\end{prop}
\begin{pf}
This is a well-known fact. 
In our case, we can show this by a direct calculation$:$ 
\begin{lem}\label{lem:calculation residue}
We have
\begin{subequations}
\begin{equation}
{\rm Res}_{\CC^3\times M/M}\left[
\begin{gathered}
1\cdot dx_1\wedge dx_2 \wedge dx_3\\
\frac{\p F_A}{\p x_1}\ \frac{\p F_A}{\p x_2}\ \frac{\p F_A}{\p x_3}
\end{gathered}
\right]=0,
\end{equation}
\begin{equation} 
{\rm Res}_{\CC^3\times M/M}\left[
\begin{gathered}
x_i^j\cdot dx_1\wedge dx_2 \wedge dx_3\\
\frac{\p F_A}{\p x_1}\ \frac{\p F_A}{\p x_2}\ \frac{\p F_A}{\p x_3}
\end{gathered}
\right]=0,\ i=1,2,3,\ j=1,\dots, a_i-1,
\end{equation}
and 
\begin{equation}
{\rm Res}_{\CC^3\times M/M}\left[
\begin{gathered}
x_1x_2x_3\cdot dx_1\wedge dx_2 \wedge dx_3\\
\frac{\p F_A}{\p x_1}\ \frac{\p F_A}{\p x_2}\ \frac{\p F_A}{\p x_3}
\end{gathered}
\right]=s_{\mu_A}^3.
\end{equation}
\end{subequations}
\end{lem}
\begin{pf}
Some elementary calculations of residues yield the statement.
\qed
\end{pf}
Now, the statement of Proposition~\ref{Gro} easily follows. 
\qed
\end{pf}
\begin{prop}
The pairing $J_{F_A}$ on $\Omega_{F_A}$ induces a non-degenerate symmetric 
$\O_{\overline{M}}$-bilinear form $\overline{J}_{F_A}$ on $\overline{\Omega}_{F_A}$ 
whose restriction to $\iota^*\left(\overline{\Omega}_{F_A}\right)\simeq\Omega_{F_A}$ 
coincides with $J_{F_A}$.
\end{prop}
\begin{pf}
Since the elements $r^{(0)}(\zeta_A), r^{(0)}(x_i^j\zeta_A), r^{(0)}(s_{\mu_A}^{-1}x_1x_2x_3\zeta_A)$ 
form an $\O_{\overline{M}}$-basis of $\overline{\Omega}_{F_A}$,
it is clear from Lemma~\ref{lem:calculation residue}.
\qed
\end{pf}
In order to define the higher residue pairing, we prepare a notation. 
For $P=\sum_{i\in\ZZ} p_iu^i\in \O_M((u))$, 
set $P^*:=\sum_{i\in\ZZ} p_i(-u)^i$ such that $(P^*)^*=P$.
\begin{defn}
An $\O_M$-bilinear form
\begin{equation}
K_{F_A}:\H_{F_A}\otimes_{\O_M}\H_{F_A}\to\O_M((u))
\end{equation}
is called the {\it higher residue pairing} if it satisfies the following properties$:$
\begin{enumerate}
\item For all $\omega_1,\omega_2\in\H_{F_A}$,  
\begin{equation}\tag{K1}\label{K1}
K_{F_A}(\omega_1,\omega_2)=(-1)^3 K_{F_A}(\omega_2,\omega_1)^*.
\end{equation}
\item For all $P\in\O_M((u))$ and $\omega_1,\omega_2\in\H_F$, 
\begin{equation}\tag{K2}\label{K2}
PK_{F_A}(\omega_1,\omega_2)=K_{F_A}(P\omega_1,\omega_2)=K_{F_A}(\omega_1,P^*\omega_2).
\end{equation}
\item For all $\omega_1, \omega_2\in\H_{F_A}^{(0)}$,  
\begin{equation}\tag{K3}\label{K3}
K_{F_A}(\omega_1,\omega_2)\in u^{3}\O_M[[u]].
\end{equation}
\item The following diagram is commutative$:$
\begin{equation*}
\begin{CD}
K_{F_A}: & \H_{F_A}^{(0)}\times\H_{F_A}^{(0)} & @>>> & u^3\O_M[[u]]\\
& @VVV & & @VV
{{\rm mod}~u^{4}\O_M[[u]]}V\\
J_{F_A}: & \Omega_{F_A}\times\Omega_{F_A} & @>>> & u^3\O_M.
\end{CD}
\end{equation*}
\item For all $\omega_1, \omega_2\in\H_F$ and $\delta\in\T_M$, 
\begin{equation}\tag{K4}\label{K4}
\delta K_{F_A}(\omega_1,\omega_2)=K_{F_A}(\nabla_\delta\omega_1,\omega_2)+
K_{F_A}(\omega_1,\nabla_\delta\omega_2).
\end{equation}
\item For all $\omega_1, \omega_2\in\H_F$,
\begin{equation}\tag{K5}\label{K5}
u\frac{d}{du}K_{F_A}(\omega_1,\omega_2)=K_{F_A}(u\nabla_\frac{d}{du}\omega_1,\omega_2)
+K_{F_A}(\omega_1,u\nabla_\frac{d}{du}\omega_2).
\end{equation}
\end{enumerate}
\end{defn}
\begin{defn}
Define $K^{(k)}_{F_A}$ for $k\in\ZZ$ by the coefficient of the expansion of $K_{F_A}$ in $u$
\begin{equation}
K_{F_A}(\omega_1,\omega_2):=
\sum_{k\in\ZZ}K^{(k)}_{F_A}(\omega_1,\omega_2)u^{k+3},
\end{equation}
and call it the $k$-th {\it higher residue pairing}.
\end{defn}
\begin{rem}
In the paper \cite{st:1}, there were errors concerning to the signs in \eqref{K1} and \eqref{K5}.
We fix them here.
\end{rem}
\begin{prop}\label{exsist-higher-residue}
There exists a unique higher residue pairing $K_{F_A}$.
\end{prop}
\begin{pf}
Since $F_A$ is tame at any point on the deformation space $M$,
namely, there are no critical points coming from infinity.
Therefore, the statement follows from the same argument in
Theorem 2.14 and Remark 2.15 in \cite{douai-sabbah:1}.
\qed
\end{pf}
We shall consider the extension $\overline{K}_{F_A}$ of $K_{F_A}$ on $\overline{\H}_{F_A}$.
\begin{prop}\label{extension higher}
The pairing $K_{F_A}$ on $\H_{F_A}^{(0)}$ induces an $\O_{\overline{M}}$-bilinear form 
\begin{equation}
\overline{K}_{F_A}:\overline{\H}_{F_A}^{(0)}\otimes_{\O_{\overline{M}}}\overline{\H}_{F_A}^{(0)}
\longrightarrow u^3\O_{\overline{M}}[[u]]
\end{equation}
whose restriction to $\iota^*\left(\overline{\H}_{F_A}^{(0)}\right)\simeq\H_{F_A}^{(0)}$ 
coincides with $K_{F_A}$.
\end{prop}
\begin{pf}
This follows from Lemma~3.4 of \cite{her:2}, where 
$\overline{M}$, $M$, $M_\infty$, $\H_{F_A}^{(0)}$, $K_{F_A}$, $\overline{\H}_{F_A}^{(0)}$ 
correspond to $X$, $Y$, $D$, $\H$, $P$, ${}_{1}\F$ in \cite{her:2}, respectively.
\qed
\end{pf}
%%%%%%%%%%%%%%%%%%%%%%%%%%%%%%%%%%%%%%%%%%%%%%%%%%%%%%%%%%%%%%%%%%%%%%%%%%%%%%
\subsection{Primitive form}
\begin{defn}\label{primitive form zeta}
An element $\zeta\in\Gamma(M,\H_{F_A}^{(0)})$ is called a {\it primitive form} for the tuple $(\H_{F_A}^{(0)},\nabla, K_{F_A})$ 
if it satisfies following five conditions;
\begin{enumerate}
\item $u\nabla_{e}\zeta=\zeta$ and $\zeta$ induces $\O_M$-isomorphism:
\begin{equation}\tag{P1}\label{P1}
\T_M[[u]]\simeq \H_{F_A}^{(0)},\quad \sum_{k=0}^\infty\delta_k u^{k}\mapsto 
\sum_{k=0}^\infty u^{k}(u\nabla_{\delta_k}\zeta).
\end{equation}
\item We have
\begin{equation}\tag{P2}\label{P2}
K_{F_A}(u\nabla_\delta\zeta,u\nabla_{\delta'}\zeta)\in \CC\cdot u^{3},\ \text{for all}\ \delta,\delta'\in\T_M.
\end{equation}
\item There exists $r\in\CC$ such that
\begin{equation}\tag{P3}\label{P3}
\nabla_{u\frac{d}{du}+E}\zeta=r\zeta.
\end{equation}
\item There exists a connection $\ns:\T_M\times \T_M\longrightarrow \T_M$ such that 
\begin{equation}\tag{P4}\label{P4}
u\nabla_\delta\nabla_{\delta'}\zeta=\nabla_{\delta\circ\delta'}\zeta
+u\nabla_{\ns_\delta\delta'}\zeta,\quad \delta,\delta'\in\T_M.
\end{equation}
\item There exists an $\O_M$-endomorphism $N:\T_M\longrightarrow \T_M$ such that 
\begin{equation}\tag{P5}\label{P5}
u\nabla_\frac{d}{du}(u\nabla_\delta\zeta)=-\nabla_{E\circ \delta}\zeta+u\nabla_{N\delta}\zeta,\quad 
\delta\in\T_M.
\end{equation}
\end{enumerate}
In particular, the constant $r$ of {\rm (P3)} is called the {\it minimal exponent}.
\end{defn}

%%%%%%%%%%%%%%%%%%%%%%%%%%%%%%%%%%%%%%%%%%%%%%%%%%%%%%%%%%%%%%%%%%%%%%%%%%%%%%
%%%%%%%%%%%%%%%%%%%%%%%%%%%%%%%%%%%%%%%%%%%%%%%%%%%%%%%%%%%%%%%%%%%%%%%%%%%%%%
\section{Main Theorem}

\begin{thm}\label{primitive form}
The element $\zeta_A:=[s_{\mu_A}^{-1}dx_1 \wedge dx_2 \wedge dx_3]\in \H_{F_A}^{(0)}$ is a primitive form 
for the tuple $(\H_{F_A}^{(0)},\nabla, K_{F_A})$ with the minimal exponent $r=1$. 
\end{thm}
Once we have a primitive form $\zeta_A$, we obtain a Frobenius structure on $M$ by the general theory 
developed by K.~Saito. 
\begin{cor}\label{cor:prim to Frob}
The primitive form $\zeta_A$ determines a Frobenius structure of rank $\mu_A$ and dimension one 
on the deformation space $M$ of the universal unfolding of $f_A$. 
More precisely, the non-degenerate symmetric bilinear form $\eta$ on $\T_M$ defined by 
\begin{equation}
\eta(\delta,\delta'):=K_{F_A}^{(0)}(u\nabla_\delta\zeta_A,u\nabla_{\delta'}\zeta)
=J_{F_A}(r^{(0)}\left(u\nabla_\delta\zeta_A),r^{(0)}(u\nabla_{\delta'}\zeta_A)\right)
,\quad \delta,\delta'\in\T_M,
\end{equation}
together with the product $\circ$ on $\T_M$, the primitive vector field $e\in\Gamma(M,\T_M)$ and 
the Euler vector field $E\in\Gamma(M,\T_M)$ define a Frobenius structure on $M$ of rank $\mu_A$ and dimension one.
\end{cor}
\begin{pf}
This is a direct consequence of Theorem~7.5 of \cite{st:1}.
\qed
\end{pf}
\begin{defn}
For simplicity, we shall denote by $M_{f_A,\zeta_A}$ the deformation space $M$ together with the Frobenius structure on $M$ 
obtained in Corollary~\ref{cor:prim to Frob} and call it {\it the Frobenius manifold constructed from the pair $(f_A,\zeta_A)$}.
\end{defn}
For a proof of Theorem~\ref{primitive form}, we first give some useful lemmas. 
\begin{lem}\label{lem:P5'}
For $\delta\in\T_M$ such that $[\delta,E]=l\delta$ for some $l\in\QQ$, we have
\begin{equation}
\nabla_{u\frac{d}{du}+E}\left(u\nabla_\delta\zeta_A\right)=(2-l)\cdot u\nabla_\delta\zeta_A,
\end{equation}
In particular, we have the following equations$:$
\begin{subequations}
\begin{equation}\label{lem:P5'-1}
\nabla_{u\frac{d}{du}+E}\left(u\nabla_e\zeta_A\right)=1\cdot u\nabla_e\zeta_A, 
\end{equation}
\begin{equation}\label{lem:P5'-2}
\nabla_{u\frac{d}{du}+E}\left(u\nabla_{\frac{\p}{\p s_{i,j}}}\zeta_A\right)
=\left(1+\frac{j}{a_i}\right)\cdot u\nabla_{\frac{\p}{\p s_{i,j}}}\zeta_A, \quad i=1,2,3,\ j=1,\dots, a_i-1.
\end{equation}
\begin{equation}\label{lem:P5'-3}
\nabla_{u\frac{d}{du}+E}\left(u\nabla_{s_{\mu_A}\frac{\p}{\p s_{\mu_A}}}\zeta_A\right)
=2\cdot u\nabla_{s_{\mu_A}\frac{\p}{\p s_{\mu_A}}}\zeta_A.
\end{equation}
\end{subequations}
\end{lem}

\begin{pf}
Note that
\begin{eqnarray*}
\lefteqn{\hspace{-50pt}(ud+dF_A)\left(\frac{x_1}{a_1}dx_2\wedge dx_3-\frac{x_2}{a_2}dx_1\wedge dx_3+\frac{x_3}{a_3}dx_1\wedge dx_2\right)}\\
&&=\left\{u\left(\sum^{3}_{i=1}\frac{1}{a_{i}}\right)+\sum_{i=1}^3\frac{1}{a_i}x_i\frac{\p F_A}{\p x_i}\right\}dx_{1}\wedge dx_{2}\wedge dx_{3}.
\end{eqnarray*}

For $\delta\in\T_M$ such that $[\delta,E]=l\delta$ for some $l\in\QQ$, one has
\begin{align*}
&\ \nabla_{u\frac{d}{du}+E}\left(u\nabla_{\delta}\zeta_A\right)\\
=&\ \nabla_{u\frac{d}{du}+E}\left[\left(\delta(F_A)s_{\mu_A}^{-1}+u\delta(s_{\mu_A}^{-1})\right)dx_{1}
\wedge dx_{2}\wedge dx_{3}\right]\\
=&\ \left[\left\{-\frac{1}{u}F_A \delta(F_A)s_{\mu_A}^{-1}-F_A \delta(s_{\mu_A}^{-1})
+u\delta(s_{\mu_A}^{-1})\right\}
dx_{1}\wedge dx_{2}\wedge dx_{3}\right]\\
&+\left[\left\{\frac{1}{u}EF_A(\delta(F_A)s_{\mu_A}^{-1}+u\delta(s_{\mu_A}^{-1}))+E(\delta(F_A)s_{\mu_A}^{-1}+
u\delta(s_{\mu_A}^{-1}))\right\}dx_{1}\wedge dx_{2}\wedge dx_{3}\right]\\
=&\ \left[\left\{\frac{1}{u}(EF_A-F_A)\delta(F_A)s_{\mu_A}^{-1}+E\delta(F_A)s_{\mu_A}^{-1}+\delta(F_A)E(s_{\mu_A}^{-1})\right\}
dx_{1}\wedge dx_{2}\wedge dx_{3}\right]\\
&+\left[\left\{
(EF_A-F_A)\delta(s_{\mu_A}^{-1})+
uE\delta(s_{\mu_A}^{-1})+u\delta(s_{\mu_A}^{-1})\right\}dx_{1}\wedge dx_{2}\wedge dx_{3}\right]\\
=&\ \left[\left\{u\cdot \frac{1}{u}(\frac{1}{a_1}+\frac{1}{a_2}+\frac{1}{a_3})
+(1-l)-\chi_A\right\}\delta(F_A)s_{\mu_A}^{-1}dx_{1}\wedge dx_{2}\wedge dx_{3}\right]\\
&+\left[\left\{u(\frac{1}{a_1}+\frac{1}{a_2}+\frac{1}{a_3})-u\chi_A-ul+u\right\}\delta(s_{\mu_A}^{-1})dx_{1}
\wedge dx_{2}\wedge dx_{3}\right]\\
=&\ (2-l)\cdot u\nabla_{\delta}\zeta_A.
\end{align*}
Therefore we have Lemma~\ref{lem:P5'}.
\qed
\end{pf}

\begin{lem}\label{lem:key}
Let $\delta_1,\delta_2,\delta_3$ be elements of $\T_M$ satisfying $[\delta_i, E]=l_i\delta_i$, $i=1,2,3$ for some $l_1,l_2,l_3\in\QQ$.
Then, we have
\begin{equation}
\deg \left(K^{(k)}_{F_A}\left(u^2\nabla_{\delta_1}\nabla_{\delta_2} \zeta_A ,u\nabla_{\delta_3}\zeta_A\right)\right)
=2-l_1-l_2-l_3-k,\quad k\in\ZZ.
\end{equation}
\end{lem}
\begin{pf}
The equations \eqref{K4} and \eqref{K5} yields
\begin{align*}
&\ \left(u\frac{d}{du}+E\right)K_{F_A}\left(u^2\nabla_{\delta_1}\nabla_{\delta_2} \zeta_A ,u\nabla_{\delta_3}\zeta_A\right)\\
=&\ K_{F_A}\left(\nabla_{u\frac{d}{du}+E}(u\nabla_{\delta_1}(u\nabla_{\delta_2} \zeta_A)) ,u\nabla_{\delta_3}\zeta_A\right)
+K_{F_A}\left(u^2\nabla_{\delta_1}\nabla_{\delta_2} \zeta_A ,\nabla_{u\frac{d}{du}+E}(u\nabla_{\delta_3}\zeta_A)\right).
\end{align*}
By the integrability of the Gau\ss--Manin connection, we have 
\[
\left[\nabla_{u\frac{d}{du}+E}, u\nabla_{\delta_1}\right]=(1-l_1)u\nabla_{\delta_1}.
\]
The statement now easily follows from the definition of $K^{(k)}_{F_A}$ and Lemma~\ref{lem:P5'}. 
\qed
\end{pf}
%%%%%%%%%%%%%%%%%%%%%%%%%%%%%%%%%%%%%%%%%%%%%%%%%%%%%%%%%%%%%%%%%%%%%%%%%%%%%%
\subsection{Condition {\rm (i)}}
Condition {\rm (i)} is satisfied since for all $\delta\in\T_M$ we have
\begin{equation}
u\nabla_{\delta}[s_{\mu_A}^{-1}dx_1 \wedge dx_2 \wedge dx_3]=
\left[\left(\delta F+us_{\mu_A}\delta(s_{\mu_A}^{-1})\right)s_{\mu_A}^{-1}dx_1 \wedge dx_2 \wedge dx_3\right].
\end{equation}
%%%%%%%%%%%%%%%%%%%%%%%%%%%%%%%%%%%%%%%%%%%%%%%%%%%%%%%%%%%%%%%%%%%%%%%%%%%%%%
\subsection{Condition {\rm (ii)}}
It is enough to show the following$:$
\begin{lem}\label{cond2}
Let $\delta_1,\delta_2$ be elements of the set of vector fields
\[
\left\{e=\frac{\p}{\p s_1}, \frac{\p}{\p s_{1,1}},\dots, \frac{\p}{\p s_{i,j}}, \dots,\frac{\p}{\p s_{3,a_3-1}},
s_{\mu_A}\frac{\p}{\p s_{\mu_A}}\right\}.
\]
Then, for all $k\ge 1$, 
$K^{(k)}_{F_A}(u\nabla_{\delta_1}\zeta_A, u\nabla_{\delta_2}\zeta_A)=0$.
\end{lem}
\begin{pf}
By Lemma~\ref{lem:key}, we have 
\[
\deg \left(K^{(k)}_{F_A}\left(u\nabla_{\delta_1}\zeta_A ,u\nabla_{\delta_2}\zeta_A\right)\right)
=1-l_1-l_2-k,\quad k\in\ZZ.
\]
Therefore, 
$
\deg \left(K^{(k)}_{F_A}\left(u\nabla_{\delta_1} \zeta_A ,u\nabla_{\delta_2}\zeta_A\right)\right)
$
is zero if $\delta_1=\delta_2=s_{\mu_A}\frac{\p}{\p s_{\mu_A}}$, $k=1$ and is negative otherwise. 
If the degree is negative, then since $K^{(k)}_{F_A}\left(u\nabla_{\delta_1} \zeta_A ,u\nabla_{\delta_2}\zeta_A\right)\in 
\O_{\overline{M}}$ by Proposition~\ref{extension higher}, it must be zero.
Since the first higher residue pairing $K^{(1)}_{F_A}$ is skew-symmetric by the equation \eqref{K1}, 
we have
\[
K_{F_A}^{(1)}\left(u\nabla_{s_{\mu_A}\frac{\p}{\p s_{\mu_A}}} \zeta ,u\nabla_{s_{\mu_A}\frac{\p}{\p s_{\mu_A}}}\zeta\right)=0.
\]
\qed
\end{pf}
%%%%%%%%%%%%%%%%%%%%%%%%%%%%%%%%%%%%%%%%%%%%%%%%%%%%%%%%%%%%%%%%%%%%%%%%%%%%%%
\subsection{Condition {\rm (iii)}}
This is already proven in Lemma~\ref{lem:P5'}.
%%%%%%%%%%%%%%%%%%%%%%%%%%%%%%%%%%%%%%%%%%%%%%%%%%%%%%%%%%%%%%%%%%%%%%%%%%%%%%
\subsection{Condition {\rm (iv)}}

We shall prepare the following useful lemma to check Condition {\rm (iv)}:
\begin{lem}\label{equivalent P3}
Under the conditions~\eqref{P1} and \eqref{P2}, the following two conditions are equivalent for $\delta_1,\delta_2\in\T_M:$
\begin{enumerate}
\item The equation~\eqref{P4} holds.
\item For all $\delta_3\in\T_M$ and $k\ge 2$, 
$K^{(k)}_{F_A}(u^2\nabla_{\delta_1}\nabla_{\delta_2}\zeta_A, u\nabla_{\delta_3}\zeta_A)=0$.
\end{enumerate}
\end{lem}
\begin{pf}
\noindent
\underline{{\rm (i)} $\Rightarrow$ {\rm (ii)}.} 
It immediately follows from Lemma~\ref{cond2}.
\underline{{\rm (ii)} $\Rightarrow$ {\rm (i)}.}
Recall that the $K^{(0)}(u\nabla_{\delta}\zeta_A,u\nabla_{\delta'}\zeta_A)=
J_{F_A}(r^{(0)}\left(u\nabla_{\delta}\zeta_A\right),r^{(0)}\left(u\nabla_{\delta'}\zeta_A\right))$ by \eqref{K3}.
Since $J_{F_A}$ is non-degenerate, the conditions \eqref{P1} and \eqref{P2} yield the statement.
\qed
\end{pf}
By Lemma~\ref{equivalent P3}, it is enough to show the following:
\begin{lem}\label{cond4}
Let $\delta_1,\delta_2,\delta_3$ be elements of the set of vector fields
\[
\left\{e=\frac{\p}{\p s_1}, \frac{\p}{\p s_{1,1}},\dots, \frac{\p}{\p s_{i,j}}, \dots,\frac{\p}{\p s_{3,a_3-1}},
s_{\mu_A}\frac{\p}{\p s_{\mu_A}}\right\}.
\]
Then, $K^{(k)}_{F_A}(u^2\nabla_{\delta_1}\nabla_{\delta_2}\zeta_A, u\nabla_{\delta_3}\zeta_A)=0$ 
for all $\delta_1,\delta_2, \delta_3\in\T_M$ and $k\ge 2$ except for one case when 
$\delta_1=\delta_2=\delta_3=s_{\mu_A}\frac{\p}{\p s_{\mu_A}}$, $k=2$. 
\end{lem}
\begin{pf}
By Lemma~\ref{lem:key}, for all $k\ge 2$, 
$\deg \left(K^{(k)}_{F_A}\left(u^2\nabla_{\delta_1}\nabla_{\delta_2} \zeta_A ,u\nabla_{\delta_3}\zeta_A\right)\right)=0$
if $\delta_1=\delta_2=\delta_3=s_{\mu_A}\frac{\p}{\p s_{\mu_A}}$, $k=2$ and is negative otherwise.  
If the degree is negative, then since $K^{(k)}_{F_A}\left(u^2\nabla_{\delta_1}\nabla_{\delta_2} \zeta_A ,u\nabla_{\delta_3}\zeta_A\right)\in \O_{\overline{M}}$ by 
Proposition~\ref{extension higher}, it must be zero. 
If $\delta_1=\delta_2=\delta_3=s_{\mu_A}\frac{\p}{\p s_{\mu_A}}$ and $k=2$, 
$K^{(2)}_{F_A}\left(u^2\nabla_{s_{\mu_A}\frac{\p}{\p s_{\mu_A}}}\nabla_{s_{\mu_A}\frac{\p}{\p s_{\mu_A}}} \zeta_A ,u\nabla_{s_{\mu_A}\frac{\p}{\p s_{\mu_A}}}\zeta_A\right)$
is a constant. Then we only have to calculate this value along the submanifold $\{{\bf s}={\bf 0}\}$. 
By Lemma~\ref{lem:P5'}, one has 
$\nabla_{u\frac{d}{du}+E}\left(u\nabla_{s_{\mu_A}\frac{\p}{\p s_{\mu_A}}}\zeta_A\right)=
2u\nabla_{s_{\mu_A}\frac{\p}{\p s_{\mu_A}}}\zeta_A$.
This implies 
\begin{align*}
\chi_A \nabla_{s_{\mu_A}\frac{\p}{\p s_{\mu_A}}} \left(u\nabla_{s_{\mu_A}\frac{\p}{\p s_{\mu_A}}}\zeta_A\right)
=&2u\nabla_{s_{\mu_A}\frac{\p}{\p s_{\mu_A}}}\zeta_A
-\nabla_{u\frac{d}{du}}\left( u\nabla_{s_{\mu_A}\frac{\p}{\p s_{\mu_A}}}\zeta_A\right)\\
&-\nabla_{E-\chi_As_{\mu_A}\frac{\p}{\p s_{\mu_A}}}\left(u\nabla_{s_{\mu_A}\frac{\p}{\p s_{\mu_A}}}\zeta_A\right).
\end{align*}
Then, by $\chi_A\ne 0$ and Lemma~\ref{cond2}, one has
\begin{align*}
&\left.K^{(2)}_{F_A}\left(u^2\nabla_{s_{\mu_A}\frac{\p}{\p s_{\mu_A}}}\nabla_{s_{\mu_A}\frac{\p}{\p s_{\mu_A}}} \zeta_A ,u\nabla_{s_{\mu_A}\frac{\p}{\p s_{\mu_A}}}\zeta_A\right)\right|_{\bf s=0}\\
=&\left.\frac{1}{\chi_A}K^{(2)}_{F_A}\left(2u^2\nabla_{s_{\mu_A}\frac{\p}{\p s_{\mu_A}}}\zeta_A,
u\nabla_{s_{\mu_A}\frac{\p}{\p s_{\mu_A}}}\zeta_A\right)\right|_{\bf s=0}\\
&-\left.\frac{1}{\chi_A}K^{(2)}_{F_A}\left(u^2\nabla_{u\frac{d}{du}}
\left( u\nabla_{s_{\mu_A}\frac{\p}{\p s_{\mu_A}}}\zeta_A\right),u\nabla_{s_{\mu_A}\frac{\p}{\p s_{\mu_A}}}\zeta_A\right)\right|_{\bf s=0}\\
=&-\left.\frac{1}{\chi_A}K^{(0)}_{F_A}\left(\nabla_{u\frac{d}{du}} 
\left(u\nabla_{s_{\mu_A}\frac{\p}{\p s_{\mu_A}}}\zeta_A\right),u\nabla_{s_{\mu_A}\frac{\p}{\p s_{\mu_A}}}\zeta_A\right)\right|_{\bf s=0}.
\end{align*}
By the symmetry of $K^{(0)}_{F_A}$ and the equation \eqref{K5}, one also has
\begin{align*}
&K^{(0)}_{F_A}\left(\nabla_{u\frac{d}{du}} \left(u\nabla_{s_{\mu_A}\frac{\p}{\p s_{\mu_A}}}\zeta_A\right),u\nabla_{s_{\mu_A}\frac{\p}{\p s_{\mu_A}}}\zeta_A\right)\\
=&\frac{1}{2}\left\{K^{(0)}_{F_A}\left(\nabla_{u\frac{d}{du}} \left(u\nabla_{s_{\mu_A}\frac{\p}{\p s_{\mu_A}}}\zeta_A\right),
u\nabla_{s_{\mu_A}\frac{\p}{\p s_{\mu_A}}}\zeta_A\right)\right.\\
&\left.+K^{(0)}_{F_A}\left(u\nabla_{s_{\mu_A}\frac{\p}{\p s_{\mu_A}}}\zeta_A,
\nabla_{u\frac{d}{du}} \left(u\nabla_{s_{\mu_A}\frac{\p}{\p s_{\mu_A}}}\zeta_A\right)\right)\right\}\\
=&\frac{3}{2}K^{(0)}_{F_A}\left(u\nabla_{s_{\mu_A}\frac{\p}{\p s_{\mu_A}}}\zeta_A, u\nabla_{s_{\mu_A}\frac{\p}{\p s_{\mu_A}}}\zeta_A\right)=0.
\end{align*}
Therefore we have Lemma~\ref{cond4}.
\qed
\end{pf}
%%%%%%%%%%%%%%%%%%%%%%%%%%%%%%%%%%%%%%%%%%%%%%%%%%%%%%%%%%%%%%%%%%%%%%%%%%%%%%
\subsection{Condition {\rm (v)}}
Under the Condition {\rm (iv)}, it is easy to see that Condition {\rm (v)} follows from Lemma~\ref{lem:P5'}.

%%%%%%%%%%%%%%%%%%%%%%%%%%%%%%%%%%%%%%%%%%%%%%%%%%%%%%%%%%%%%%%%%%%%%%%%%%%%%%
%%%%%%%%%%%%%%%%%%%%%%%%%%%%%%%%%%%%%%%%%%%%%%%%%%%%%%%%%%%%%%%%%%%%%%%%%%%%%%
\section{Mirror Symmetry}
From now on, we shall denote by $\F_{f_A,\zeta_A}$ the Frobenius potential for the Frobenius structure 
on $M$ constructed from the pair $(f_A,\zeta_A)$ in Corollary~\ref{cor:prim to Frob}. 
\begin{thm}\label{thm:satisfies IST1}
For the Frobenius manifold $M_{f_A,\zeta_A}$ constructed from the pair $(f_A,\zeta_A)$, 
there exists flat coordinates $t_1$, $t_{1,1},\dots ,t_{3,a_3-1}$, $t_{\mu_A}$ satisfying the following conditions$:$
\begin{enumerate}
\item 
The unit vector field $e$ and the Euler vector field $E$ are given by
\[
e=\frac{\p}{\p t_1},\ E=t_1\frac{\p}{\p t_1}+\sum_{i=1}^3\sum_{j=1}^{a_i-1}\frac{a_i-j}{a_i}t_{i,j}\frac{\p}{\p t_{i,j}}
+\chi_A\frac{\p}{\p t_{\mu_A}}.
\]
\item 
The non-degenerate symmetric bilinear form $\eta$ on $\T_M$ defined by 
\begin{equation}
\eta(\delta,\delta'):=K_{F_A}^{(0)}(u\nabla_\delta\zeta_A,u\nabla_{\delta'}\zeta),\quad \delta,\delta'\in\T_M,
\end{equation}
satisfies
\begin{align*}
&\ \eta\left(\frac{\p}{\p t_1}, \frac{\p}{\p t_{\mu_A}}\right)=
\eta\left(\frac{\p}{\p t_{\mu_A}}, \frac{\p}{\p t_1}\right)=1,\\ 
&\ \eta\left(\frac{\p}{\p t_{i_1,j_1}}, \frac{\p}{\p t_{i_2,j_2}}\right)=
\begin{cases}
\frac{1}{a_{i_1}}\quad i_1=i_2\text{ and }j_2=a_{i_1}-j_1,\\
0 \quad \text{otherwise}.
\end{cases}
\end{align*}
\item 
The Frobenius potential $\F_{f_A,\zeta_A}$ satisfies $E\F_{f_A,\zeta_A}|_{t_{1}=0}=2\F_{f_A,\zeta_A}|_{t_{1}=0}$,
\[
\left.\F_{f_A,\zeta_A}\right|_{t_1=0}\in\CC\left[t_{1,1}, \dots, t_{1,a_1-1},
t_{2,1}, \dots, t_{2,a_2-1},t_{3,1}, \dots, t_{3,a_3-1},e^{t_{\mu_A}}\right].
\]
\item Assume the condition {\rm (iii)}. The restriction of the Frobenius potential $\F_{f_A,\zeta_A}$ to 
the submanifold $\{t_1=e^{t_{\mu_A}}=0\}$ is given as 
\[
\F_{f_A,\zeta_A}|_{t_1=e^{t_{\mu_A}}=0}=\G^{(1)}+\G^{(2)}+\G^{(3)}, 
\]
where $\G^{(i)}\in \CC[t_{i,1},\dots, t_{i,a_i-1}]$, $i=1,2,3$.
\item 
Assume the condition {\rm (iii)}. In the frame $\frac{\p}{\p t_1}$, $\frac{\p}{\p t_{1,1}},\dots, 
\frac{\p}{\p t_{3,a_3-1}}$, $\frac{\p}{\p t_{\mu_A}}$ of $\T_M$,
the product $\circ$ can be extended to the limit $t_1=t_{1,1}=\dots=t_{3,a_3-1}=e^{t_{\mu_A}}=0$.
The $\CC$-algebra obtained in this limit is isomorphic to 
\[
\CC[x_1,x_2,x_3]\left/\left(x_1x_2,x_2x_3,x_3x_1,a_1x_1^{a_1}-a_2x_2^{a_2},
a_2x_2^{a_2}-a_3x_3^{a_3}\right)\right.,
\]
where $\p/\p t_{1,1},\p/\p t_{2,1},\p/\p t_{3,1}$ are mapped to
$x_1,x_2,x_3$, respectively.
\item The term 
\[
\begin{cases}
e^{t_{\mu_A}}\quad \textit{if}\quad a_{1}=a_{2}=a_{3}=1,\\
t_{3,1}e^{t_{\mu_A}}\quad \textit{if}\quad 1=a_{1}=a_{2}<a_{3},\\
t_{2,1}t_{3,1}e^{t_{\mu_A}}\quad \textit{if}\quad 1=a_{1}<a_{2},\\
t_{1,1}t_{2,1}t_{3,1}e^{t_{\mu_A}}\quad \textit{if}\quad a_{1}=2,
\end{cases}
\]
occurs with the coeffcient $1$ in $\F_{f_A,\zeta_A}$. 
\end{enumerate}
\end{thm}
We shall check all conditions in the above theorem.
%%%%%%%%%%%%%%%%%%%%%%%%%%%%%%%%%%%%%%%%%%%%%%%%%%%%%%%%%%%%%%%%%%%%%%%%%%%%%%
\subsection{Condition {\rm (i)}}
\begin{lem}\label{lem:0-th residue}
We have
\begin{equation}
K_{F_A}^{(0)}(\zeta_A,u\nabla_{\frac{\p}{\p s_1}}\zeta_A)=0,\ 
K_{F_A}^{(0)}(\zeta_A,u\nabla_{\frac{\p}{\p s_{i,j}}}\zeta_A)=0,\ i=1,2,3,j=1,\dots, a_i-1,
\end{equation}
and 
\begin{equation}
K_{F_A}^{(0)}(\zeta_A,u\nabla_{s_{\mu_A}\frac{\p}{\p s_{\mu_A}}}\zeta_A)=1.
\end{equation}
\end{lem}
\begin{pf}
The statement follows from Lemma~\ref{lem:calculation residue}.
\qed
\end{pf}
\begin{lem}\label{lem:flat fcn}
The one form $\theta\in\Gamma(M,\Omega_M^1)$ defined by 
\begin{multline}
\theta: =K_{F_A}^{(0)}(\zeta_A,u\nabla_{\frac{\p}{\p s_1}}\zeta_A)ds_1\\
+\sum_{i=1}^3\sum_{j=1}^{a_i-1}K_{F_A}^{(0)}(\zeta_A,u\nabla_{\frac{\p}{\p s_{i,j}}}\zeta_A)ds_{i,j}
+K_{F_A}^{(0)}(\zeta_A,u\nabla_{\frac{\p}{\p s_{\mu_A}}}\zeta_A)ds_{\mu_A},
\end{multline}
is a closed form which is independent from the choice of coordinates on $M$.
Moreover, there exists a flat coordinate $t$ such that $\theta=dt$.
\end{lem}
\begin{pf}
See section 3.3 3) of \cite{S1202-Saito}.
\qed
\end{pf}
Combining these two Lemmas, we have $dt=ds_{\mu_A}/s_{\mu_A}$. 
Therefore, we can choose $t_{\mu_A}:=\log s_{\mu_A}$ as a flat coordinate. 
Since $\zeta_A, x_i^j\zeta_A (i=1,2,3,j=1,\dots, a_i-1), s_{\mu_A}^{-1}x_1x_2x_3\zeta_A$ form 
a $\CC[[u]]$-basis of $\overline{\H}_{F_A}^{(0)}$ at $({\bf 0},0)\in\overline{M}$, by Lemma~\ref{lem:P5'},
one can choose other flat coordinates $t_1, t_{1,1}, \dots, t_{3,a_3-1}$ such that
\[
\left.t_1\right|_{({\bf s},s_{\mu_A})=({\bf 0},0)}=\left.t_{1,1}\right|_{({\bf s},s_{\mu_A})=({\bf 0},0)}=\dots=
\left.t_{3,a_3-1}\right|_{({\bf s},s_{\mu_A})=({\bf 0},0)}=0
\]
together with the following normalization$;$
\begin{subequations}\label{eq:s-t}
\begin{align}
&\left.\frac{\p t_1}{\p s_1}\right|_{({\bf s},s_{\mu_A})=({\bf 0},0)}=1,\ 
\left.\frac{\p t_1}{\p s_{i,j}}\right|_{({\bf s},s_{\mu_A})=({\bf 0},0)}=0,\ 
\left.\frac{\p t_1}{\p s_{\mu_A}}\right|_{({\bf s},s_{\mu_A})=({\bf 0},0)}=0,\\ 
&\left.\frac{\p t_{i,j}}{\p s_1}\right|_{({\bf s},s_{\mu_A})=({\bf 0},0)}=0,\ 
\left.\frac{\p t_{i,j}}{\p s_{i',j'}}\right|_{({\bf s},s_{\mu_A})=({\bf 0},0)}=
\delta_{ii'}\delta_{jj'},\ \left.\frac{\p t_{i,j}}{\p s_{\mu_A}}\right|_{({\bf s},s_{\mu_A})=({\bf 0},0)}=0,
\end{align}
\end{subequations}
where $\delta_{ii'}$ and $\delta_{jj'}$ are Kronecker's deltas.
In particular, the flat coordinates $t_1$, $t_{1,1}$, $\dots$, $t_{3,a_3-1}$, $t_{\mu_A}$ satisfy 
\[
e=\frac{\p}{\p t_1},\ E=t_1\frac{\p}{\p t_1}+\sum_{i=1}^3\sum_{j=1}^{a_i-1}\frac{a_i-j}{a_i}t_{i,j}\frac{\p}{\p t_{i,j}}
+\chi_A\frac{\p}{\p t_{\mu_A}},
\]
which is the condition {\rm (i)}.
%%%%%%%%%%%%%%%%%%%%%%%%%%%%%%%%%%%%%%%%%%%%%%%%%%%%%%%%%%%%%%%%%%%%%%%%%%%%%%
\subsection{Condition {\rm (ii)}}
\begin{lem}
We have
\begin{equation}
K_{F_A}^{(0)}(\zeta_A,u\nabla_{\frac{\p}{\p t_1}}\zeta_A)=0,\ 
K_{F_A}^{(0)}(\zeta_A,u\nabla_{\frac{\p}{\p t_{i,j}}}\zeta_A)=0,\ i=1,2,3,j=1,\dots, a_i-1,
\end{equation}
and 
\begin{equation}
K_{F_A}^{(0)}(\zeta_A,u\nabla_{\frac{\p}{\p t_{\mu_A}}}\zeta_A)=1.
\end{equation}
\end{lem}
\begin{pf}
By Lemma~\ref{lem:flat fcn}, we have
\[
dt_{\mu_A} =K_{F_A}^{(0)}(\zeta_A,u\nabla_{\frac{\p}{\p t_1}}\zeta_A)dt_1
+\sum_{i=1}^3\sum_{j=1}^{a_i-1}K_{F_A}^{(0)}(\zeta_A,u\nabla_{\frac{\p}{\p t_{i,j}}}\zeta_A)dt_{i,j}
+K_{F_A}^{(0)}(\zeta_A,u\nabla_{\frac{\p}{\p t_{\mu_A}}}\zeta_A)dt_{\mu_A}.
\]
The statement follows.
\qed
\end{pf}
Note that the pairings to consider are constant since we are dealing with flat coordinates.
Therefore, we can evaluate them at $({\bf 0},0)\in\overline{M}$.
Moreover, by the normalization~\eqref{eq:s-t}, we have
\begin{align*}
&\ \lim_{s_{\mu_A}\to 0}\left.K_{F_A}^{(0)}\left(u\nabla_{\frac{\p}{\p t_{i,j}}}\zeta_A,u\nabla_{\frac{\p}{\p t_{i',j'}}}\zeta_A\right)
\right|_{{\bf s}=0}\\
=&\ \left.\overline{K}_{F_A}^{(0)}\left(u\overline{\nabla}_{\frac{\p}{\p t_{i,j}}}\zeta_A,u\overline{\nabla}_{\frac{\p}{\p t_{i',j'}}}\zeta_A\right)\right|_{({\bf s},s_{\mu_A})=({\bf 0},0)}\\
=&\ \left.\overline{K}_{F_A}^{(0)}\left(u\overline{\nabla}_{\frac{\p}{\p s_{i,j}}}\zeta_A,u\overline{\nabla}_{\frac{\p}{\p s_{i',j'}}}\zeta_A\right)
\right|_{({\bf s},s_{\mu_A})=({\bf 0},0)}\\
=&\ \lim_{s_{\mu_A}\to 0}\left.K_{F_A}^{(0)}\left(u\nabla_{\frac{\p}{\p s_{i,j}}}\zeta_A,u\nabla_{\frac{\p}{\p s_{i',j'}}}\zeta_A\right)
\right|_{{\bf s}=0}\\
=&\ \lim_{s_{\mu_A}\to 0}\left.K_{F_A}^{(0)}\left(\zeta_A,u\nabla_{\frac{\p}{\p s_{i,j}}\circ \frac{\p}{\p s_{i',j'}}}\zeta_A\right)
\right|_{{\bf s}=0}.
\end{align*}
The statement follows from Lemma~\ref{lem:limit of Jac}.
%%%%%%%%%%%%%%%%%%%%%%%%%%%%%%%%%%%%%%%%%%%%%%%%%%%%%%%%%%%%%%%%%%%%%%%%%%%%%%
\subsection{Condition {\rm (iii)}}
Note that we have an $\O_{\overline{M}}$-isomorphism
\begin{equation}
\T_{\overline{M}}\left(-\log M_\infty\right)\simeq \O_{\overline{M}}\frac{\p}{\p t_1}
\bigoplus_{\substack{
1\le i\le 3, \\
1\le j\le a_{i}-1}} 
\O_{\overline{M}}\frac{\p}{\p t_{i,j}}
\bigoplus \O_{\overline{M}}\frac{\p}{\p t_{\mu_A}}.
\end{equation}
The condition {\rm (iii)} follows from this together with Proposition~\ref{prop:KS-extension} 
and the fact that all flat coordinates have positive degrees.
%%%%%%%%%%%%%%%%%%%%%%%%%%%%%%%%%%%%%%%%%%%%%%%%%%%%%%%%%%%%%%%%%%%%%%%%%%%%%%
\subsection{Condition {\rm (iv)}}
Recall that the ideal in the equation~\eqref{eq:Jacobi-ring extension} restricted to $M_\infty$
is given by
\[
\left(x_2x_3,x_3x_1,x_1x_2, H_1({\bf x},{\bf s}),H_2({\bf x},{\bf s})\right), 
\]
where
\[
H_{i}({\bf x}, {\bf s}):=a_{i}x_{i}^{a_i}-a_{i+1}x_{i+1}^{a_{i+1}}+
\sum_{j=1}^{a_i-1}j\cdot s_{i,j}\cdot x_i^j-\sum_{j=1}^{a_{i+1}-1}j\cdot s_{i+1,j}\cdot x_{i+1}^j,\ i=1,2.
\]
In particular, we have 
\begin{equation}\label{eq:cond iv}
\left.\eta\left(\frac{\p}{\p s_{i,j}},\frac{\p}{\p s_{i',j'}}\right)\right|_{s_{\mu_A}=0}=0,\quad \text{if}\ i\ne i'.
\end{equation}
Note that the connection $\ns$ on $\T_M$ in \eqref{P4} is the Levi--Civita connection with respect to $\eta$
(see Proposition 7.9 and Proposition 7.16 of \cite{st:1}) and that $t_{i,j}$ is the solution of the system 
\begin{align*}
&\ns^*dt_{i,j}=0,\ \left.t_{i,j}\right|_{({\bf s},s_{\mu_A})=({\bf 0},0)}=0,\\ 
&\left.\frac{\p t_{i,j}}{\p s_1}\right|_{({\bf s},s_{\mu_A})=({\bf 0},0)}=0,\ 
\left.\frac{\p t_{i,j}}{\p s_{i',j'}}\right|_{({\bf s},s_{\mu_A})=({\bf 0},0)}=\delta_{ii'}\delta_{jj'},\ 
\left.\frac{\p t_{i,j}}{\p s_{\mu_A}}\right|_{({\bf s},s_{\mu_A})=({\bf 0},0)}=0,
\end{align*}
where $\ns^*$ is a connection on $\Omega_M^1$ dual to $\ns$.
Therefore, by \eqref{eq:cond iv}, we have
\begin{equation}\label{eq:normalization 2}
\left.\frac{\p t_{i,j}}{\p s_{i',j'}}\right|_{s_{\mu_A}=0}=0,\quad\text{if}\  i\ne i'.
\end{equation}
The third derivatives of the Frobenius potential with respect to flat coordinates 
are given by residues. For example, we have
\begin{equation}\label{3-point}
\frac{\p^3 \F_{f_A,\zeta_A}}{\p t_{i_1,j_1}\p t_{i_2,j_2}\p t_{i_3,j_3}}=-e^{-2 t_{\mu_A}}{\rm Res}_{\CC^3\times M/M}\left[
\begin{gathered}
\frac{\p F_A}{\p t_{i_1,j_1}}\frac{\p F_A}{\p t_{i_2,j_2}}\frac{\p F_A}{\p t_{i_3,j_3}} dx_1\wedge dx_2 \wedge dx_3\\
\frac{\p F_A}{\p x_1}\ \frac{\p F_A}{\p x_2}\ \frac{\p F_A}{\p x_3}
\end{gathered}
\right].
\end{equation}
Therefore, by using the above description of the ideal, the normalization~\eqref{eq:s-t} and \eqref{eq:normalization 2},
we can show that  
\[
\left.\lim_{e^{t_{\mu_A}}\to 0}\frac{\p^n \F_{f_A,\zeta_A}}{\p t_{i_1,j_1}\cdots \p t_{i_n,j_n}}
\right|_{t_1=t_{1,1}=\dots =t_{3,a_3-1}=0} \ne 0\ 
\text{only if}\ i_1=\dots=i_n,
\]
by induction on $n$.
%%%%%%%%%%%%%%%%%%%%%%%%%%%%%%%%%%%%%%%%%%%%%%%%%%%%%%%%%%%%%%%%%%%%%%%%%%%%%%
\subsection{Condition {\rm (v)}}
Taking into account the normalization~\eqref{eq:s-t},
the condition {\rm (v)} easily follows from Lemma~\ref{lem:limit of Jac}. 
%%%%%%%%%%%%%%%%%%%%%%%%%%%%%%%%%%%%%%%%%%%%%%%%%%%%%%%%%%%%%%%%%%%%%%%%%%%%%%
\subsection{Condition {\rm (vi)}}
Note that the coefficient of the term $e^{t_{\mu_A}}\ (\text{if}\ a_{1}=a_{2}=a_{3}=1)$, 
$t_{3,1}e^{t_{\mu_A}}\ (\text{if}\ 1=a_{1}=a_{2}<a_{3})$, $t_{2,1}t_{3,1}e^{t_{\mu_A}}\ (\text{if}\ 1=a_{1}<a_{2})$, 
$t_{1,1}t_{2,1}t_{3,1}e^{t_{\mu_A}}\ (\text{if}\ a_{1}\ge 2)$ is given by the limit  
\[
\lim_{e^{t_{\mu_A}} \to 0}\left(e^{-t_{\mu_A}}\cdot
\left.\frac{\p^3 \F_{f_A,\zeta_A}}{\p t_{\mu_A}\p t_{\mu_A}\p t_{\mu_A}}\right|_{t_1=0} \right) \ 
\text{if}\ a_{1}=a_{2}=a_{3}=1, 
\]
\[
\lim_{e^{t_{\mu_A}} \to 0}\left(e^{-t_{\mu_A}}\cdot
\left.\frac{\p^3 \F_{f_A,\zeta_A}}{\p t_{\mu_A}\p t_{\mu_A}\p t_{3,1}}\right|_{t_1=t_{3,1}=\dots =t_{3,a_3-1}=0} \right) \ 
\text{if}\ 1=a_{1}=a_{2}<a_{3},
\]
\[
\lim_{e^{t_{\mu_A}} \to 0}\left(e^{-t_{\mu_A}}\cdot
\left.\frac{\p^3 \F_{f_A,\zeta_A}}{\p t_{\mu_A}\p t_{2,1}\p t_{3,1}}\right|_{t_1=t_{2,1}=\dots =t_{3,a_3-1}=0} \right) \ 
\text{if}\ 1=a_{1}<a_{2},
\]
\[
\lim_{e^{t_{\mu_A}} \to 0}\left(e^{-t_{\mu_A}}\cdot
\left.\frac{\p^3 \F_{f_A,\zeta_A}}{\p t_{1,1}\p t_{2,1}\p t_{3,1}}\right|_{t_1=t_{1,1}=\dots =t_{3,a_3-1}=0} \right) \ 
\text{if}\ a_{1}\ge 2.
\]
Since the third derivatives of the Frobenius potential with respect to flat coordinates 
are given by residues of the form \eqref{3-point}, it is easy to see that the above limits are always reduced to the following limit 
\[
\lim_{s_{\mu_A} \to 0}\left(\left.s_{\mu_A}^{-1}\cdot s_{\mu_A}^{-2} {\rm Res}_{\CC^3\times M/M}\left[
\begin{gathered}
x_1x_2x_3 dx_1\wedge dx_2 \wedge dx_3\\
\frac{\p F_A}{\p x_1}\ \frac{\p F_A}{\p x_2}\ \frac{\p F_A}{\p x_3}
\end{gathered}
\right]\right|_{{\bf s}=0}\right)
\]
by using the relations in $p_*\O_\C$ and the normalization~\eqref{eq:s-t}, which is one.
%%%%%%%%%%%%%%%%%%%%%%%%%%%%%%%%%%%%%%%%%%%%%%%%%%%%%%%%%%%%%%%%%%%%%%%%%%%%%%
\subsection{Mirror isomorphism}
This theorem enables us to simplify the proofs given by Milanov--Tseng \cite{mt:1} and Rossi \cite{r:1}.
\begin{cor}\label{cor:classical mirror}
There exists an isomorphism of Frobenius manifolds between the one constructed from the Gromov--Witten theory 
for $\PP^{1}_A$ and the one constructed from the pair $(f_A,\zeta_A)$. 
\end{cor}
%%%%%%%%%%%%%%%%%%%%%%%%%%%%%%%%%%%%%%%%%%%%%%%%%%%%%%%%%%%%%%%%%%%%%%%%%%%%%%
%%%%%%%%%%%%%%%%%%%%%%%%%%%%%%%%%%%%%%%%%%%%%%%%%%%%%%%%%%%%%%%%%%%%%%%%%%%%%%
\section{Another proof of our main theorem}

Denote by $(\eta_A,\circ_A,e_A, E_A)$ the Frobenius structure for a triplet
of positive integers $A=(a_1,a_2,a_3)$ with $\chi_A>0$ uniquely determined by Theorem~3.1 of \cite{ist:1}.
In particular, we shall denote by 
$({\bf t}, t_{\mu_A}):=(t_1,t_2, \dots, t_{\mu_A-1},t_{\mu_A}):=(t_1,t_{1,1},\dots ,t_{3,a_3-1}, t_{\mu_A})$
the flat coordinates on $M$
and also denote by $\p_{i}$ the basis of $\T_{M}$ corresponding to the coordinate $t_{i}$
for simplicity.
In this section, we shall show the way to determine not only the Frobenius structure $(\eta, \circ, e,E)$ on $M$ 
constructed from the pair $(f_A,\zeta_A)$ but also
the primitive form $\zeta_A$ itself from the Frobenius structure $(\eta_A,\circ_A,e_A, E_A)$.
It is very important to note that the Frobenius potential $\F_A$ for $(\eta_A,\circ_A,e_A, E_A)$ 
is isomorphic to the one on $M$ constructed from the Gromov-Witten theory 
of the orbifold $\PP^{1}_A$.
%%%%%%%%%%%%%%%%%%%%%%%%%%%%%%%%%%%%%%%%%%%%%%%%%%%%%%%%%%%%%%%%%%%%%%%%%%%%%%
\subsection{Key lemma}
We can determine a primitive form by the following Lemma, which is essentially Lemma 3.3 of \cite{S1202-Saito}.
\begin{lem}\label{Key lemma}
Suppose that there exists a local coordinates $(t_1,t_{1,1},\dots ,t_{3,a_3-1}, t_{\mu_A})$ on $M$ 
and an element $\zeta\in\Gamma(M,\H_{F_A}^{(0)})$ satisfying the following properties$:$ 
\begin{enumerate}
\item 
The unit vector field $e$ and the Euler vector field $E$ are given by
\[
e=\frac{\p}{\p t_1},\ E=t_1\frac{\p}{\p t_1}+\sum_{i=1}^3\sum_{j=1}^{a_i-1}\frac{a_i-j}{a_i}t_{i,j}\frac{\p}{\p t_{i,j}}
+\chi_A\frac{\p}{\p t_{\mu_A}}.
\]
\item
$\zeta$ satisfies
\[
u\nabla_e\zeta=\zeta,\quad 
\nabla_{u\frac{d}{du}+E}\zeta=\zeta.
\]
\item
We have an $\O_M[[u]]$-isomorphism
\[
\H_{F_A}^{(0)}\simeq \O_M[[u]] \cdot u\nabla_\frac{\p}{\p t_1}\zeta
\bigoplus_{\substack{
1\le i\le 3, \\
1\le j\le a_{i}-1}} 
\O_M[[u]]\cdot u\nabla_\frac{\p}{\p t_{i,j}}\zeta
\bigoplus \O_M[[u]]\cdot u\nabla_\frac{\p}{\p t_{\mu_A}}\zeta.
\]
\item
Let $V$ be the $\CC$-vector space given by
\[
V:=\CC\frac{\p}{\p t_1}
\bigoplus_{\substack{
1\le i\le 3, \\
1\le j\le a_{i}-1}} 
\CC\frac{\p}{\p t_{i,j}}
\bigoplus \CC\frac{\p}{\p t_{\mu_A}}.
\]
We have 
\[
K_{F_A}^{(0)}\left(u\nabla_{\delta_1}\zeta,u\nabla_{\delta_2}\zeta\right)\in\CC,\quad\text{for all}\quad \delta_1,\delta_2\in V.
\]
\item 
We have
\[
u\nabla_{\delta_1}\nabla_{\delta_2}\zeta=\nabla_{\delta_1\circ\delta_2}\zeta,\quad\text{for all}\quad \delta_1,\delta_2\in V.
\]
\end{enumerate}
Then, the $\O_M$-bilinear form 
\begin{equation}
K'_{F_A}:\H_{F_A}\otimes_{\O_M}\H_{F_A}\longrightarrow \O_M((u))
\end{equation}
defined by 
\begin{equation}
K'_{F_A}(u^{i_1+1}\nabla_{\delta_1}\zeta,u^{i_2+1}\nabla_{\delta_2}\zeta)
:=(-1)^{i_2}\cdot K_{F_A}^{(0)}
\left(u\nabla_{\delta_1}\zeta,u\nabla_{\delta_2}\zeta\right)
\cdot u^{i_1+i_2+3},\ 
\delta_1,\delta_2\in \T_M,
\end{equation}
coincides with $K_{F_A}$. 
In particular, $\zeta$ is a primitive form for the tuple $(\H_{F_A}^{(0)},\nabla, K_{F_A})$ with the minimal exponent $r=1$. 
\end{lem}
\begin{pf}
The pairing $K'_{F_A}$ obviously satisfies the same conditions \eqref{K1}, \eqref{K2}, \eqref{K3} 
as the higher residue pairing $K_{F_A}$ does. 
We can easily show by a direct calculation that the conditions {\rm (iv)} and {\rm (v)} imply 
the condition \eqref{K4} for $K'_{F_A}$.
We have  
\[
\nabla_{u\frac{d}{du}+E}\left(u\nabla_{\frac{\p}{\p t_{i,j}}}\zeta\right)
=\left(\frac{j}{a_i}+1\right)\cdot u\nabla_{\frac{\p}{\p t_{i,j}}}\zeta, \quad i=1,2,3,\ j=1,\dots, a_i-1,
\]
\[
\nabla_{u\frac{d}{du}+E}\left(u\nabla_{\frac{\p}{\p t_{\mu_A}}}\zeta\right)=2\cdot u\nabla_{\frac{\p}{\p t_{\mu_A}}}\zeta,
\]
which implies the condition \eqref{P5} for $\zeta$. 
By using this, we can also show the condition \eqref{K5} for $K'_{F_A}$.
Therefore, by the uniquness of the higher residue pairings in Proposition~\ref{exsist-higher-residue}, we have
$K'_{F_A}=K_{F_A}$. 
$\zeta$ obviously satisfies the conditions \eqref{P1}, \eqref{P3} and \eqref{P4} by the assumption.
Since 
\[
K_{F_A}(u\nabla_{\delta_1}\zeta,u\nabla_{\delta_2}\zeta)=
K'_{F_A}(u\nabla_{\delta_1}\zeta,u\nabla_{\delta_2}\zeta)
:=K_{F_A}^{(0)}
\left(u\nabla_{\delta_1}\zeta,u\nabla_{\delta_2}\zeta\right)
\cdot u^{3},\ 
\delta_1,\delta_2\in \T_M,
\]
$\zeta$ also satisfies the condition \eqref{P2}.
The statement follows.
\qed
\end{pf}
%%%%%%%%%%%%%%%%%%%%%%%%%%%%%%%%%%%%%%%%%%%%%%%%%%%%%%%%%%%%%%%%%%%%%%%%%%%%%%
\subsection{From the Frobenius potential $\F_A$ to the primitive form $\zeta_A$}
In our case, the conditions for $\zeta=\zeta_A$ in Lemma~\ref{Key lemma} can be reduced to the 
assumption in the following lemma:
\begin{lem}\label{phi-psi}
Assume that we have a homogeneous change of coordinates 
$({\bf t}, e^{t_{\mu_A}})\mapsto ({\bf s},s_{\mu_A})=
 ({\bf s}({\bf t}, e^{t_{\mu_A}}),e^{t_{\mu_A}})$ on $\overline{M}$, 
\begin{subequations}\label{eq:t-s degree}
\begin{align}
E_As_1({\bf t}, e^{t_{\mu_A}})=&s_1({\bf t}, e^{t_{\mu_A}}),\\
E_As_{i,j}({\bf t}, e^{t_{\mu_A}})=&\frac{a_i-j}{a_i}s_{i,j}({\bf t}, e^{t_{\mu_A}}),\ 
i=1,2,3,\ j=1,\dots, j=a_i-1,
\end{align}
\end{subequations}
such that there exist functions 
\[
\phi_{ij}^a({\bf x};{\bf t},e^{t_{\mu_A}}), \psi_a({\bf x};{\bf t},e^{t_{\mu_A}})\in 
\Gamma(M,\O_M)[x_1,x_2,x_3],\ a=1,2,3, \ i,j=1,\dots, \mu_A
\]
satisfying the following equations$:$
\begin{subequations}
\begin{multline}
\frac{\partial F_A({\bf x};{\bf s}({\bf t},e^{t_{\mu_A}}),e^{t_{\mu_A}})}{\partial t_i}\cdot 
\frac{\partial F_A({\bf x};{\bf s}({\bf t},e^{t_{\mu_A}}),e^{t_{\mu_A}})}{\partial t_j}
-\sum_{k=1}^{\mu_A} C_{i j}^{k} \cdot \frac{\partial F_A({\bf x};{\bf s}({\bf t},e^{t_{\mu_A}}),e^{t_{\mu_A}})}{\partial t_k}\\
=\sum_{a=1}^{3} \phi_{i j}^{a}({\bf x};{\bf t},e^{t_{\mu_A}}) 
\cdot \frac{\partial F_A({\bf x};{\bf s}({\bf t},e^{t_{\mu_A}}),e^{t_{\mu_A}})}{\partial x_a},
\end{multline}
\begin{multline}
\frac{\partial^2 F_A({\bf x};{\bf s}({\bf t},e^{t_{\mu_A}}),e^{t_{\mu_A}})}{\partial t_i\partial t_j}
+\eta_{i j}-\delta_{i \mu}\cdot \frac{\partial F_A({\bf x};{\bf s}({\bf t},e^{t_{\mu_A}}),e^{t_{\mu_A}})}{\partial t_j}\\
-\delta_{j \mu}\cdot \frac{\partial F_A({\bf x};{\bf s}({\bf t},e^{t_{\mu_A}}),e^{t_{\mu_A}})}{\partial t_i}
-\sum_{a=1}^{3}\frac{\partial \phi_{i j}^{a}({\bf x};{\bf t},e^{t_{\mu_A}})}{\partial x_a}\\
=\delta_{i \mu}\cdot \delta_{j \mu}\cdot \displaystyle \sum_{a=1}^{3} \psi^a({\bf x};{\bf t},e^{t_{\mu_A}})
\cdot \frac{\partial F_A({\bf x};{\bf s}({\bf t},e^{t_{\mu_A}}),e^{t_{\mu_A}})}{\partial x_a},
\end{multline}
\begin{equation}
\sum_{a=1}^{3}\frac{\partial \psi_a ({\bf x};{\bf t},e^{t_{\mu_A}})}{\partial x_a}=1,
\end{equation}
\end{subequations}
where $(\eta^{ij}):=(\eta_A(\p_{i},\p_{j}))^{-1}$ and $C^{k}_{ij}=\displaystyle\sum_{l=1}^{\mu_A}\eta^{lk}\p_{l}\p_{i}\p_{j}\F_{A}$.
Then, we have $e=e_A$, $E=E_A$, where $e$ and $E$ 
are the primitive vector field and the Euler vector field defined in Definition~\ref{defn:vector fields}, 
and we have the equations
\[
u\nabla_{\frac{\p}{\p t_i}} \nabla_{\frac{\p}{\p t_j}}\zeta_A=
\displaystyle\sum_{k=1}^{\mu_A} C_{i j}^{k}\nabla_{\frac{\p}{\p t_k}} \zeta_A,\ i,j=1,\dots,\mu_A.
\]
In particular, the product $\circ_A$ on $\T_M$ coincides with the product $\circ$ on $\T_M$ 
defined in Definition~\ref{defn:Jacobi ring}.
\end{lem}
\begin{pf}
Note that the degrees of $t_1,\dots, t_{\mu_A-1}$ and $e^{t_{\mu_A}}$ are all positive 
since we assume $\chi_A>0$.
Therefore, the homogeneity~\eqref{eq:t-s degree} implies that
\[
e_A=\frac{\p}{\p s_1},\ E_A=s_1\frac{\p}{\p s_1}+\sum_{i=1}^3\sum_{j=1}^{a_i-1}\frac{a_i-j}{a_i}s_{i,j}\frac{\p}{\p s_{i,j}}
+\chi_A\frac{\p}{\p t_{\mu_A}}.
\]
Hence, we have $e=e_A$ and $E=E_A$. 
Based on the definition of the Gau\ss--Manin connection $\nabla$, one can show that 
the equation $u\nabla_{\p_i} \nabla_{\p_j}\zeta_A=\displaystyle\sum_{k=1}^{\mu_A} C_{i j}^{k}\nabla_{\p_k} \zeta_A$ is satisfied 
if and only if  
\begin{multline}
\left[\left(\p_{i}F_A \p_{j}F_A-\sum^{\mu}_{k=1}C^{k}_{ij} \p_{k} F_A \right)e^{-t_{\mu_A}}dx_{1}\wedge dx_{2}\wedge dx_{3}\right]\\ 
+u\left[(\p_{i}\p_{j}F_A+\eta_{ij} -\delta_{i \mu}\p_{j}F_A-\delta_{j \mu}\p_{i}F_A)e^{-t_{\mu_A}}dx_{1}\wedge dx_{2}\wedge dx_{3}\right]\\
+u^{2}\left[\delta_{i \mu} \delta_{j \mu} e^{-t_{\mu_A}}dx_{1}\wedge dx_{2}\wedge dx_{3}\right]
=0.
\end{multline}
The rest is straightforward.
\qed
\end{pf}
In our case, Lemma~\ref{phi-psi} also enables us to compare $\eta_A$ and $\eta$ defined in Corollary~\ref{cor:prim to Frob}.
\begin{lem}
Suppose that Lemma~\ref{phi-psi} holds.
Then, we have the equation 
\[
K_{F_A}^{(0)}\left(u\nabla_{\p_i}\zeta_A,u\nabla_{\p_j}\zeta_A\right)=\eta_A(\p_i,\p_j),\quad i,j=1,\dots,\mu_A.
\]
\end{lem}
\begin{pf}
It is easy to calculate that $K_{F_A}^{(0)}\left(u\nabla_{\p_i}\zeta_A,u\nabla_{\p_j}\zeta_A\right)
=\displaystyle\sum_{k=1}^{\mu_A} C_{i j}^{k}K_{F_A}^{(0)}\left(\zeta_A,u\nabla_{\p_k}\zeta_A\right)$.
We can check that $K_{F_A}^{(0)}\left(\zeta_A,u\nabla_{\p_k}\zeta_A\right)=0$ if $k\ne \mu_A$ and 
$K_{F_A}^{(0)}\left(\zeta_A,u\nabla_{\p_{\mu_A}}\zeta_A\right)=1$ by Lemma~\ref{lem:0-th residue} and 
by the first statement of Lemma~\ref{lem:flat fcn}.
\qed
\end{pf}
The functions $\phi$, $\psi$ and the homgenous change of coordinates 
in Lemma~\ref{phi-psi} should exist for any type $A$ with $\chi_A>0$
due to the existence of the primitive form by Theorem~\ref{primitive form}.
We would, conversely, reconstruct
the primitive form from the Gromov-Witten potential $\F_{0}^{\PP^{1}_A}$
via the existence of them.
Indeed, for a given $A$, they can be figured out by hand. 
We hope to give a proof of the existence of them for any $A$ with $\chi_A>0$ elsewhere in the near future. 
In the present paper, we shall consider particular cases by giving all data explicitly in Appendix:  
\begin{thm}\label{thm:another}
Let $A$ be one of $(2,3,3)$, $(2,3,4)$, $(2,3,5)$.
Then, the assumption in Lemma~\ref{phi-psi} is satisfied.
Therefore, it turns out that $\zeta_A=[e^{-t_{\mu_A}}dx_1 \wedge dx_2 \wedge dx_3]\in \H_{F_A}^{(0)}$ is a primitive form 
for the tuple $(\H_{F_A}^{(0)},\nabla, K_{F_A})$ with the minimal exponent $r=1$. 
\end{thm}
%%%%%%%%%%%%%%%%%%%%%%%%%%%%%%%%%%%%%%%%%%%%%%%%%%%%%%%%%%%%%%%%%%%%%%%%%%%%%%%
%%%%%%%%%%%%%%%%%%%%%%%%%%%%%%%%%%%%%%%%%%%%%%%%%%%%%%%%%%%%%%%%%%%%%%%%%%%%%%%
\appendix
\section{Data for the proof of Theorem~\ref{thm:another}}
%%%%%%%%%%%%%%%%%%%%%%%%%%%%%%%%%%%%%%%%%%%%%%%%%%%%%%%%%%%%%%%%%%%%%%%%%%%%%%%
\subsection{Frobenius potential for $A=(2,3,3)$}
\
{\tiny
\begin{align*}
&\F_{(2,3,3)}(t_1,t_{1,1},t_{2,1},t_{2,2},t_{3,1},t_{3,2},t_{\mu_A})\\
&=\frac{{t_1}^2 {t_{\mu_A}}}{2}+{t_1} \left (\frac{{t_{1,1}}^2}{4}+\frac{{t_{2,1}} {t_{2,2}}}{3}+\frac{{t_{3,1}} {t_{3,2}}}{3}\right)+\frac{1}{18} \left({t_{2,1}}^3+{t_{3,1}}^3\right)-\left(\frac{{t_{1,1}}^4}{96}+\frac{1}{36}\left({t_{2,1}}^2 {t_{2,2}}^2+{t_{3,1}}^2 {t_{3,2}}^2\right)\right)\\
&+\frac{1}{648}\left({t_{2,1}} {t_{2,2}}^4+{t_{3,1}} {t_{3,2}}^4\right)-\frac{1}{19440}\left({t_{2,2}}^6+{t_{3,2}}^6\right)+e^{{t_{\mu_A}}} \left({t_{1,1}} {t_{2,1}} {t_{3,1}}+\frac{1}{6}{t_{1,1}}\left( {t_{2,2}}^2 {t_{3,1}}+{t_{2,1}} {t_{3,2}}^2\right)+\frac{1}{36} {t_{1,1}} {t_{2,2}}^2 {t_{3,2}}^2\right)\\
&+e^{2 {t_{\mu_A}}} \left(\frac{1}{2} {t_{1,1}}^2 {t_{2,2}} {t_{3,2}}+\frac{1}{2}\left({t_{2,2}} {t_{3,1}}^2+{t_{2,1}}^2 {t_{3,2}}\right)+\frac{1}{6} \left({t_{2,1}} {t_{2,2}}^2 {t_{3,2}}+{t_{2,2}} {t_{3,1}} {t_{3,2}}^2\right)+\frac{1}{72}\left({t_{2,2}}^4 {t_{3,2}}+{t_{2,2}} {t_{3,2}}^4\right)\right)\\
&+e^{3 {t_{\mu_A}}} \left(\frac{{t_{1,1}}^3}{3}+{t_{1,1}} \left({t_{2,1}} {t_{2,2}}+ {t_{3,1}} {t_{3,2}}\right)+\frac{1}{6} t_{1,1}\left({t_{2,2}}^3+{t_{1,1}} {t_{3,2}}^3\right)\right)\\
&+e^{4 {t_{\mu_A}}} \left({t_{2,1}} {t_{3,1}}+\frac{1}{6}\left({t_{2,2}}^2 {t_{3,1}}+{t_{2,1}} {t_{3,2}}^2\right)+\frac{5 {t_{2,2}}^2 {t_{3,2}}^2}{18}\right)+e^{5 {t_{\mu_A}}} {t_{1,1}} {t_{2,2}} {t_{3,2}}+e^{6 {t_{\mu_A}}} \left(\frac{{t_{1,1}}^2}{2}+\frac{1}{6}\left({t_{2,2}}^3+{t_{3,2}}^3\right)\right)\\
&+\frac{1}{2} e^{8 {t_{\mu_A}}} {t_{2,2}} {t_{3,2}}+\frac{e^{12 {t_{\mu_A}}}}{12}.
\end{align*}
}
%%%%%%%%%%%%%%%%%%%%%%%%%%%%%%%%%%%%%%%%%%%%%%%%%%%%%%%%%%%%%%%%%%%%%%%%%%%%%%%
\subsection{Coordinate change for $A=(2,3,3)$}
\
{\tiny
\begin{align*}
s_1({\bf t},e^{t_{\mu_A}})=&t_1-2 e^{2 t_{\mu_A}} t_{2,2} t_{3,2}-6 e^{3 t_{\mu_A}} t_{1,1}+18 e^{6 t_{\mu_A}},\ 
s_{1,1}({\bf t},e^{t_{\mu_A}})=t_{1,1}-8 e^{3 t_{\mu_A}},\\
s_{2,1}({\bf t},e^{t_{\mu_A}})=&t_{2,1}+\frac{t_{2,2}^2}{6}-3 e^{2 t_{\mu_A}} t_{3,2},\
s_{2,2}({\bf t},e^{t_{\mu_A}})=t_{2,2},\\
s_{3,1}({\bf t},e^{t_{\mu_A}})=&t_{3,1}+\frac{t_{3,2}^2}{6}-3 e^{2 t_{\mu_A}} t_{2,2},\
s_{3,2}({\bf t},e^{t_{\mu_A}})=t_{3,2}.
\end{align*}
}
%%%%%%%%%%%%%%%%%%%%%%%%%%%%%%%%%%%%%%%%%%%%%%%%%%%%%%%%%%%%%%%%%%%%%%%%%%%%%%%
\subsection{$\phi$ for $A=(2,3,3)$}
\
{\tiny
\[
\phi ^a_{ij}=0,\ a=1,2,3,i=1,j=1,\dots,7.
\]
\[
\phi ^1_{22}=\frac{x_1}{2},\
\phi ^2_{22}=0,\
\phi ^3_{22}=0.\quad 
\phi ^1_{23}=0,\
\phi ^2_{23}=0,\
\phi ^3_{23}=-e^{t_{\mu_A}}.\\
\]
\[
\phi ^1_{24}=e^{2 t_{\mu_A}} (-2 t_{3,2}-3 x_3),\
\phi ^2_{24}=-3 e^{3 t_{\mu_A}},\
\phi ^3_{24}=e^{t_{\mu_A}} \left(-\frac{t_{2,2}}{3}-x_2\right).
\]
\[
\phi ^1_{25}=0,\
\phi ^2_{25}=-e^{t_{\mu_A}},\
\phi ^3_{25}=0.
\]
\[
\phi ^1_{26}=e^{2 t_{\mu_A}} (-2 t_{2,2}-3 x_2),\
\phi ^2_{26}=e^{t_{\mu_A}} \left(-\frac{t_{3,2}}{3}-x_3\right),\
\phi ^3_{26}=-3 e^{3 t_{\mu_A}}.
\]
\[
\phi ^1_{27}=e^{2 t_{\mu_A}} (-4 t_{2,2} t_{3,2}-6 t_{3,2} x_2-6 t_{2,2} x_3-9 x_2 x_3)+e^{3 t_{\mu_A}} (-9 t_{1,1}-12 x_1)+72 e^{6 t_{\mu_A}},
\]
\[
\phi ^2_{27}=e^{t_{\mu_A}} \left(-t_{3,1}-\frac{t_{3,2}^2}{6}-2 t_{3,2} x_3-3 x_3^2\right)-3 e^{3 t_{\mu_A}} t_{2,2},\
\phi ^3_{27}=-x_1 x_3-6 e^{3 t_{\mu_A}} t_{3,2}.
\]
\[
\phi ^1_{33}=0,\
\phi ^2_{33}=0,\
\phi ^3_{33}=0.\quad
\phi ^1_{34}=3 e^{3 t_{\mu_A}},\
\phi ^2_{34}=\frac{x_2}{3},\
\phi ^3_{34}=0.
\]
\[
\phi ^1_{35}=-e^{t_{\mu_A}},\
\phi ^2_{35}=0,\
\phi ^3_{35}=0.\quad 
\phi ^1_{36}=e^{t_{\mu_A}} \left(-\frac{t_{3,2}}{3}-x_3\right),\
\phi ^2_{36}=-2 e^{2 t_{\mu_A}},\
\phi ^3_{36}=0.
\]
\[
\phi ^1_{37}=-x_1 x_2+e^{3 t_{\mu_A}} (-6 t_{2,2}-18 x_2),\
\phi ^2_{37}=e^{2 t_{\mu_A}} (-4 t_{3,2}-6 x_3),\
\phi ^3_{37}=e^{t_{\mu_A}} (-t_{1,1}-2 x_1)-4 e^{4 t_{\mu_A}}.
\]
\[
\phi ^1_{44}=-\frac{1}{3} x_1 x_2-2 e^{3 t_{\mu_A}} t_{2,2},\
\phi ^2_{44}=\frac{x_2^2}{3}+e^{2 t_{\mu_A}} \left(-\frac{4 t_{3,2}}{3}-2 x_3\right),
\phi ^3_{44}=e^{t_{\mu_A}} \left(-\frac{t_{1,1}}{3}-\frac{2 x_1}{3}\right)+\frac{8}{3} e^{4 t_{\mu_A}}.
\]
\[
\phi ^1_{45}=e^{t_{\mu_A}} \left(-\frac{t_{2,2}}{3}-x_2\right),\
\phi ^2_{45}=0,\
\phi ^3_{45}=-2 e^{2 t_{\mu_A}}.
\]
\[
\phi ^1_{46}=e^{t_{\mu_A}} \left(-\frac{1}{9} t_{2,2} t_{3,2}-\frac{1}{3} t_{3,2} x_2-\frac{1}{3} t_{2,2} x_3-x_2 x_3\right)-e^{2 t_{\mu_A}} t_{1,1}-e^{5 t_{\mu_A}},
\]
\[
\phi ^2_{46}=e^{2 t_{\mu_A}} \left(-\frac{2 t_{2,2}}{3}-x_2\right),\
\phi ^3_{46}=e^{2 t_{\mu_A}} \left(-\frac{2 t_{3,2}}{3}-x_3\right).
\]
\[
\phi ^1_{47}=-\frac{1}{3} t_{2,2} x_1 x_2-x_1 x_2^2+e^{2 t_{\mu_A}} (-2 t_{1,1} t_{3,2}-3 t_{1,1} x_3-3 x_1 x_3)-2 e^{3 t_{\mu_A}} t_{2,2}^2+e^{5 t_{\mu_A}} (10 t_{3,2}+42 x_3),
\]
\[
\phi ^2_{47}=e^{2 t_{\mu_A}} \left(-\frac{4}{3} t_{2,2} t_{3,2}-2 t_{3,2} x_2-2 t_{2,2} x_3\right)+e^{3 t_{\mu_A}} (-9 t_{1,1}-6 x_1)+36 e^{6 t_{\mu_A}},
\]
\[
\phi ^3_{47}=e^{t_{\mu_A}} \left(-\frac{1}{3} t_{1,1} t_{2,2}-\frac{2}{3} t_{2,2} x_1-t_{1,1} x_2-2 x_1 x_2\right)+e^{2 t_{\mu_A}} \left(-2 t_{3,1}-\frac{t_{3,2}^2}{3}\right)+e^{4 t_{\mu_A}} \left(\frac{14 t_{2,2}}{3}+14 x_2\right).
\]
\[
\phi ^1_{55}=0,\
\phi ^2_{55}=0,\
\phi ^3_{55}=0.\quad 
\phi ^1_{56}=3 e^{3 t_{\mu_A}},\
\phi ^2_{56}=0,\
\phi ^3_{56}=\frac{x_3}{3}.
\]
\[
\phi ^1_{57}=-x_1 x_3+e^{3 t_{\mu_A}} (-6 t_{3,2}-18 x_3),\
\phi ^2_{57}=e^{t_{\mu_A}} (-t_{1,1}-2 x_1)-4 e^{4 t_{\mu_A}},\
\phi ^3_{57}=e^{2 t_{\mu_A}} (-4 t_{2,2}-6 x_2).
\]
\[
\phi ^1_{66}=-\frac{1}{3} x_1 x_3-2 e^{3 t_{\mu_A}} t_{3,2},\
\phi ^2_{66}=e^{t_{\mu_A}} \left(-\frac{t_{1,1}}{3}-\frac{2 x_1}{3}\right)+\frac{8}{3} e^{4 t_{\mu_A}},
\phi ^3_{66}=\frac{x_3^2}{3}+e^{2 t_{\mu_A}} \left(-\frac{4 t_{2,2}}{3}-2 x_2\right).
\]
\[
\phi ^1_{67}=-\frac{1}{3} t_{3,2} x_1 x_3-x_1 x_3^2+e^{2 t_{\mu_A}} (-2 t_{1,1} t_{2,2}-3 t_{1,1} x_2-3 x_1 x_2)-2 e^{3 t_{\mu_A}} t_{3,2}^2+e^{5 t_{\mu_A}} (10 t_{2,2}+42 x_2),
\]
\[
\phi ^2_{67}=e^{t_{\mu_A}} \left(-\frac{1}{3} t_{1,1} t_{3,2}-\frac{2}{3} t_{3,2} x_1-t_{1,1} x_3-2 x_1 x_3\right)+e^{2 t_{\mu_A}} \left(-2 t_{2,1}-\frac{t_{2,2}^2}{3}\right)+e^{4 t_{\mu_A}} \left(\frac{14 t_{3,2}}{3}+14 x_3\right),
\]
\[
\phi ^3_{67}=e^{2 t_{\mu_A}} \left(-\frac{4}{3} t_{2,2} t_{3,2}-2 t_{3,2} x_2-2 t_{2,2} x_3\right)+e^{3 t_{\mu_A}} (-9 t_{1,1}-6 x_1)+36 e^{6 t_{\mu_A}}.
\]

\begin{align*}
\phi ^1_{77}=&-e^{-t_{\mu_A}} x_1^2 x_2 x_3+e^{2 t_{\mu_A}} \left(-4 t_{1,1} t_{2,2} t_{3,2}-6 t_{1,1} t_{3,2} x_2-12 t_{3,2} x_1 x_2-6 t_{1,1} t_{2,2} x_3\right.\\
&\left.-12 t_{2,2} x_1 x_3-9 t_{1,1} x_2 x_3-54 x_1 x_2 x_3\right)+e^{3 t_{\mu_A}} \left(-9 t_{1,1}^2+24 t_{2,2}^2 x_2+72 t_{2,2} x_2^2\right.\\
&\left.+54 x_2^3+24 t_{3,2}^2 x_3+72 t_{3,2} x_3^2+54 x_3^3\right)+e^{5 t_{\mu_A}} \left(8 t_{2,2} t_{3,2}+120 t_{3,2} x_2+120 t_{2,2} x_3\right.\\
&\left.+180 x_2 x_3\right)+e^{6 t_{\mu_A}} \left(252 t_{1,1}+288 x_1\right)-1440 e^{9 t_{\mu_A}},\\
\phi ^2_{77}=&e^{t_{\mu_A}} \left(-t_{1,1} t_{3,1}-\frac{1}{6} t_{1,1} t_{3,2}^2-2 t_{3,1} x_1-\frac{1}{3} t_{3,2}^2 x_1-2 t_{1,1} t_{3,2} x_3-4 t_{3,2} x_1 x_3\right.\\
&-3 t_{1,1} x_3^2-6 x_1 x_3^2\biggr)+e^{2 t_{\mu_A}} \left(-4 t_{2,1} t_{3,2}-\frac{2}{3} t_{2,2}^2 t_{3,2}-6 t_{2,1} x_3-t_{2,2}^2 x_3\right.\\
&+12 t_{2,2} x_2 x_3+18 x_2^2 x_3\biggr)+e^{3 t_{\mu_A}} (-21 t_{1,1} t_{2,2}-18 t_{2,2} x_1-18 t_{1,1} x_2)\\
&+e^{4 t_{\mu_A}} \left(-52 t_{3,1}+\frac{10 t_{3,2}^2}{3}-14 t_{3,2} x_3-48 x_3^2\right)+e^{6 t_{\mu_A}} (204 t_{2,2}+144 x_2),\\
\phi ^3_{77}=&-t_{1,1} x_1 x_3-2 x_1^2 x_3+e^{2 t_{\mu_A}} \left(-4 t_{2,2} t_{3,1}-\frac{2}{3} t_{2,2} t_{3,2}^2-6 t_{3,1} x_2-t_{3,2}^2 x_2\right.\\
&+12 t_{3,2} x_2 x_3+18 x_2 x_3^2\biggr)+e^{3 t_{\mu_A}} (-24 t_{1,1} t_{3,2}-24 t_{3,2} x_1-18 t_{1,1} x_3-52 x_1 x_3)\\
&+e^{4 t_{\mu_A}} \left(12 t_{2,2}^2+90 t_{2,2} x_2+108 x_2^2\right)+e^{6 t_{\mu_A}} (48 t_{3,2}+144 x_3).
\end{align*}
}
%%%%%%%%%%%%%%%%%%%%%%%%%%%%%%%%%%%%%%%%%%%%%%%%%%%%%%%%%%%%%%%%%%%%%%%%%%%%%%%
\subsection{$\psi$ for $A=(2,3,3)$}
\
{\tiny
\[
\psi ^1=x_1+18 e^{3 t_{\mu_A}},\ \psi ^2=0,\ \psi ^3=0.
\]
}
%%%%%%%%%%%%%%%%%%%%%%%%%%%%%%%%%%%%%%%%%%%%%%%%%%%%%%%%%%%%%%%%%%%%%%%%%%%%%%%
\subsection{Frobenius potential for $A=(2,3,4)$}
\ 
{\tiny
\begin{align*}
&\mathcal{F}_{(2,3,4)} (t_1,t_{1,1},t_{2,1},t_{2,2},t_{3,1},t_{3,2},t_{3,3},t_{\mu_A})\\
&=\frac{{t_1}^2 {t_{\mu_A}}}{2}+{t_1} \left(\frac{{t_{1,1}}^2}{4}+\frac{{t_{2,1}} {t_{2,2}}}{3}+\frac{{t_{3,2}}^2}{8}+\frac{{t_{3,1}} {t_{3,3}}}{4}\right)+\left(\frac{{t_{2,1}}^3}{18}+\frac{{t_{3,1}}^2 {t_{3,2}}}{8}\right)\\
&+\left(-\frac{{t_{1,1}}^4}{96}-\frac{{t_{2,1}}^2 {t_{2,2}}^2}{36}-\frac{{t_{3,2}}^4}{192}-\frac{1}{32} {t_{3,1}} {t_{3,2}}^2 {t_{3,3}}-\frac{{t_{3,1}}^2 {t_{3,3}}^2}{64}\right)+\left(\frac{{t_{2,1}} {t_{2,2}}^4}{648}+\frac{{t_{3,2}}^3 {t_{3,3}}^2}{384}+\frac{1}{384} {t_{3,1}} {t_{3,2}} {t_{3,3}}^3\right)\\
&+\left(-\frac{{t_{2,2}}^6}{19440}-\frac{{t_{3,2}}^2 {t_{3,3}}^4}{3072}-\frac{{t_{3,1}} {t_{3,3}}^5}{30720}\right)+\frac{{t_{3,2}} {t_{3,3}}^6}{73728}-\frac{{t_{3,3}}^8}{4128768}+e^{{t_{\mu_A}}} \left({t_{1,1}} {t_{2,1}} {t_{3,1}} +\frac{1}{6} {t_{1,1}} {t_{2,2}}^2 {t_{3,1}}+\frac{1}{4} {t_{1,1}} {t_{2,1}} {t_{3,2}} {t_{3,3}}\right.\\
&\left. +\frac{1}{24} {t_{1,1}} {t_{2,2}}^2 {t_{3,2}} {t_{3,3}}+\frac{1}{96} {t_{1,1}} {t_{2,1}} {t_{3,3}}^3+\frac{1}{576} {t_{1,1}} {t_{2,2}}^2 {t_{3,3}}^3\right)+e^{2 {t_{\mu_A}}} \left(\frac{{t_{2,2}} {t_{3,1}}^2}{2}+\frac{{t_{2,1}}^2 {t_{3,2}}}{2}+\frac{1}{2} {t_{1,1}}^2 {t_{2,2}} {t_{3,2}}\right.\\
&+\frac{1}{6} {t_{2,1}} {t_{2,2}}^2 {t_{3,2}}+\frac{{t_{2,2}}^4 {t_{3,2}}}{72}+\frac{1}{4} {t_{2,2}} {t_{3,1}} {t_{3,2}} {t_{3,3}}+\frac{{t_{2,1}}^2 {t_{3,3}}^2}{8}+\frac{1}{8} {t_{1,1}}^2 {t_{2,2}} {t_{3,3}}^2+\frac{1}{24} {t_{2,1}} {t_{2,2}}^2 {t_{3,3}}^2+\frac{{t_{2,2}}^4 {t_{3,3}}^2}{288}\\
&\left.+\frac{1}{32} {t_{2,2}} {t_{3,2}}^2 {t_{3,3}}^2+\frac{1}{96} {t_{2,2}} {t_{3,1}} {t_{3,3}}^3+\frac{1}{384} {t_{2,2}} {t_{3,2}} {t_{3,3}}^4+\frac{{t_{2,2}} {t_{3,3}}^6}{18432}\right)+e^{3 {t_{\mu_A}}} \left({t_{1,1}} {t_{3,1}} {t_{3,2}}+\frac{{t_{1,1}}^3 {t_{3,3}}}{3}\right.\\
&\left.+{t_{1,1}} {t_{2,1}} {t_{2,2}} {t_{3,3}}+\frac{1}{6} {t_{1,1}} {t_{2,2}}^3 {t_{3,3}}+\frac{1}{4} {t_{1,1}} {t_{3,2}}^2 {t_{3,3}}+\frac{1}{4} {t_{1,1}} {t_{3,1}} {t_{3,3}}^2+\frac{7}{96} {t_{1,1}} {t_{3,2}} {t_{3,3}}^3+\frac{{t_{1,1}} {t_{3,3}}^5}{384}\right)\\
&+e^{4 {t_{\mu_A}}} \left({t_{1,1}}^2 {t_{2,1}}+\frac{{t_{2,1}}^2 {t_{2,2}}}{2}+\frac{2 {t_{1,1}}^2 {t_{2,2}}^2}{3}+\frac{{t_{2,1}} {t_{2,2}}^3}{6}+\frac{{t_{2,2}}^5}{72}+\frac{{t_{2,2}}^2 {t_{3,2}}^2}{4}+{t_{2,1}} {t_{3,1}} {t_{3,3}}+\frac{1}{6} {t_{2,2}}^2 {t_{3,1}} {t_{3,3}}\right.\\
&\left.+\frac{1}{4} {t_{2,1}} {t_{3,2}} {t_{3,3}}^2+\frac{1}{6} {t_{2,2}}^2 {t_{3,2}} {t_{3,3}}^2+\frac{{t_{2,1}} {t_{3,3}}^4}{96}+\frac{5 {t_{2,2}}^2 {t_{3,3}}^4}{288}\right)+e^{5 {t_{\mu_A}}} \left({t_{1,1}} {t_{2,2}} {t_{3,1}}+\frac{5}{4} {t_{1,1}} {t_{2,2}} {t_{3,2}} {t_{3,3}}+\frac{25}{96} {t_{1,1}} {t_{2,2}} {t_{3,3}}^3\right)\\
&+e^{6 {t_{\mu_A}}} \left(\frac{{t_{3,1}}^2}{2}+\frac{{t_{1,1}}^2 {t_{3,2}}}{2}+{t_{2,1}} {t_{2,2}} {t_{3,2}}+\frac{{t_{2,2}}^3 {t_{3,2}}}{6}+\frac{{t_{3,2}}^3}{6}+\frac{{t_{3,1}} {t_{3,2}} {t_{3,3}}}{4}+\frac{5 {t_{1,1}}^2 {t_{3,3}}^2}{8}+\frac{1}{4} {t_{2,1}} {t_{2,2}} {t_{3,3}}^2\right.\\
&\left.+\frac{5 {t_{2,2}}^3 {t_{3,3}}^2}{24}+\frac{5 {t_{3,2}}^2 {t_{3,3}}^2}{32}+\frac{{t_{3,1}} {t_{3,3}}^3}{96}+\frac{13 {t_{3,2}} {t_{3,3}}^4}{384}+\frac{49 {t_{3,3}}^6}{18432}\right)+e^{7 {t_{\mu_A}}} \left({t_{1,1}} {t_{2,1}} {t_{3,3}}+\frac{7}{6} {t_{1,1}} {t_{2,2}}^2 {t_{3,3}}\right)\\
&+e^{8 {t_{\mu_A}}} \left(\frac{{t_{2,1}}^2}{2}+{t_{1,1}}^2 {t_{2,2}}+\frac{{t_{2,1}} {t_{2,2}}^2}{6}+\frac{5 {t_{2,2}}^4}{36}+\frac{1}{2} {t_{2,2}} {t_{3,2}} {t_{3,3}}^2+\frac{{t_{2,2}} {t_{3,3}}^4}{8}\right)+e^{9 {t_{\mu_A}}} \left({t_{1,1}} {t_{3,2}} {t_{3,3}}+\frac{{t_{1,1}} {t_{3,3}}^3}{4}\right)\\
&+e^{10 {t_{\mu_A}}} \left(\frac{{t_{2,2}}^2 {t_{3,2}}}{2}+\frac{5 {t_{2,2}}^2 {t_{3,3}}^2}{8}\right)+e^{11 {t_{\mu_A}}} {t_{1,1}} {t_{2,2}} {t_{3,3}}+e^{12 {t_{\mu_A}}} \left(\frac{{t_{1,1}}^2}{2}+\frac{{t_{2,2}}^3}{6}+\frac{{t_{3,2}}^2}{4}+\frac{{t_{3,2}} {t_{3,3}}^2}{8}+\frac{19 {t_{3,3}}^4}{192}\right)\\
&+\frac{1}{2} e^{14 {t_{\mu_A}}} {t_{2,2}} {t_{3,3}}^2+\frac{1}{4} e^{16 {t_{\mu_A}}} {t_{2,2}}^2+\frac{1}{6} e^{18 {t_{\mu_A}}} {t_{3,3}}^2+\frac{e^{24 {t_{\mu_A}}}}{24}.
\end{align*}
}
%%%%%%%%%%%%%%%%%%%%%%%%%%%%%%%%%%%%%%%%%%%%%%%%%%%%%%%%%%%%%%%%%%%%%%%%%%%%%%%
\subsection{Coordinate change for $A=(2,3,4)$}
\
{\tiny
\begin{align*}
s_1({\bf t})=&t_1+e^{2 t_{\mu_A}} \left(-2 t_{2,2} t_{3,2}-\frac{1}{2} t_{2,2} t_{3,3}^2\right)-6 e^{3 t_{\mu_A}} t_{1,1} t_{3,3}-12 e^{4 t_{\mu_A}} t_{2,1}\\
&+e^{6 t_{\mu_A}} \left(30 t_{3,2}+\frac{51 t_{3,3}^2}{2}\right)+52 e^{8 t_{\mu_A}}
t_{2,2}-702 e^{12 t_{\mu_A}},\\
s_{1,1}({\bf t})=&t_{1,1}-8 e^{3 t_{\mu_A}} t_{3,3},\\
s_{2,1}({\bf t})=&t_{2,1}+\frac{t_{2,2}^2}{6}+e^{2 t_{\mu_A}} \left(-3 t_{3,2}-\frac{3 t_{3,3}^2}{4}\right)-15 e^{4 t_{\mu_A}} t_{2,2}+240 e^{8 t_{\mu_A}},\
s_{2,2}({\bf t})=t_{2,2}-27 e^{4 t_{\mu_A}},\\
s_{3,1}({\bf t})=&t_{3,1}+\frac{1}{4} t_{3,2} t_{3,3}+\frac{t_{3,3}^3}{96}-3 e^{2 t_{\mu_A}} t_{2,2} t_{3,3}-8 e^{3 t_{\mu_A}} t_{1,1}+45 e^{6 t_{\mu_A}} t_{3,3},\\
s_{3,2}({\bf t})=&t_{3,2}+\frac{t_{3,3}^2}{4}-4 e^{2 t_{\mu_A}} t_{2,2}+28 e^{6 t_{\mu_A}},\ s_{3,3}({\bf t})=t_{3,3}.
\end{align*}
}
%%%%%%%%%%%%%%%%%%%%%%%%%%%%%%%%%%%%%%%%%%%%%%%%%%%%%%%%%%%%%%%%%%%%%%%%%%%%%%%
\subsection{$\phi$ for $A=(2,3,4)$}
\
{\tiny
\[
\phi ^a_{ij}=0,\ a=1,2,3,i=1,j=1,\dots,8.
\]
\[
\phi^1_{22}=\frac{x_1}{2},\
\phi^2_{22}=16 e^{4 t_{\mu_A}},\
\phi^3_{22}=0.\quad
\phi^1_{23}=8 e^{4 t_{\mu_A}},\
\phi^2_{23}=0,\
\phi^3_{23}=-e^{t_{\mu_A}}.
\]
\[
\phi^1_{24}=e^{2 t_{\mu_A}} \left(-2 t_{3,2}-\frac{t_{3,3}^2}{2}-3 t_{3,3} x_3-4 x_3^2\right)+e^{4 t_{\mu_A}} \left(\frac{8 t_{2,2}}{3}-4 x_2\right)+40 e^{8 t_{\mu_A}},
\]
\[
\phi^2_{24}=e^{3 t_{\mu_A}} (-3 t_{3,3}-4 x_3),\
\phi^3_{24}=e^{t_{\mu_A}} \left(-\frac{t_{2,2}}{3}-x_2\right)+7 e^{5 t_{\mu_A}}.
\]
\[
\phi^1_{25}=0,\
\phi^2_{25}=-e^{t_{\mu_A}},\
\phi^3_{25}=0.
\]
\[
\phi^1_{26}=e^{2 t_{\mu_A}} (-2 t_{2,2}-3 x_2)+30 e^{6 t_{\mu_A}},\
\phi^2_{26}=e^{t_{\mu_A}} \left(-\frac{t_{3,3}}{4}-x_3\right),\
\phi^3_{26}=-3 e^{3 t_{\mu_A}}.
\]
\[
\phi^1_{27}=e^{2 t_{\mu_A}} \left(-t_{2,2} t_{3,3}-\frac{3}{2} t_{3,3} x_2-2 t_{2,2} x_3-3 x_2 x_3\right)+e^{3 t_{\mu_A}} (-3 t_{1,1}-4 x_1)+e^{6 t_{\mu_A}} (39 t_{3,3}+54 x_3),
\]
\[
\phi^2_{27}=e^{t_{\mu_A}} \left(-\frac{t_{3,2}}{4}-\frac{t_{3,3}^2}{32}-\frac{1}{2} t_{3,3} x_3-x_3^2\right)-e^{3 t_{\mu_A}} t_{2,2}+-e^{7 t_{\mu_A}},\phi^3_{27}=e^{3 t_{\mu_A}} \left(-\frac{3 t_{3,3}}{2}-2 x_3\right).
\]
\begin{align*}
\phi^1_{28}=&e^{2 t_{\mu_A}} \left(-4 t_{2,2} t_{3,2}-t_{2,2} t_{3,3}^2-6 t_{3,2} x_2-\frac{3}{2} t_{3,3}^2 x_2-6 t_{2,2} t_{3,3} x_3-9 t_{3,3} x_2 x_3\right.\\
&-8 t_{2,2} x_3^2-12 x_2 x_3^2\biggr)+e^{3 t_{\mu_A}} (-9 t_{1,1} t_{3,3}-12 t_{3,3} x_1-12 t_{1,1} x_3-16 x_1 x_3)\\
&+e^{6 t_{\mu_A}} \left(84 t_{3,2}+93 t_{3,3}^2+294 t_{3,3} x_3+264 x_3^2\right)+e^{8 t_{\mu_A}} (32 t_{2,2}+48 x_2)-384 e^{12 t_{\mu_A}},\\
\phi^2_{28}=&e^{t_{\mu_A}} \left(-t_{3,1}-\frac{1}{4} t_{3,2} t_{3,3}-\frac{t_{3,3}^3}{96}-2 t_{3,2} x_3-\frac{1}{2} t_{3,3}^2 x_3-3 t_{3,3} x_3^2-4 x_3^3\right)\\
&-3 e^{3 t_{\mu_A}} t_{2,2} t_{3,3}+e^{4 t_{\mu_A}} (-8 t_{1,1}-32 x_1)+e^{7 t_{\mu_A}} (209 t_{3,3}+304 x_3),\\
\phi^3_{28}=&-x_1 x_3+e^{3 t_{\mu_A}} \left(-6 t_{3,2}-\frac{3 t_{3,3}^2}{2}\right)+e^{5 t_{\mu_A}} (-4 t_{2,2}+12 x_2)-96 e^{9 t_{\mu_A}}.
\end{align*}
\[
\phi^1_{33}=0,\
\phi^2_{33}=0,\
\phi^3_{33}=0.\quad
\phi^1_{34}=e^{3 t_{\mu_A}} (3 t_{3,3}+4 x_3),\
\phi^2_{34}=\frac{x_2}{3}+8 e^{4 t_{\mu_A}},\
\phi^3_{34}=0.
\]
\[
\phi^1_{35}=-e^{t_{\mu_A}},\
\phi^2_{35}=0,\
\phi^3_{35}=0.\quad
\phi^1_{36}=e^{t_{\mu_A}} \left(-\frac{t_{3,3}}{4}-x_3\right),\
\phi^2_{36}=-2 e^{2 t_{\mu_A}},\
\phi^3_{36}=0.
\]
\[
\phi^1_{37}=e^{t_{\mu_A}} \left(-\frac{t_{3,2}}{4}-\frac{t_{3,3}^2}{32}-\frac{1}{2} t_{3,3} x_3-x_3^2\right)+e^{3 t_{\mu_A}} (-t_{2,2}-6 x_2)+63 e^{7 t_{\mu_A}},
\]
\[
\phi^2_{37}=e^{2 t_{\mu_A}} (-t_{3,3}-2 x_3),\
\phi^3_{37}=-4 e^{4 t_{\mu_A}}.
\]
\[
\phi^1_{38}=-x_1 x_2+e^{3 t_{\mu_A}} (-6 t_{2,2} t_{3,3}-18 t_{3,3} x_2-8 t_{2,2} x_3-24 x_2 x_3)+e^{7 t_{\mu_A}} (246 t_{3,3}+264 x_3),
\]
\[
\phi^2_{38}=e^{2 t_{\mu_A}} \left(-4 t_{3,2}-t_{3,3}^2-6 t_{3,3} x_3-8 x_3^2\right)-36 e^{4 t_{\mu_A}} x_2+416 e^{8 t_{\mu_A}},
\]
\[
\phi^3_{38}=e^{t_{\mu_A}} (-t_{1,1}-2 x_1)-4 e^{4 t_{\mu_A}} t_{3,3}.
\]
\[
\phi^1_{44}=-\frac{1}{3} x_1 x_2+e^{3 t_{\mu_A}} \left(-2 t_{2,2} t_{3,3}-\frac{8}{3} t_{2,2} x_3\right)+e^{7 t_{\mu_A}} (18 t_{3,3}+24 x_3),
\]
\[
\phi^2_{44}=\frac{x_2^2}{3}+e^{2 t_{\mu_A}} \left(-\frac{4 t_{3,2}}{3}-\frac{t_{3,3}^2}{3}-2 t_{3,3} x_3-\frac{8 x_3^2}{3}\right)-4 e^{4 t_{\mu_A}} x_2+\frac{32}{3} e^{8 t_{\mu_A}},
\]
\[
\phi^3_{44}=e^{t_{\mu_A}} \left(-\frac{t_{1,1}}{3}-\frac{2 x_1}{3}\right)+e^{4 t_{\mu_A}} \left(\frac{8 t_{3,3}}{3}+4 x_3\right).
\]
\[
\phi^1_{45}=e^{t_{\mu_A}} \left(-\frac{t_{2,2}}{3}-x_2\right)+15 e^{5 t_{\mu_A}},\
\phi^2_{45}=0,\
\phi^3_{45}=-2 e^{2 t_{\mu_A}}.
\]
\[
\phi^1_{46}=e^{t_{\mu_A}} \left(-\frac{1}{12} t_{2,2} t_{3,3}-\frac{1}{4} t_{3,3} x_2-\frac{1}{3} t_{2,2} x_3-x_2 x_3\right)-e^{2 t_{\mu_A}} t_{1,1}+e^{5 t_{\mu_A}} \left(\frac{11 t_{3,3}}{4}+3 x_3\right),
\]
\[
\phi^2_{46}=e^{2 t_{\mu_A}} \left(-\frac{2 t_{2,2}}{3}-x_2\right)+6 e^{6 t_{\mu_A}},\
\phi^3_{46}=e^{2 t_{\mu_A}} \left(-\frac{t_{3,3}}{2}-x_3\right).
\]
\begin{align*}
\phi^1_{47}=&e^{t_{\mu_A}} \left(-\frac{1}{12} t_{2,2} t_{3,2}-\frac{1}{96} t_{2,2} t_{3,3}^2-\frac{1}{4} t_{3,2} x_2-\frac{1}{32} t_{3,3}^2 x_2-\frac{1}{6} t_{2,2} t_{3,3} x_3
-\frac{1}{2} t_{3,3} x_2 x_3-\frac{1}{3} t_{2,2} x_3^2-x_2 x_3^2\right)\\
&+e^{3 t_{\mu_A}} \left(-\frac{1}{3} t_{2,2}^2+t_{2,2} x_2\right)+e^{5 t_{\mu_A}} \left(\frac{31 t_{3,2}}{4}+\frac{31 t_{3,3}^2}{32}+\frac{31}{2} t_{3,3} x_3+23 x_3^2\right)+e^{7 t_{\mu_A}} (-12 t_{2,2}-3 x_2)+31 e^{11 t_{\mu_A}},
\end{align*}
\[
\phi^2_{47}=e^{2 t_{\mu_A}} \left(-\frac{1}{3} t_{2,2} t_{3,3}-\frac{1}{2} t_{3,3} x_2-\frac{2}{3} t_{2,2} x_3\right)+e^{3 t_{\mu_A}} (-3 t_{1,1}-2 x_1)+e^{6 t_{\mu_A}} (15 t_{3,3}+14 x_3),
\]
\[
\phi^3_{47}=e^{2 t_{\mu_A}} \left(-\frac{t_{3,2}}{2}-\frac{t_{3,3}^2}{16}-\frac{1}{2} t_{3,3} x_3-x_3^2\right)+e^{4 t_{\mu_A}} \left(\frac{2 t_{2,2}}{3}+2 x_2\right)-18 e^{8 t_{\mu_A}}.
\]
\begin{align*}
\phi^1_{48}=&-\frac{1}{3} t_{2,2} x_1 x_2-x_1 x_2^2+e^{2 t_{\mu_A}} \left(-2 t_{1,1} t_{3,2}-\frac{1}{2} t_{1,1} t_{3,3}^2-3 t_{1,1} t_{3,3} x_3-3 t_{3,3} x_1 x_3-4 t_{1,1} x_3^2-4 x_1 x_3^2\right)\\
&+e^{3 t_{\mu_A}} \left(-2 t_{2,2}^2 t_{3,3}-\frac{8}{3} t_{2,2}^2 x_3\right)+e^{4 t_{\mu_A}} (-8 t_{1,1} t_{2,2}-12 t_{1,1} x_2+27 x_1 x_2)\\
&+e^{5 t_{\mu_A}} \left(10 t_{3,2} t_{3,3}+\frac{5 t_{3,3}^3}{2}-24 t_{3,2} x_3+36 t_{3,3}^2 x_3+56 t_{3,3} x_3^2\right)\\
&+e^{7 t_{\mu_A}} (56 t_{2,2} t_{3,3}+36 t_{3,3} x_2+160 t_{2,2} x_3+48 x_2 x_3)+280 e^{8 t_{\mu_A}} t_{1,1}+e^{11 t_{\mu_A}} (-1130 t_{3,3}-1464 x_3),\\
\phi^2_{48}=&e^{2 t_{\mu_A}} \left(-\frac{4}{3} t_{2,2} t_{3,2}-\frac{1}{3} t_{2,2} t_{3,3}^2-2 t_{3,2} x_2-\frac{1}{2} t_{3,3}^2 x_2-2 t_{2,2} t_{3,3} x_3-\frac{8}{3} t_{2,2} x_3^2\right)\\
&+e^{3 t_{\mu_A}} (-9 t_{1,1} t_{3,3}-6 t_{3,3} x_1-12 t_{1,1} x_3-8 x_1 x_3)+e^{4 t_{\mu_A}} \left(-8 t_{2,1}-\frac{4 t_{2,2}^2}{3}-8 t_{2,2} x_2-36 x_2^2\right)\\
&+e^{6 t_{\mu_A}} \left(84 t_{3,2}+57 t_{3,3}^2+210 t_{3,3} x_3+216 x_3^2\right)+e^{8 t_{\mu_A}} \left(-\frac{88 t_{2,2}}{3}+532 x_2\right)-1824 e^{12 t_{\mu_A}},\\
\phi^3_{48}=&e^{t_{\mu_A}} \left(-\frac{1}{3} t_{1,1} t_{2,2}-\frac{2}{3} t_{2,2} x_1-t_{1,1} x_2-2 x_1 x_2\right)+e^{2 t_{\mu_A}} \left(-2 t_{3,1}-\frac{1}{2} t_{3,2} t_{3,3}-\frac{t_{3,3}^3}{48}\right)\\
&+e^{4 t_{\mu_A}} \left(\frac{14}{3} t_{2,2} t_{3,3}+14 t_{3,3} x_2+8 t_{2,2} x_3\right)+e^{5 t_{\mu_A}} (19 t_{1,1}+54 x_1)+e^{8 t_{\mu_A}} (-270 t_{3,3}-168 x_3).
\end{align*}
\[
\phi^1_{55}=0,\
\phi^2_{55}=0,\
\phi^3_{55}=0.\quad
\phi^1_{56}=3 e^{3 t_{\mu_A}},\
\phi^2_{56}=0,\
\phi^3_{56}=0.
\]
\[
\phi^1_{57}=\frac{3}{2} e^{3 t_{\mu_A}} t_{3,3},\
\phi^2_{57}=8 e^{4 t_{\mu_A}},\
\phi^3_{57}=\frac{x_3}{4}.
\]
\[
\phi^1_{58}=-x_1 x_3+e^{3 t_{\mu_A}} \left(-6 t_{3,2}-\frac{3 t_{3,3}^2}{2}-18 t_{3,3} x_3-24 x_3^2\right)+e^{5 t_{\mu_A}} (12 t_{2,2}-36 x_2)+336 e^{9 t_{\mu_A}},
\]
\[
\phi^2_{58}=e^{t_{\mu_A}} (-t_{1,1}-2 x_1)+e^{4 t_{\mu_A}} (-4 t_{3,3}-48 x_3),
\]
\[
\phi^3_{58}=e^{2 t_{\mu_A}} (-4 t_{2,2}-6 x_2)+36 e^{6 t_{\mu_A}}.
\]
\[
\phi^1_{66}=e^{3 t_{\mu_A}} \left(\frac{3 t_{3,3}}{2}+6 x_3\right),\
\phi^2_{66}=12 e^{4 t_{\mu_A}},\
\phi^3_{66}=\frac{x_3}{4}.
\]
\[
\phi^1_{67}=-\frac{1}{4} x_1 x_3+e^{3 t_{\mu_A}} \left(-\frac{9 t_{3,2}}{4}-\frac{9 t_{3,3}^2}{32}-\frac{3}{2} t_{3,3} x_3-3 x_3^2\right)+e^{5 t_{\mu_A}} (7 t_{2,2}+6 x_2)-57 e^{9 t_{\mu_A}},
\]
\[
\phi^2_{67}=e^{t_{\mu_A}} \left(-\frac{t_{1,1}}{4}-\frac{x_1}{2}\right)+e^{4 t_{\mu_A}} (t_{3,3}+2 x_3),
\]
\[
\phi^3_{67}=\frac{x_3^2}{4}+e^{2 t_{\mu_A}} \left(-t_{2,2}-\frac{3 x_2}{2}\right)+15 e^{6 t_{\mu_A}}.
\]
\begin{align*}
\phi^1_{68}=&-\frac{1}{4} t_{3,3} x_1 x_3-x_1 x_3^2+e^{2 t_{\mu_A}} (-2 t_{1,1} t_{2,2}-3 t_{1,1} x_2-3 x_1 x_2)\\
&+e^{3 t_{\mu_A}} \left(-\frac{3}{2} t_{3,2} t_{3,3}-\frac{3 t_{3,3}^3}{8}+6 t_{3,2} x_3-3 t_{3,3}^2 x_3-6 t_{3,3} x_3^2\right)\\
&+e^{5 t_{\mu_A}} (13 t_{2,2} t_{3,3}+33 t_{3,3} x_2-12 t_{2,2} x_3+36 x_2 x_3)+18 e^{6 t_{\mu_A}} t_{1,1}+e^{9 t_{\mu_A}} (-150 t_{3,3}-120 x_3),\\
\phi^2_{68}=&e^{t_{\mu_A}} \left(-\frac{1}{4} t_{1,1} t_{3,3}-\frac{1}{2} t_{3,3} x_1-t_{1,1} x_3-2 x_1 x_3\right)+e^{2 t_{\mu_A}} \left(-2 t_{2,1}-\frac{t_{2,2}^2}{3}\right)\\
&+e^{4 t_{\mu_A}} \left(6 t_{3,2}+\frac{t_{3,3}^2}{2}+2 t_{3,3} x_3-24 x_3^2\right)+e^{6 t_{\mu_A}} (54 t_{2,2}+108 x_2)-1056 e^{10 t_{\mu_A}},
\end{align*}
\[
\phi^3_{68}=e^{2 t_{\mu_A}} \left(-t_{2,2} t_{3,3}-\frac{3}{2} t_{3,3} x_2-2 t_{2,2} x_3\right)+e^{3 t_{\mu_A}} (-9 t_{1,1}-6 x_1)+e^{6 t_{\mu_A}} (45 t_{3,3}+42 x_3).
\]
\begin{align*}
\phi^1_{77}=&-\frac{1}{16} t_{3,3} x_1 x_3-\frac{1}{4} x_1 x_3^2+e^{2 t_{\mu_A}} \left(-\frac{1}{2} t_{1,1} t_{2,2}-\frac{3}{4} t_{1,1} x_2-\frac{3}{2} x_1 x_2\right)+e^{3 t_{\mu_A}} \left(-\frac{3}{32} t_{3,3}^3+\frac{3}{8} t_{3,3}^2 x_3+\frac{3}{2} t_{3,3} x_3^2\right)\\
&+e^{5 t_{\mu_A}} (2 t_{2,2} t_{3,3}+12 t_{3,3} x_2+8 t_{2,2} x_3+12 x_2 x_3)+e^{6 t_{\mu_A}} \left(\frac{51 t_{1,1}}{2}+32 x_1\right)+e^{9 t_{\mu_A}} (-198 t_{3,3}-216 x_3),\\
\phi^2_{77}=&e^{t_{\mu_A}} \left(-\frac{1}{16} t_{1,1} t_{3,3}-\frac{1}{8} t_{3,3} x_1-\frac{1}{4} t_{1,1} x_3-\frac{1}{2} x_1 x_3\right)+e^{2 t_{\mu_A}} \left(-\frac{t_{2,1}}{2}-\frac{t_{2,2}^2}{12}\right)\\
&+e^{4 t_{\mu_A}} \left(-\frac{t_{3,2}}{2}+\frac{5 t_{3,3}^2}{8}+4 t_{3,3} x_3+4 x_3^2\right)-\frac{1}{2} e^{6 t_{\mu_A}} t_{2,2},\\
\phi^3_{77}=&\frac{1}{16} t_{3,3} x_3^2+\frac{x_3^3}{4}+e^{2 t_{\mu_A}} \left(-\frac{1}{4} t_{2,2} t_{3,3}-\frac{3}{8} t_{3,3} x_2-t_{2,2} x_3\right)+e^{3 t_{\mu_A}} (-3 t_{1,1}-3 x_1)+e^{6 t_{\mu_A}} \left(\frac{75 t_{3,3}}{4}+19 x_3\right).
\end{align*}
\begin{align*}
\phi^1_{78}=&-\frac{1}{4} t_{3,2} x_1 x_3-\frac{1}{32} t_{3,3}^2 x_1 x_3-\frac{1}{2} t_{3,3} x_1 x_3^2-x_1 x_3^3+e^{2 t_{\mu_A}} \left(-t_{1,1} t_{2,2} t_{3,3}-\frac{3}{2} t_{1,1} t_{3,3} x_2-\frac{3}{2} t_{3,3} x_1 x_2\right.\\
&-2 t_{1,1} t_{2,2} x_3-3 t_{2,2} x_1 x_3-3 t_{1,1} x_2 x_3-12 x_1 x_2 x_3\Bigr)+e^{3 t_{\mu_A}} \left(-3 t_{1,1}^2-\frac{3 t_{3,2}^2}{2}-\frac{9}{16} t_{3,2} t_{3,3}^2\right.\\
&\left.-\frac{3 t_{3,3}^4}{64}+8 t_{2,2}^2 x_2+24 t_{2,2} x_2^2+18 x_2^3-\frac{3}{2} t_{3,2} t_{3,3} x_3+\frac{3}{16} t_{3,3}^3 x_3+\frac{3}{4} t_{3,3}^2 x_3^2\right)\\
&+e^{5 t_{\mu_A}} \left(5 t_{2,2} t_{3,2}+\frac{47}{8} t_{2,2} t_{3,3}^2+3 t_{3,2} x_2+\frac{183}{8} t_{3,3}^2 x_2+52 t_{2,2} t_{3,3} x_3+78 t_{3,3} x_2 x_3+52 t_{2,2} x_3^2+60 x_2 x_3^2\right)\\
&+e^{6 t_{\mu_A}} (93 t_{1,1} t_{3,3}+96 t_{3,3} x_1+114 t_{1,1} x_3+249 x_1 x_3)+e^{7 t_{\mu_A}} \left(-84 t_{2,2}^2-540 t_{2,2} x_2-648 x_2^2\right)\\
&+e^{9 t_{\mu_A}} \left(102 t_{3,2}-582 t_{3,3}^2-1674 t_{3,3} x_3-1128 x_3^2\right)+e^{11 t_{\mu_A}} (2780 t_{2,2}+7284 x_2)-26064 e^{15 t_{\mu_A}},\\
\phi^2_{78}=&e^{t_{\mu_A}} \left(-\frac{1}{4} t_{1,1} t_{3,2}-\frac{1}{32} t_{1,1} t_{3,3}^2-\frac{1}{2} t_{3,2} x_1-\frac{1}{16} t_{3,3}^2 x_1-\frac{1}{2} t_{1,1} t_{3,3} x_3-t_{3,3} x_1 x_3-t_{1,1} x_3^2-2 x_1 x_3^2\right)\\
&+e^{2 t_{\mu_A}} \left(-t_{2,1} t_{3,3}-\frac{1}{6} t_{2,2}^2 t_{3,3}-2 t_{2,1} x_3-\frac{1}{3} t_{2,2}^2 x_3+4 t_{2,2} x_2 x_3+6 x_2^2 x_3\right)+e^{3 t_{\mu_A}} (-7 t_{1,1} t_{2,2}-6 t_{2,2} x_1-6 t_{1,1} x_2)\\
&+e^{4 t_{\mu_A}} \left(-16 t_{3,1}-2 t_{3,2} t_{3,3}+\frac{11 t_{3,3}^3}{24}-26 t_{3,2} x_3+2 t_{3,3}^2 x_3-28 t_{3,3} x_3^2-64 x_3^3\right)\\
&+e^{6 t_{\mu_A}} (95 t_{2,2} t_{3,3}+102 t_{3,3} x_2+198 t_{2,2} x_3-108 x_2 x_3)+e^{7 t_{\mu_A}} (413 t_{1,1}+498 x_1)+e^{10 t_{\mu_A}} (-3228 t_{3,3}-2288 x_3),
\end{align*}
\begin{align*}
\phi^3_{78}=&e^{2 t_{\mu_A}} \left(-t_{2,2} t_{3,2}-\frac{1}{8} t_{2,2} t_{3,3}^2-\frac{3}{2} t_{3,2} x_2-\frac{3}{16} t_{3,3}^2 x_2-t_{2,2} t_{3,3} x_3-2 t_{2,2} x_3^2\right)\\
&+e^{3 t_{\mu_A}} \left(-\frac{9}{2} t_{1,1} t_{3,3}-3 t_{3,3} x_1-6 t_{1,1} x_3-16 x_1 x_3\right)+e^{4 t_{\mu_A}} \left(4 t_{2,2}^2+30 t_{2,2} x_2+36 x_2^2\right)\\
&+e^{6 t_{\mu_A}} \left(-15 t_{3,2}+\frac{105 t_{3,3}^2}{8}+57 t_{3,3} x_3+42 x_3^2\right)+e^{8 t_{\mu_A}} (-240 t_{2,2}-666 x_2)+3036 e^{12 t_{\mu_A}}.
\end{align*}
\begin{align*}
\phi^1_{88}=&-e^{-t_{\mu_A}} x_1^2 x_2 x_3+e^{2 t_{\mu_A}} \left(-4 t_{1,1} t_{2,2} t_{3,2}-t_{1,1} t_{2,2} t_{3,3}^2-6 t_{1,1} t_{3,2} x_2-\frac{3}{2} t_{1,1} t_{3,3}^2 x_2-12 t_{3,2} x_1 x_2-3 t_{3,3}^2 x_1 x_2\right.\\
&-6 t_{1,1} t_{2,2} t_{3,3} x_3-12 t_{2,2} t_{3,3} x_1 x_3-9 t_{1,1} t_{3,3} x_2 x_3-54 t_{3,3} x_1 x_2 x_3-8 t_{1,1} t_{2,2} x_3^2-16 t_{2,2} x_1 x_3^2-12 t_{1,1} x_2 x_3^2\\
&-72 x_1 x_2 x_3^2\biggr)+e^{3 t_{\mu_A}} \biggl(-9 t_{1,1}^2 t_{3,3}+24 t_{2,2}^2 t_{3,3} x_2+72 t_{2,2} t_{3,3} x_2^2+54 t_{3,3} x_2^3-12 t_{1,1}^2 x_3+24 t_{3,2}^2 x_3\\
&+12 t_{3,2} t_{3,3}^2 x_3+\frac{3}{2} t_{3,3}^4 x_3-24 t_{1,1} x_1 x_3+32 t_{2,2}^2 x_2 x_3+96 t_{2,2} x_2^2 x_3+72 x_2^3 x_3+72 t_{3,2} t_{3,3} x_3^2\\
&+18 t_{3,3}^3 x_3^2+96 t_{3,2} x_3^3+78 t_{3,3}^2 x_3^3+144 t_{3,3} x_3^4+96 x_3^5\biggr)+e^{5 t_{\mu_A}} \left(8 t_{2,2} t_{3,2} t_{3,3}+2 t_{2,2} t_{3,3}^3+120 t_{3,2} t_{3,3} x_2+30 t_{3,3}^3 x_2\right.\\
&\left.-96 t_{2,2} t_{3,2} x_3+96 t_{2,2} t_{3,3}^2 x_3+288 t_{3,2} x_2 x_3+252 t_{3,3}^2 x_2 x_3+160 t_{2,2} t_{3,3} x_3^2+672 t_{3,3} x_2 x_3^2+576 x_2 x_3^3\right)\\
&+e^{4 t_{\mu_A}} \left(-16 t_{1,1} t_{2,2}^2+24 t_{1,1} t_{2,2} x_2+24 t_{2,2} x_1 x_2+72 t_{1,1} x_2^2-72 x_1 x_2^2\right)\\
&+e^{6 t_{\mu_A}} \left(-36 t_{1,1} t_{3,2}+243 t_{1,1} t_{3,3}^2+288 t_{3,3}^2 x_1+714 t_{1,1} t_{3,3} x_3+1404 t_{3,3} x_1 x_3+504 t_{1,1} x_3^2+1104 x_1 x_3^2\right)\\
&+e^{7 t_{\mu_A}} \left(-304 t_{2,2}^2 t_{3,3}-2352 t_{2,2} t_{3,3} x_2-2520 t_{3,3} x_2^2-192 t_{2,2}^2 x_3-2880 t_{2,2} x_2 x_3-2592 x_2^2 x_3\right)\\
&+e^{8 t_{\mu_A}} (224 t_{1,1} t_{2,2}-1608 t_{1,1} x_2+816 x_1 x_2)+e^{9 t_{\mu_A}} \left(216 t_{3,2} t_{3,3}-1386 t_{3,3}^3-1824 t_{3,2} x_3-6912 t_{3,3}^2 x_3\right.\\
&\left.-11808 t_{3,3} x_3^2-7680 x_3^3\right)+e^{11 t_{\mu_A}} (13808 t_{2,2} t_{3,3}+35400 t_{3,3} x_2+17472 t_{2,2} x_3\\
&+33120 x_2 x_3)+10224 e^{12 t_{\mu_A}} t_{1,1}+e^{15 t_{\mu_A}} (-160704 t_{3,3}-149376 x_3),
\end{align*}
\begin{align*}
\phi^2_{88}=&-t_{1,1} x_1 x_2-2 x_1^2 x_2,+e^{2 t_{\mu_A}} \left(-4 t_{2,1} t_{3,2}-\frac{2}{3} t_{2,2}^2 t_{3,2}-t_{2,1} t_{3,3}^2-\frac{1}{6} t_{2,2}^2 t_{3,3}^2-6 t_{2,1} t_{3,3} x_3-t_{2,2}^2 t_{3,3} x_3\right.\\
&\left.+12 t_{2,2} t_{3,3} x_2 x_3+18 t_{3,3} x_2^2 x_3-8 t_{2,1} x_3^2-\frac{4}{3} t_{2,2}^2 x_3^2+16 t_{2,2} x_2 x_3^2+24 x_2^2 x_3^2\right)+e^{3 t_{\mu_A}} (-24 t_{1,1} t_{2,2} t_{3,3}\\
&-24 t_{2,2} t_{3,3} x_1-18 t_{1,1} t_{3,3} x_2-52 t_{3,3} x_1 x_2-32 t_{1,1} t_{2,2} x_3-32 t_{2,2} x_1 x_3-144 x_1 x_2 x_3)\\
&+e^{4 t_{\mu_A}} \left(-48 t_{1,1}^2-16 t_{2,1} t_{2,2}-\frac{8 t_{2,2}^3}{3}+12 t_{3,2}^2+6 t_{3,2} t_{3,3}^2+\frac{3 t_{3,3}^4}{4}-48 t_{1,1} x_1\right.\\
&\left.+90 t_{3,2} t_{3,3} x_3+\frac{45}{2} t_{3,3}^3 x_3+120 t_{3,2} x_3^2+138 t_{3,3}^2 x_3^2+288 t_{3,3} x_3^3+192 x_3^4\right)\\
&+e^{6 t_{\mu_A}} \left(204 t_{2,2} t_{3,2}+99 t_{2,2} t_{3,3}^2+432 t_{3,2} x_2+252 t_{3,3}^2 x_2+298 t_{2,2} t_{3,3} x_3-324 t_{3,3} x_2 x_3+312 t_{2,2} x_3^2-432 x_2 x_3^2\right)\\
&+e^{7 t_{\mu_A}} (1872 t_{1,1} t_{3,3}+2808 t_{3,3} x_1+2112 t_{1,1} x_3+2208 x_1 x_3)+e^{8 t_{\mu_A}} \left(2096 t_{2,1}-\frac{536 t_{2,2}^2}{3}+1728 t_{2,2} x_2+3888 x_2^2\right)\\
&+e^{10 t_{\mu_A}} \left(-8976 t_{3,2}-13428 t_{3,3}^2-20928 t_{3,3} x_3-10752 x_3^2\right)+e^{12 t_{\mu_A}} (-10704 t_{2,2}-68256 x_2)+295680 e^{16 t_{\mu_A}},
\end{align*}
\begin{align*}
\phi^3_{88}=&e^{t_{\mu_A}} \left(-t_{1,1} t_{2,1}-\frac{1}{6} t_{1,1} t_{2,2}^2-2 t_{2,1} x_1-\frac{1}{3} t_{2,2}^2 x_1-2 t_{1,1} t_{2,2} x_2-4 t_{2,2} x_1 x_2-3 t_{1,1} x_2^2-6 x_1 x_2^2\right)+e^{2 t_{\mu_A}} \Biggl(-4 t_{2,2} t_{3,1}\\
&\left.-t_{2,2} t_{3,2} t_{3,3}-\frac{1}{24} t_{2,2} t_{3,3}^3-6 t_{3,1} x_2-\frac{3}{2} t_{3,2} t_{3,3} x_2-\frac{1}{16} t_{3,3}^3 x_2+12 t_{3,2} x_2 x_3+3 t_{3,3}^2 x_2 x_3+18 t_{3,3} x_2 x_3^2+24 x_2 x_3^3\right)\\
&+e^{3 t_{\mu_A}} \left(-21 t_{1,1} t_{3,2}-\frac{21}{4} t_{1,1} t_{3,3}^2-18 t_{3,2} x_1-\frac{9}{2} t_{3,3}^2 x_1-18 t_{1,1} t_{3,3} x_3-24 t_{1,1} x_3^2\right)+e^{4 t_{\mu_A}} \biggl(-52 t_{2,1} t_{3,3}\\
&\left.+\frac{10}{3} t_{2,2}^2 t_{3,3}-14 t_{2,2} t_{3,3} x_2-48 t_{3,3} x_2^2-48 t_{2,1} x_3+8 t_{2,2}^2 x_3-48 t_{2,2} x_2 x_3\right)+e^{5 t_{\mu_A}} (23 t_{1,1} t_{2,2}+78 t_{2,2} x_1\\
&+30 t_{1,1} x_2-36 x_1 x_2)+e^{6 t_{\mu_A}} \left(-36 t_{3,1}+195 t_{3,2} t_{3,3}+\frac{405 t_{3,3}^3}{8}+144 t_{3,2} x_3+180 t_{3,3}^2 x_3+192 t_{3,3} x_3^2\right)\\
&+e^{8 t_{\mu_A}} (-156 t_{2,2} t_{3,3}+1314 t_{3,3} x_2+48 t_{2,2} x_3+336 x_2 x_3)+e^{9 t_{\mu_A}} (240 t_{1,1}+1152 x_1)+e^{12 t_{\mu_A}} (-9300 t_{3,3}-4464 x_3).
\end{align*}
}
%%%%%%%%%%%%%%%%%%%%%%%%%%%%%%%%%%%%%%%%%%%%%%%%%%%%%%%%%%%%%%%%%%%%%%%%%%%%%%%
\subsection{$\psi$ for $A=(2,3,4)$}
\
{\tiny
\[
\psi ^1=x_1+e^{3 t_{\mu_A}} (18 t_{3,3}+48 x_3),\ \psi ^2=-48 e^{4 t_{\mu_A}}, \ \psi ^3=0.
\]
}
%%%%%%%%%%%%%%%%%%%%%%%%%%%%%%%%%%%%%%%%%%%%%%%%%%%%%%%%%%%%%%%%%%%%%%%%%%%%%%%
\subsection{Frobenius potential for $A=(2,3,5)$}
\
{\tiny
\begin{align*}
&\mathcal{F}(t_1,t_{1,1},t_{2,1},t_{2,2},t_{3,1},t_{3,2},t_{3,3},t_{3,4},t_{\mu_A})\\
&=\frac{{t_1}^2 {t_{\mu_A}}}{2}+{t_1} \left(\frac{{t_{1,1}}^2}{4}+\frac{{t_{2,1}} {t_{2,2}}}{3}+\frac{{t_{3,2}} {t_{3,3}}}{5}+\frac{{t_{3,1}} {t_{3,4}}}{5}\right)+\left(\frac{{t_{2,1}}^3}{18}+\frac{{t_{3,1}} {t_{3,2}}^2}{10}+\frac{{t_{3,1}}^2 {t_{3,3}}}{10}\right)\\
&+\left(-\frac{{t_{1,1}}^4}{96}-\frac{{t_{2,1}}^2 {t_{2,2}}^2}{36}-\frac{{t_{3,2}}^2 {t_{3,3}}^2}{50}-\frac{{t_{3,1}} {t_{3,3}}^3}{150}-\frac{{t_{3,2}}^3 {t_{3,4}}}{150}-\frac{1}{25} {t_{3,1}} {t_{3,2}} {t_{3,3}} {t_{3,4}}-\frac{{t_{3,1}}^2 {t_{3,4}}^2}{100}\right)\\
&+\left(\frac{{t_{2,1}} {t_{2,2}}^4}{648}+\frac{{t_{3,3}}^5}{3000}+\frac{1}{250} {t_{3,2}} {t_{3,3}}^3 {t_{3,4}}+\frac{1}{250} {t_{3,2}}^2 {t_{3,3}} {t_{3,4}}^2+\frac{1}{500} {t_{3,1}} {t_{3,3}}^2 {t_{3,4}}^2+\frac{1}{750} {t_{3,1}} {t_{3,2}} {t_{3,4}}^3\right)\\
&+\left(-\frac{{t_{2,2}}^6}{19440}-\frac{3 {t_{3,3}}^4 {t_{3,4}}^2}{10000}-\frac{{t_{3,2}} {t_{3,3}}^2 {t_{3,4}}^3}{1500}-\frac{{t_{3,2}}^2 {t_{3,4}}^4}{7500}-\frac{{t_{3,1}} {t_{3,3}} {t_{3,4}}^4}{15000}\right)+\left(\frac{7 {t_{3,3}}^3 {t_{3,4}}^4}{150000}+\frac{11 {t_{3,2}} {t_{3,3}} {t_{3,4}}^5}{375000}+\frac{{t_{3,1}} {t_{3,4}}^6}{2250000}\right)\\
&+\left(-\frac{13 {t_{3,3}}^2 {t_{3,4}}^6}{4500000}-\frac{{t_{3,2}} {t_{3,4}}^7}{3150000}\right)+\frac{7 {t_{3,3}} {t_{3,4}}^8}{90000000}-\frac{7 {t_{3,4}}^{10}}{8100000000}+e^{{t_{\mu_A}}} \left({t_{1,1}} {t_{2,1}} {t_{3,1}}+\frac{1}{6} {t_{1,1}} {t_{2,2}}^2 {t_{3,1}}+\frac{1}{10} {t_{1,1}} {t_{2,1}} {t_{3,3}}^2\right.\\
&+\frac{1}{60} {t_{1,1}} {t_{2,2}}^2 {t_{3,3}}^2+\frac{1}{5} {t_{1,1}} {t_{2,1}} {t_{3,2}} {t_{3,4}}+\frac{1}{30} {t_{1,1}} {t_{2,2}}^2 {t_{3,2}} {t_{3,4}}+\frac{1}{50} {t_{1,1}} {t_{2,1}} {t_{3,3}} {t_{3,4}}^2+\frac{1}{300} {t_{1,1}} {t_{2,2}}^2 {t_{3,3}} {t_{3,4}}^2\\
&\left.+\frac{{t_{1,1}} {t_{2,1}} {t_{3,4}}^4}{3000}+\frac{{t_{1,1}} {t_{2,2}}^2 {t_{3,4}}^4}{18000}\right)+e^{2 {t_{\mu_A}}} \left(\frac{{t_{2,2}} {t_{3,1}}^2}{2}+\frac{{t_{2,1}}^2 {t_{3,2}}}{2}+\frac{1}{2} {t_{1,1}}^2 {t_{2,2}} {t_{3,2}}+\frac{1}{6} {t_{2,1}} {t_{2,2}}^2 {t_{3,2}}+\frac{{t_{2,2}}^4 {t_{3,2}}}{72}\right.\\
&+\frac{1}{10} {t_{2,2}} {t_{3,1}} {t_{3,3}}^2+\frac{{t_{2,2}} {t_{3,3}}^4}{200}+\frac{1}{5} {t_{2,2}} {t_{3,1}} {t_{3,2}} {t_{3,4}}+\frac{1}{5} {t_{2,1}}^2 {t_{3,3}} {t_{3,4}}+\frac{1}{5} {t_{1,1}}^2 {t_{2,2}} {t_{3,3}} {t_{3,4}}+\frac{1}{15} {t_{2,1}} {t_{2,2}}^2 {t_{3,3}} {t_{3,4}}\\
&+\frac{1}{180} {t_{2,2}}^4 {t_{3,3}} {t_{3,4}}+\frac{1}{50} {t_{2,2}} {t_{3,2}} {t_{3,3}}^2 {t_{3,4}}+\frac{1}{50} {t_{2,2}} {t_{3,2}}^2 {t_{3,4}}^2+\frac{1}{50} {t_{2,2}} {t_{3,1}} {t_{3,3}} {t_{3,4}}^2+\frac{1}{500} {t_{2,2}} {t_{3,3}}^3 {t_{3,4}}^2+\frac{{t_{2,1}}^2 {t_{3,4}}^3}{75}\\
&+\frac{1}{75} {t_{1,1}}^2 {t_{2,2}} {t_{3,4}}^3+\frac{1}{225} {t_{2,1}} {t_{2,2}}^2 {t_{3,4}}^3+\frac{{t_{2,2}}^4 {t_{3,4}}^3}{2700}+\frac{1}{250} {t_{2,2}} {t_{3,2}} {t_{3,3}} {t_{3,4}}^3+\frac{{t_{2,2}} {t_{3,1}} {t_{3,4}}^4}{3000}+\frac{7 {t_{2,2}} {t_{3,3}}^2 {t_{3,4}}^4}{30000}\\
&\left.+\frac{{t_{2,2}} {t_{3,2}} {t_{3,4}}^5}{15000}+\frac{{t_{2,2}} {t_{3,3}} {t_{3,4}}^6}{150000}+\frac{{t_{2,2}} {t_{3,4}}^8}{18000000}\right)+e^{3 {t_{\mu_A}}} \left({t_{1,1}} {t_{3,1}} {t_{3,2}}+\frac{{t_{1,1}}^3 {t_{3,3}}}{3}+{t_{1,1}} {t_{2,1}} {t_{2,2}} {t_{3,3}}+\frac{1}{6} {t_{1,1}} {t_{2,2}}^3 {t_{3,3}}\right.\\
&+\frac{1}{10} {t_{1,1}} {t_{3,2}} {t_{3,3}}^2+\frac{1}{5} {t_{1,1}} {t_{3,2}}^2 {t_{3,4}}+\frac{2}{5} {t_{1,1}} {t_{3,1}} {t_{3,3}} {t_{3,4}}+\frac{1}{25} {t_{1,1}} {t_{3,3}}^3 {t_{3,4}}+\frac{{t_{1,1}}^3 {t_{3,4}}^2}{10}+\frac{3}{10} {t_{1,1}} {t_{2,1}} {t_{2,2}} {t_{3,4}}^2\\
&+\frac{1}{20} {t_{1,1}} {t_{2,2}}^3 {t_{3,4}}^2+\frac{1}{10} {t_{1,1}} {t_{3,2}} {t_{3,3}} {t_{3,4}}^2+\frac{2}{75} {t_{1,1}} {t_{3,1}} {t_{3,4}}^3+\frac{4}{375} {t_{1,1}} {t_{3,3}}^2 {t_{3,4}}^3+\frac{17 {t_{1,1}} {t_{3,2}} {t_{3,4}}^4}{3000}+\frac{{t_{1,1}} {t_{3,3}} {t_{3,4}}^5}{1500}\\
&\left.+\frac{{t_{1,1}} {t_{3,4}}^7}{112500}\right)+e^{4 {t_{\mu_A}}} \left(\frac{{t_{2,2}}^2 {t_{3,2}}^2}{4}+{t_{2,1}} {t_{3,1}} {t_{3,3}}+\frac{1}{6} {t_{2,2}}^2 {t_{3,1}} {t_{3,3}}+\frac{{t_{2,1}} {t_{3,3}}^3}{10}+\frac{{t_{2,2}}^2 {t_{3,3}}^3}{60}+{t_{1,1}}^2 {t_{2,1}} {t_{3,4}}\right.\\
&+\frac{1}{2} {t_{2,1}}^2 {t_{2,2}} {t_{3,4}}+\frac{2}{3} {t_{1,1}}^2 {t_{2,2}}^2 {t_{3,4}}+\frac{1}{6} {t_{2,1}} {t_{2,2}}^3 {t_{3,4}}+\frac{{t_{2,2}}^5 {t_{3,4}}}{72}+\frac{1}{5} {t_{2,1}} {t_{3,2}} {t_{3,3}} {t_{3,4}}+\frac{7}{30} {t_{2,2}}^2 {t_{3,2}} {t_{3,3}} {t_{3,4}}\\
&+\frac{3}{10} {t_{2,1}} {t_{3,1}} {t_{3,4}}^2+\frac{1}{20} {t_{2,2}}^2 {t_{3,1}} {t_{3,4}}^2+\frac{1}{20} {t_{2,1}} {t_{3,3}}^2 {t_{3,4}}^2+\frac{29}{600} {t_{2,2}}^2 {t_{3,3}}^2 {t_{3,4}}^2+\frac{3}{50} {t_{2,1}} {t_{3,2}} {t_{3,4}}^3+\frac{7}{300} {t_{2,2}}^2 {t_{3,2}} {t_{3,4}}^3\\
&\left.+\frac{19 {t_{2,1}} {t_{3,3}} {t_{3,4}}^4}{3000}+\frac{23 {t_{2,2}}^2 {t_{3,3}} {t_{3,4}}^4}{3600}+\frac{{t_{2,1}} {t_{3,4}}^6}{10000}+\frac{7 {t_{2,2}}^2 {t_{3,4}}^6}{36000}\right)+e^{5 {t_{\mu_A}}} \left({t_{1,1}} {t_{2,1}}^2+{t_{1,1}}^3 {t_{2,2}}+\frac{4}{3} {t_{1,1}} {t_{2,1}} {t_{2,2}}^2\right.\\
&+\frac{7 {t_{1,1}} {t_{2,2}}^4}{36}+{t_{1,1}} {t_{2,2}} {t_{3,2}} {t_{3,3}}+{t_{1,1}} {t_{2,2}} {t_{3,1}} {t_{3,4}}+\frac{1}{2} {t_{1,1}} {t_{2,2}} {t_{3,3}}^2 {t_{3,4}}+\frac{1}{2} {t_{1,1}} {t_{2,2}} {t_{3,2}} {t_{3,4}}^2+\frac{1}{6} {t_{1,1}} {t_{2,2}} {t_{3,3}} {t_{3,4}}^3\\
&\left.+\frac{1}{120} {t_{1,1}} {t_{2,2}} {t_{3,4}}^5\right)+e^{6 {t_{\mu_A}}}\left({t_{1,1}}^2 {t_{3,1}}+{t_{2,1}} {t_{2,2}} {t_{3,1}}+\frac{{t_{2,2}}^3 {t_{3,1}}}{6}+\frac{{t_{3,2}}^3}{6}+\frac{3 {t_{1,1}}^2 {t_{3,3}}^2}{5}+\frac{1}{10} {t_{2,1}} {t_{2,2}} {t_{3,3}}^2\right.\\
&+\frac{11 {t_{2,2}}^3 {t_{3,3}}^2}{60}+\frac{{t_{3,1}}^2 {t_{3,4}}}{2}+\frac{7}{10} {t_{1,1}}^2 {t_{3,2}} {t_{3,4}}+\frac{6}{5} {t_{2,1}} {t_{2,2}} {t_{3,2}} {t_{3,4}}+\frac{1}{5} {t_{2,2}}^3 {t_{3,2}} {t_{3,4}}+\frac{1}{5} {t_{3,2}}^2 {t_{3,3}} {t_{3,4}}+\frac{1}{10} {t_{3,1}} {t_{3,3}}^2 {t_{3,4}}\\
&+\frac{{t_{3,3}}^4 {t_{3,4}}}{200}+\frac{1}{5} {t_{3,1}} {t_{3,2}} {t_{3,4}}^2+\frac{13}{25} {t_{1,1}}^2 {t_{3,3}} {t_{3,4}}^2+\frac{21}{50} {t_{2,1}} {t_{2,2}} {t_{3,3}} {t_{3,4}}^2+\frac{17}{100} {t_{2,2}}^3 {t_{3,3}} {t_{3,4}}^2+\frac{1}{10} {t_{3,2}} {t_{3,3}}^2 {t_{3,4}}^2\\
&+\frac{{t_{3,2}}^2 {t_{3,4}}^3}{30}+\frac{1}{50} {t_{3,1}} {t_{3,3}} {t_{3,4}}^3+\frac{19 {t_{3,3}}^3 {t_{3,4}}^3}{1500}+\frac{22 {t_{1,1}}^2 {t_{3,4}}^4}{375}+\frac{27 {t_{2,1}} {t_{2,2}} {t_{3,4}}^4}{1000}+\frac{39 {t_{2,2}}^3 {t_{3,4}}^4}{2000}+\frac{11}{750} {t_{3,2}} {t_{3,3}} {t_{3,4}}^4\\
&\left.+\frac{{t_{3,1}} {t_{3,4}}^5}{3000}+\frac{71 {t_{3,3}}^2 {t_{3,4}}^5}{30000}+\frac{19 {t_{3,2}} {t_{3,4}}^6}{45000}+\frac{67 {t_{3,3}} {t_{3,4}}^7}{450000}+\frac{521 {t_{3,4}}^9}{162000000}\right)+e^{7 {t_{\mu_A}}} \biggl({t_{1,1}} {t_{2,1}} {t_{3,2}}+\frac{7}{6} {t_{1,1}} {t_{2,2}}^2 {t_{3,2}}\\
&\left.+\frac{7}{5} {t_{1,1}} {t_{2,1}} {t_{3,3}} {t_{3,4}}+\frac{49}{30} {t_{1,1}} {t_{2,2}}^2 {t_{3,3}} {t_{3,4}}+\frac{49}{150} {t_{1,1}} {t_{2,1}} {t_{3,4}}^3+\frac{343}{900} {t_{1,1}} {t_{2,2}}^2 {t_{3,4}}^3\right)+e^{8 {t_{\mu_A}}}\biggl({t_{2,2}} {t_{3,1}} {t_{3,2}}+2 {t_{1,1}}^2 {t_{2,2}} {t_{3,3}}
\end{align*}
\begin{align*}
&+{t_{2,1}} {t_{2,2}}^2 {t_{3,3}}+\frac{{t_{2,2}}^4 {t_{3,3}}}{6}+\frac{3}{5} {t_{2,2}} {t_{3,2}} {t_{3,3}}^2+\frac{1}{5} {t_{2,2}} {t_{3,2}}^2 {t_{3,4}}+\frac{2}{5} {t_{2,2}} {t_{3,1}} {t_{3,3}} {t_{3,4}}+\frac{6}{25} {t_{2,2}} {t_{3,3}}^3 {t_{3,4}}+\frac{{t_{2,1}}^2 {t_{3,4}}^2}{2}\\
&+\frac{8}{5} {t_{1,1}}^2 {t_{2,2}} {t_{3,4}}^2+\frac{7}{15} {t_{2,1}} {t_{2,2}}^2 {t_{3,4}}^2+\frac{17 {t_{2,2}}^4 {t_{3,4}}^2}{90}+\frac{2}{5} {t_{2,2}} {t_{3,2}} {t_{3,3}} {t_{3,4}}^2+\frac{2}{75} {t_{2,2}} {t_{3,1}} {t_{3,4}}^3+\frac{18}{125} {t_{2,2}} {t_{3,3}}^2 {t_{3,4}}^3\\
&\left.+\frac{19}{375} {t_{2,2}} {t_{3,2}} {t_{3,4}}^4+\frac{2}{75} {t_{2,2}} {t_{3,3}} {t_{3,4}}^5+\frac{34 {t_{2,2}} {t_{3,4}}^7}{28125}\right)+e^{9 {t_{\mu_A}}} \biggr({t_{1,1}} {t_{3,1}} {t_{3,3}}+\frac{{t_{1,1}} {t_{3,3}}^3}{10}+{t_{1,1}}^3 {t_{3,4}}+3 {t_{1,1}} {t_{2,1}} {t_{2,2}} {t_{3,4}}\\
&+\frac{3}{2} {t_{1,1}} {t_{2,2}}^3 {t_{3,4}}+\frac{6}{5} {t_{1,1}} {t_{3,2}} {t_{3,3}} {t_{3,4}}+\frac{3}{10} {t_{1,1}} {t_{3,1}} {t_{3,4}}^2+\frac{9}{20} {t_{1,1}} {t_{3,3}}^2 {t_{3,4}}^2+\frac{9}{25} {t_{1,1}} {t_{3,2}} {t_{3,4}}^3+\frac{153 {t_{1,1}} {t_{3,3}} {t_{3,4}}^4}{1000}\\
&\left.+\frac{81 {t_{1,1}} {t_{3,4}}^6}{10000}\right)+e^{10 {t_{\mu_A}}} \biggl({t_{1,1}}^2 {t_{2,1}}+{t_{2,1}}^2 {t_{2,2}}+\frac{13 {t_{1,1}}^2 {t_{2,2}}^2}{6}+\frac{{t_{2,1}} {t_{2,2}}^3}{3}+\frac{23 {t_{2,2}}^5}{180}+{t_{2,1}} {t_{3,2}} {t_{3,3}}+\frac{1}{6} {t_{2,2}}^2 {t_{3,2}} {t_{3,3}}\\
&+{t_{2,1}} {t_{3,1}} {t_{3,4}}+\frac{1}{6} {t_{2,2}}^2 {t_{3,1}} {t_{3,4}}+\frac{1}{2} {t_{2,1}} {t_{3,3}}^2 {t_{3,4}}+\frac{7}{12} {t_{2,2}}^2 {t_{3,3}}^2 {t_{3,4}}+\frac{1}{2} {t_{2,1}} {t_{3,2}} {t_{3,4}}^2+\frac{7}{12} {t_{2,2}}^2 {t_{3,2}} {t_{3,4}}^2\\
&\left.+\frac{1}{6} {t_{2,1}} {t_{3,3}} {t_{3,4}}^3+\frac{19}{36} {t_{2,2}}^2 {t_{3,3}} {t_{3,4}}^3+\frac{{t_{2,1}} {t_{3,4}}^5}{120}+\frac{43 {t_{2,2}}^2 {t_{3,4}}^5}{720}\right)+e^{11 {t_{\mu_A}}} \left({t_{1,1}} {t_{2,2}} {t_{3,1}}+\frac{11}{10} {t_{1,1}} {t_{2,2}} {t_{3,3}}^2\right.\\
&\left.+\frac{11}{5} {t_{1,1}} {t_{2,2}} {t_{3,2}} {t_{3,4}}+\frac{121}{50} {t_{1,1}} {t_{2,2}} {t_{3,3}} {t_{3,4}}^2+\frac{1331 {t_{1,1}} {t_{2,2}} {t_{3,4}}^4}{3000}\right)+e^{12 {t_{\mu_A}}} \left(\frac{{t_{3,1}}^2}{2}+{t_{1,1}}^2 {t_{3,2}}+\frac{{t_{2,2}}^3 {t_{3,2}}}{2}+\frac{{t_{3,1}} {t_{3,3}}^2}{10}\right.\\
&+\frac{53 {t_{3,3}}^4}{600}+\frac{{t_{3,1}} {t_{3,2}} {t_{3,4}}}{5}+\frac{7}{5} {t_{1,1}}^2 {t_{3,3}} {t_{3,4}}+{t_{2,1}} {t_{2,2}} {t_{3,3}} {t_{3,4}}+\frac{41}{30} {t_{2,2}}^3 {t_{3,3}} {t_{3,4}}+\frac{1}{50} {t_{3,2}} {t_{3,3}}^2 {t_{3,4}}+\frac{27 {t_{3,2}}^2 {t_{3,4}}^2}{100}\\
&+\frac{1}{50} {t_{3,1}} {t_{3,3}} {t_{3,4}}^2+\frac{51 {t_{3,3}}^3 {t_{3,4}}^2}{500}+\frac{62 {t_{1,1}}^2 {t_{3,4}}^3}{75}+\frac{3}{10} {t_{2,1}} {t_{2,2}} {t_{3,4}}^3+\frac{53 {t_{2,2}}^3 {t_{3,4}}^3}{100}+\frac{51}{250} {t_{3,2}} {t_{3,3}} {t_{3,4}}^3\\
&\left.+\frac{{t_{3,1}} {t_{3,4}}^4}{3000}+\frac{2557 {t_{3,3}}^2 {t_{3,4}}^4}{30000}+\frac{67 {t_{3,2}} {t_{3,4}}^5}{5000}+\frac{717 {t_{3,3}} {t_{3,4}}^6}{50000}+\frac{5117 {t_{3,4}}^8}{6000000}\right)+e^{13 {t_{\mu_A}}} \left({t_{1,1}} {t_{2,1}} {t_{3,3}}+\frac{13}{6} {t_{1,1}} {t_{2,2}}^2 {t_{3,3}}\right.\\
&\left.+\frac{13}{10} {t_{1,1}} {t_{2,1}} {t_{3,4}}^2+\frac{169}{60} {t_{1,1}} {t_{2,2}}^2 {t_{3,4}}^2\right)+e^{14 {t_{\mu_A}}} \left(\frac{{t_{2,2}} {t_{3,2}}^2}{2}+\frac{{t_{2,1}}^2 {t_{3,4}}}{2}+\frac{7}{2} {t_{1,1}}^2 {t_{2,2}} {t_{3,4}}+\frac{7}{6} {t_{2,1}} {t_{2,2}}^2 {t_{3,4}}\right.\\
&\left.+\frac{49 {t_{2,2}}^4 {t_{3,4}}}{72}+\frac{7}{5} {t_{2,2}} {t_{3,2}} {t_{3,3}} {t_{3,4}}+\frac{49}{50} {t_{2,2}} {t_{3,3}}^2 {t_{3,4}}^2+\frac{49}{150} {t_{2,2}} {t_{3,2}} {t_{3,4}}^3+\frac{343}{750} {t_{2,2}} {t_{3,3}} {t_{3,4}}^4+\frac{2401 {t_{2,2}} {t_{3,4}}^6}{45000}\right)\\
&+e^{15 {t_{\mu_A}}} \left(\frac{{t_{1,1}}^3}{3}+2 {t_{1,1}} {t_{2,1}} {t_{2,2}}+\frac{4 {t_{1,1}} {t_{2,2}}^3}{3}+{t_{1,1}} {t_{3,3}}^2 {t_{3,4}}+{t_{1,1}} {t_{3,2}} {t_{3,4}}^2+{t_{1,1}} {t_{3,3}} {t_{3,4}}^3+\frac{7 {t_{1,1}} {t_{3,4}}^5}{60}\right)\\
&+e^{16 {t_{\mu_A}}} ({t_{2,2}}^2 {t_{3,3}}^2+{t_{2,1}} {t_{3,2}} {t_{3,4}}+\frac{2}{3} {t_{2,2}}^2 {t_{3,2}} {t_{3,4}}+\frac{2}{5} {t_{2,1}} {t_{3,3}} {t_{3,4}}^2+\frac{28}{15} {t_{2,2}}^2 {t_{3,3}} {t_{3,4}}^2+\frac{2 {t_{2,1}} {t_{3,4}}^4}{75}\\
&\left.+\frac{148 {t_{2,2}}^2 {t_{3,4}}^4}{225}\right)+e^{17 {t_{\mu_A}}} \left({t_{1,1}} {t_{2,2}} {t_{3,2}}+\frac{17}{5} {t_{1,1}} {t_{2,2}} {t_{3,3}} {t_{3,4}}+\frac{289}{150} {t_{1,1}} {t_{2,2}} {t_{3,4}}^3\right)+e^{18 {t_{\mu_A}}} \Biggl({t_{1,1}}^2 {t_{3,3}}+{t_{2,1}} {t_{2,2}} {t_{3,3}}\\
&+\frac{{t_{2,2}}^3 {t_{3,3}}}{2}+\frac{{t_{3,2}} {t_{3,3}}^2}{2}+\frac{{t_{3,3}}^3 {t_{3,4}}}{5}+\frac{9 {t_{1,1}}^2 {t_{3,4}}^2}{5}+\frac{3}{10} {t_{2,1}} {t_{2,2}} {t_{3,4}}^2+\frac{33 {t_{2,2}}^3 {t_{3,4}}^2}{20}+\frac{3}{10} {t_{3,2}} {t_{3,3}} {t_{3,4}}^2+\frac{3 {t_{3,3}}^2 {t_{3,4}}^3}{10}\\
&\left.+\frac{9 {t_{3,2}} {t_{3,4}}^4}{200}+\frac{63 {t_{3,3}} {t_{3,4}}^5}{500}+\frac{81 {t_{3,4}}^7}{5000}\right)+e^{19 {t_{\mu_A}}} \left({t_{1,1}} {t_{2,1}} {t_{3,4}}+\frac{19}{6} {t_{1,1}} {t_{2,2}}^2 {t_{3,4}}\right)+e^{20 {t_{\mu_A}}} \left(\frac{{t_{2,1}}^2}{2}+{t_{1,1}}^2 {t_{2,2}}\right.\\
&\left.+\frac{{t_{2,1}} {t_{2,2}}^2}{6}+\frac{37 {t_{2,2}}^4}{72}+\frac{1}{2} {t_{2,2}} {t_{3,3}}^2 {t_{3,4}}+\frac{1}{2} {t_{2,2}} {t_{3,2}} {t_{3,4}}^2+\frac{3}{2} {t_{2,2}} {t_{3,3}} {t_{3,4}}^3+\frac{43 {t_{2,2}} {t_{3,4}}^5}{120}\right)+e^{21 {t_{\mu_A}}} \Biggl({t_{1,1}} {t_{3,2}} {t_{3,4}}\\
&\left.+\frac{7}{5} {t_{1,1}} {t_{3,3}} {t_{3,4}}^2+\frac{49 {t_{1,1}} {t_{3,4}}^4}{150}\right)+e^{22 {t_{\mu_A}}} \left(\frac{{t_{2,2}}^2 {t_{3,2}}}{2}+\frac{11}{5} {t_{2,2}}^2 {t_{3,3}} {t_{3,4}}+\frac{121 {t_{2,2}}^2 {t_{3,4}}^3}{75}\right)+e^{23 {t_{\mu_A}}} \biggl({t_{1,1}} {t_{2,2}} {t_{3,3}}\\
&\left.+\frac{23}{10} {t_{1,1}} {t_{2,2}} {t_{3,4}}^2\right)+e^{24 {t_{\mu_A}}} \left(\frac{{t_{3,2}}^2}{4}+{t_{1,1}}^2 {t_{3,4}}+{t_{2,2}}^3 {t_{3,4}}+\frac{{t_{3,2}} {t_{3,3}} {t_{3,4}}}{5}+\frac{27 {t_{3,3}}^2 {t_{3,4}}^2}{50}+\frac{{t_{3,2}} {t_{3,4}}^3}{75}\right.\\
&\left.+\frac{229 {t_{3,3}} {t_{3,4}}^4}{750}+\frac{977 {t_{3,4}}^6}{11250}\right)+e^{25 {t_{\mu_A}}} {t_{1,1}} {t_{2,2}}^2+e^{26 {t_{\mu_A}}} \left(\frac{{t_{2,2}} {t_{3,3}}^2}{2}+\frac{13}{10} {t_{2,2}} {t_{3,3}} {t_{3,4}}^2+\frac{169 {t_{2,2}} {t_{3,4}}^4}{200}\right)\\
&+e^{27 {t_{\mu_A}}} \left({t_{1,1}} {t_{3,3}} {t_{3,4}}+\frac{3 {t_{1,1}} {t_{3,4}}^3}{10}\right)+\frac{7}{4} e^{28 {t_{\mu_A}}} {t_{2,2}}^2 {t_{3,4}}^2+e^{29 {t_{\mu_A}}} {t_{1,1}} {t_{2,2}} {t_{3,4}}+e^{30 {t_{\mu_A}}} \left(\frac{{t_{1,1}}^2}{2}+\frac{{t_{2,2}}^3}{3}+\frac{{t_{3,3}} {t_{3,4}}^3}{3}\right.\\
&\left.+\frac{{t_{3,4}}^5}{5}\right)+e^{32 {t_{\mu_A}}} \left({t_{2,2}} {t_{3,3}} {t_{3,4}}+\frac{4 {t_{2,2}} {t_{3,4}}^3}{5}\right)+\frac{1}{2} e^{34 {t_{\mu_A}}} {t_{2,2}}^2 {t_{3,4}}+e^{36 {t_{\mu_A}}} \left(\frac{{t_{3,3}}^2}{6}+\frac{{t_{3,3}} {t_{3,4}}^2}{10}+\frac{53 {t_{3,4}}^4}{200}\right)\\
&+\frac{1}{2} e^{38 {t_{\mu_A}}} {t_{2,2}} {t_{3,4}}^2+\frac{1}{4} e^{40 {t_{\mu_A}}} {t_{2,2}}^2+\frac{1}{6} e^{42 {t_{\mu_A}}} {t_{3,4}}^3+\frac{1}{8} e^{48 {t_{\mu_A}}} {t_{3,4}}^2+\frac{e^{60 {t_{\mu_A}}}}{60}.
\end{align*}
}
%%%%%%%%%%%%%%%%%%%%%%%%%%%%%%%%%%%%%%%%%%%%%%%%%%%%%%%%%%%%%%%%%%%%%%%%%%%%%%%
\subsection{Coordinate change for $A=(2,3,5)$}
\
{\tiny
\begin{align*}
s_1({\bf t})=&t_1+e^{2 t_{\mu_A}} \left(-2 t_{2,2} t_{3,2}-\frac{4}{5} t_{2,2} t_{3,3} t_{3,4}-\frac{4}{75} t_{2,2} t_{3,4}^3\right)+e^{3 t_{\mu_A}} \left(-6 t_{1,1} t_{3,3}-\frac{9}{5} t_{1,1} t_{3,4}^2\right)-12 e^{4 t_{\mu_A}} t_{2,1} t_{3,4}\\
&-10 e^{5 t_{\mu_A}} t_{1,1} t_{2,2}+e^{6 t_{\mu_A}} \left(-60 t_{3,1}+12 t_{3,3}^2+18 t_{3,2} t_{3,4}+\frac{108}{5} t_{3,3} t_{3,4}^2+\frac{12 t_{3,4}^4}{5}\right)\\
&+e^{8 t_{\mu_A}} \left(264 t_{2,2} t_{3,3}+\frac{656}{5} t_{2,2} t_{3,4}^2\right)+918 e^{9 t_{\mu_A}} t_{1,1} t_{3,4}+e^{10 t_{\mu_A}} \left(1380 t_{2,1}+90 t_{2,2}^2\right)\\
&+e^{12 t_{\mu_A}} \left(-252 t_{3,2}-\frac{33204}{5} t_{3,3} t_{3,4}-\frac{66768 t_{3,4}^3}{25}\right)-36566 e^{14 t_{\mu_A}} t_{2,2} t_{3,4}-72390 e^{15 t_{\mu_A}} t_{1,1}\\
&+e^{18 t_{\mu_A}} \left(334020 t_{3,3}+689508 t_{3,4}^2\right)+2297860 e^{20 t_{\mu_A}} t_{2,2}-60218748 e^{24 t_{\mu_A}} t_{3,4}+1594559250 e^{30 t_{\mu_A}},\\
s_{1,1}({\bf t})=&t_{1,1}+e^{3 t_{\mu_A}} \left(-8 t_{3,3}-\frac{12 t_{3,4}^2}{5}\right)-24 e^{5 t_{\mu_A}} t_{2,2}+1000 e^{9 t_{\mu_A}} t_{3,4}-59400 e^{15 t_{\mu_A}},\\
s_{2,1}({\bf t})=&t_{2,1}+\frac{t_{2,2}^2}{6}+e^{2 t_{\mu_A}} \left(-3 t_{3,2}-\frac{6}{5} t_{3,3} t_{3,4}-\frac{2 t_{3,4}^3}{25}\right)-15 e^{4 t_{\mu_A}} t_{2,2} t_{3,4}-24 e^{5 t_{\mu_A}} t_{1,1}\\
&+e^{8 t_{\mu_A}} \left(225 t_{3,3}+\frac{615 t_{3,4}^2}{2}\right)+1725 e^{10 t_{\mu_A}} t_{2,2}-56640 e^{14 t_{\mu_A}} t_{3,4}+2508000 e^{20 t_{\mu_A}},\\
s_{2,2}({\bf t})=&t_{2,2}-27 e^{4 t_{\mu_A}} t_{3,4}+2376 e^{10 t_{\mu_A}},\\
s_{3,1}({\bf t})=&t_{3,1}+\frac{t_{3,3}^2}{10}+\frac{1}{5} t_{3,2} t_{3,4}+\frac{1}{50} t_{3,3} t_{3,4}^2+\frac{t_{3,4}^4}{3000}+e^{2 t_{\mu_A}} \left(-3 t_{2,2} t_{3,3}-\frac{9}{10} t_{2,2} t_{3,4}^2\right)-8 e^{3 t_{\mu_A}} t_{1,1} t_{3,4}\\
&+e^{4 t_{\mu_A}} \left(-15 t_{2,1}+\frac{5 t_{2,2}^2}{2}\right)+e^{6 t_{\mu_A}} \left(-80 t_{3,2}+13 t_{3,3} t_{3,4}+\frac{341 t_{3,4}^3}{30}\right)+522 e^{8 t_{\mu_A}} t_{2,2} t_{3,4}\\
&+1360 e^{9 t_{\mu_A}} t_{1,1}+e^{12 t_{\mu_A}} \left(1872 t_{3,3}-\frac{23967 t_{3,4}^2}{5}\right)-56115 e^{14 t_{\mu_A}} t_{2,2}+387405 e^{18 t_{\mu_A}} t_{3,4}+3926880 e^{24 t_{\mu_A}},\\
s_{3,2}({\bf t})=&-4 e^{2 t_{\mu_A}} t_{2,2} t_{3,4},t_{3,2}+\frac{2}{5} t_{3,3} t_{3,4}+\frac{2 t_{3,4}^3}{75}-10 e^{3 t_{\mu_A}} t_{1,1}+e^{6 t_{\mu_A}} \left(-120 t_{3,3}-8 t_{3,4}^2\right)\\
&+865 e^{8 t_{\mu_A}} t_{2,2}+6496 e^{12 t_{\mu_A}} t_{3,4}-971000 e^{18 t_{\mu_A}},\\
s_{3,3}({\bf t})=&t_{3,3}+\frac{3 t_{3,4}^2}{10}-5 e^{2 t_{\mu_A}} t_{2,2}-165 e^{6 t_{\mu_A}} t_{3,4}+23745 e^{12 t_{\mu_A}},\\
s_{3,4}({\bf t})=&t_{3,4}-250 e^{6 t_{\mu_A}}.
\end{align*}
}
%%%%%%%%%%%%%%%%%%%%%%%%%%%%%%%%%%%%%%%%%%%%%%%%%%%%%%%%%%%%%%%%%%%%%%%%%%%%%%%
\subsection{$\phi$ for $A=(2,3,5)$}
\
{\tiny
\[
\phi ^a_{ij}=0,\ a=1,2,3,i=1,j=1,\dots,9.
\]
{\tiny
\[
\phi^1_{22}=\frac{x_1}{2}+e^{5 t_{\mu_A}} (40 t_{2,2}+60 x_2)+e^{9 t_{\mu_A}} (-1464 t_{3,4}-480 x_3)+160320 e^{15 t_{\mu_A}},
\]
\[
\phi^2_{22}=e^{4 t_{\mu_A}} (16 t_{3,4}+20 x_3)-3680 e^{10 t_{\mu_A}},\
\phi^3_{22}=168 e^{6 t_{\mu_A}}.
\]
\[
\phi^1_{23}=e^{4 t_{\mu_A}} (8 t_{3,4}+10 x_3)-1720 e^{10 t_{\mu_A}},\
\phi^2_{23}=35 e^{5 t_{\mu_A}},\
\phi^3_{23}=-e^{t_{\mu_A}}.
\]

\begin{align*}
&\phi^1_{24}=e^{2 t_{\mu_A}} \left(-2 t_{3,2}-\frac{4}{5} t_{3,3} t_{3,4}-\frac{4 t_{3,4}^3}{75}-3 t_{3,3} x_3-\frac{9}{10} t_{3,4}^2 x_3-4 t_{3,4} x_3^2-5 x_3^3\right)+e^{4 t_{\mu_A}} \left(\frac{8}{3} t_{2,2} t_{3,4}-4 t_{3,4} x_2\right.\\
&\left.+\frac{25}{3} t_{2,2} x_3-5 x_2 x_3\right)+e^{5 t_{\mu_A}} (15 t_{1,1}-12 x_1)+e^{8 t_{\mu_A}} \left(208 t_{3,3}+\frac{512 t_{3,4}^2}{5}+519 t_{3,4} x_3+880 x_3^2\right)\\
&+e^{10 t_{\mu_A}} \left(-\frac{1900 t_{2,2}}{3}+605 x_2\right)+e^{14 t_{\mu_A}} (-25642 t_{3,4}-56985 x_3)+1658800 e^{20 t_{\mu_A}},
\end{align*}
\[
\phi^2_{24}=e^{3 t_{\mu_A}} \left(-3 t_{3,3}-\frac{9 t_{3,4}^2}{10}-4 t_{3,4} x_3-5 x_3^2\right)+e^{5 t_{\mu_A}} \left(\frac{35 t_{2,2}}{3}-8 x_2\right)+e^{9 t_{\mu_A}} (324 t_{3,4}+655 x_3)-25215 e^{15 t_{\mu_A}},
\]
\[
\phi^3_{24}=e^{t_{\mu_A}} \left(-\frac{t_{2,2}}{3}-x_2\right)+e^{5 t_{\mu_A}} (7 t_{3,4}+10 x_3)-1091 e^{11 t_{\mu_A}}.\quad \phi^1_{25}=24 e^{6 t_{\mu_A}},\
\phi^2_{25}=-e^{t_{\mu_A}},\
\phi^3_{25}=0.
\]
\[
\phi^1_{26}=e^{2 t_{\mu_A}} (-2 t_{2,2}-3 x_2)+e^{6 t_{\mu_A}} \left(\frac{174 t_{3,4}}{5}-6 x_3\right)-2592 e^{12 t_{\mu_A}},
\]
\[
\phi^2_{26}=e^{t_{\mu_A}} \left(-\frac{t_{3,4}}{5}-x_3\right)+68 e^{7 t_{\mu_A}},\
\phi^3_{26}=-3 e^{3 t_{\mu_A}}.
\]
\begin{align*}
\phi^1_{27}&=e^{2 t_{\mu_A}} \left(-\frac{4}{5} t_{2,2} t_{3,4}-\frac{6}{5} t_{3,4} x_2-2 t_{2,2} x_3-3 x_2 x_3\right)+e^{3 t_{\mu_A}} (-3 t_{1,1}-4 x_1)+e^{6 t_{\mu_A}} \left(\frac{144 t_{3,3}}{5}+\frac{492 t_{3,4}^2}{25}\right.\\
&\left.+\frac{258}{5} t_{3,4} x_3+24 x_3^2\right)+e^{8 t_{\mu_A}} (176 t_{2,2}+264 x_2)+e^{12 t_{\mu_A}} \left(-\frac{29544 t_{3,4}}{5}-5382 x_3\right)+335304 e^{18 t_{\mu_A}},\\
\end{align*}
\[
\phi^2_{27}=e^{t_{\mu_A}} \left(-\frac{t_{3,3}}{5}-\frac{t_{3,4}^2}{50}-\frac{2}{5} t_{3,4} x_3-x_3^2\right)-e^{3 t_{\mu_A}} t_{2,2}+e^{7 t_{\mu_A}} \left(\frac{131 t_{3,4}}{5}+88 x_3\right)-1756 e^{13 t_{\mu_A}},
\]
\[
\phi^3_{27}=e^{3 t_{\mu_A}} \left(-\frac{6 t_{3,4}}{5}-2 x_3\right)+111 e^{9 t_{\mu_A}}.
\]

\begin{align*}
&\phi^1_{28}=e^{2 t_{\mu_A}} \left(-\frac{4}{5} t_{2,2} t_{3,3}-\frac{4}{25} t_{2,2} t_{3,4}^2-\frac{6}{5} t_{3,3} x_2-\frac{6}{25} t_{3,4}^2 x_2-\frac{6}{5} t_{2,2} t_{3,4} x_3-\frac{9}{5} t_{3,4} x_2 x_3-2 t_{2,2} x_3^2-3 x_2 x_3^2\right)\\
&+e^{3 t_{\mu_A}} \left(-\frac{9}{5} t_{1,1} t_{3,4}-\frac{12}{5} t_{3,4} x_1-3 t_{1,1} x_3-4 x_1 x_3\right)+e^{6 t_{\mu_A}} \left(\frac{54 t_{3,2}}{5}+\frac{744}{25} t_{3,3} t_{3,4}+\frac{864 t_{3,4}^3}{125}+\frac{153}{5} t_{3,3} x_3\right.\\
&+\frac{2091}{50} t_{3,4}^2 x_3+\frac{402}{5} t_{3,4} x_3^2+39 x_3^3\biggr)+e^{8 t_{\mu_A}} \left(\frac{568}{5} t_{2,2} t_{3,4}+\frac{852}{5} t_{3,4} x_2+267 t_{2,2} x_3+459 x_2 x_3\right)+e^{9 t_{\mu_A}} (207 t_{1,1}+500 x_1)\\
&+e^{12 t_{\mu_A}} \left(-\frac{21744 t_{3,3}}{5}-\frac{105504 t_{3,4}^2}{25}-\frac{73461}{5} t_{3,4} x_3-11712 x_3^2\right)+e^{14 t_{\mu_A}} (-12738 t_{2,2}-23040 x_2)\\
&+e^{18 t_{\mu_A}} \left(\frac{3865152 t_{3,4}}{5}+1136145 x_3\right)-41422320 e^{24 t_{\mu_A}},\\
&\phi^2_{28}=e^{t_{\mu_A}} \left(-\frac{t_{3,2}}{5}-\frac{1}{25} t_{3,3} t_{3,4}-\frac{t_{3,4}^3}{750}-\frac{2}{5} t_{3,3} x_3-\frac{2}{25} t_{3,4}^2 x_3-\frac{3}{5} t_{3,4} x_3^2-x_3^3\right)-\frac{3}{5} e^{3 t_{\mu_A}} t_{2,2} t_{3,4}+e^{4 t_{\mu_A}} (-2 t_{1,1}-8 x_1)\\
&+e^{7 t_{\mu_A}} \left(\frac{401 t_{3,3}}{5}+\frac{1037 t_{3,4}^2}{50}+\frac{644}{5} t_{3,4} x_3+243 x_3^2\right)+e^{9 t_{\mu_A}} (108 t_{2,2}+216 x_2)+e^{13 t_{\mu_A}} \left(-\frac{46348 t_{3,4}}{5}-19920 x_3\right)\\
&+762485 e^{19 t_{\mu_A}},\\
&\phi^3_{28}=e^{3 t_{\mu_A}} \left(-\frac{6 t_{3,3}}{5}-\frac{6 t_{3,4}^2}{25}-\frac{6}{5} t_{3,4} x_3-2 x_3^2\right)+e^{5 t_{\mu_A}} (-t_{2,2}+3 x_2)+e^{9 t_{\mu_A}} \left(\frac{213 t_{3,4}}{5}+202 x_3\right)-3456 e^{15 t_{\mu_A}}.\\
&\phi^1_{29}=e^{2 t_{\mu_A}} \left(-4 t_{2,2} t_{3,2}-\frac{8}{5} t_{2,2} t_{3,3} t_{3,4}-\frac{8}{75} t_{2,2} t_{3,4}^3-6 t_{3,2} x_2-\frac{12}{5} t_{3,3} t_{3,4} x_2-\frac{4}{25} t_{3,4}^3 x_2-6 t_{2,2} t_{3,3} x_3\right.\\
&-\frac{9}{5} t_{2,2} t_{3,4}^2 x_3-9 t_{3,3} x_2 x_3-\frac{27}{10} t_{3,4}^2 x_2 x_3-8 t_{2,2} t_{3,4} x_3^2-12 t_{3,4} x_2 x_3^2-10 t_{2,2} x_3^3-15 x_2 x_3^3\biggr)+e^{3 t_{\mu_A}} \Biggl(-9 t_{1,1} t_{3,3}\\
&-\frac{27}{10} t_{1,1} t_{3,4}^2-12 t_{3,3} x_1-\frac{18}{5} t_{3,4}^2 x_1-12 t_{1,1} t_{3,4} x_3-16 t_{3,4} x_1 x_3-15 t_{1,1} x_3^2-20 x_1 x_3^2\biggr)+e^{4 t_{\mu_A}} \left(10 t_{2,2}^2 x_3+15 t_{2,2} x_2 x_3\right)\\
&+e^{5 t_{\mu_A}} (15 t_{1,1} t_{2,2}-60 t_{2,2} x_1-96 x_1 x_2)+e^{6 t_{\mu_A}} \left(72 t_{3,3}^2+84 t_{3,2} t_{3,4}+\frac{384}{5} t_{3,3} t_{3,4}^2+\frac{218 t_{3,4}^4}{25}-60 t_{3,2} x_3\right.\\
&+270 t_{3,3} t_{3,4} x_3+\frac{433}{5} t_{3,4}^3 x_3+120 t_{3,3} x_3^2+300 t_{3,4}^2 x_3^2+330 t_{3,4} x_3^3\biggr)+e^{8 t_{\mu_A}} \Biggl(1040 t_{2,2} t_{3,3}+344 t_{2,2} t_{3,4}^2+1632 t_{3,3} x_2\\
\end{align*}

\begin{align*}
&+\frac{2688}{5} t_{3,4}^2 x_2+1800 t_{2,2} t_{3,4} x_3+3051 t_{3,4} x_2 x_3+3200 t_{2,2} x_3^2+4620 x_2 x_3^2\biggr)+e^{9 t_{\mu_A}} (1275 t_{1,1} t_{3,4}+4500 t_{3,4} x_1\\
&+2700 t_{1,1} x_3+3440 x_1 x_3)+e^{10 t_{\mu_A}} \left(-60 t_{2,2}^2+2610 t_{2,2} x_2+4320 x_2^2\right)+e^{12 t_{\mu_A}} \left(-6144 t_{3,2}-\frac{236208}{5} t_{3,3} t_{3,4}\right.\\
&\left.-\frac{349576 t_{3,4}^3}{25}-33936 t_{3,3} x_3-\frac{396234}{5} t_{3,4}^2 x_3-130848 t_{3,4} x_3^2-32160 x_3^3\right)+e^{14 t_{\mu_A}} (-153740 t_{2,2} t_{3,4}-303600 t_{3,4} x_2\\
&-288330 t_{2,2} x_3-447885 x_2 x_3)+e^{15 t_{\mu_A}} (-135885 t_{1,1}-445500 x_1)+e^{18 t_{\mu_A}} \left(3793104 t_{3,3}+\frac{30387216 t_{3,4}^2}{5}\right.\\
&+16278702 t_{3,4} x_3+9527640 x_3^2\Biggr)+e^{20 t_{\mu_A}} (13301840 t_{2,2}+25757760 x_2)+e^{24 t_{\mu_A}} (-828787440 t_{3,4}-884891040 x_3)\\
&+34686244200 e^{30 t_{\mu_A}},\\
&\phi^2_{29}=-x_1 x_2+e^{3 t_{\mu_A}} \left(-6 t_{2,2} t_{3,3}-\frac{9}{5} t_{2,2} t_{3,4}^2-8 t_{2,2} t_{3,4} x_3-10 t_{2,2} x_3^2\right)+e^{4 t_{\mu_A}} (-16 t_{1,1} t_{3,4}-32 t_{3,4} x_1-20 t_{1,1} x_3\\
&-40 x_1 x_3)+e^{5 t_{\mu_A}} \left(20 t_{2,1}+\frac{10 t_{2,2}^2}{3}\right)+e^{7 t_{\mu_A}} \left(96 t_{3,2}+\frac{1462}{5} t_{3,3} t_{3,4}+\frac{1969 t_{3,4}^3}{25}+304 t_{3,3} x_3+\frac{2256}{5} t_{3,4}^2 x_3+882 t_{3,4} x_3^2\right.\\
&+540 x_3^3\Biggr)+e^{9 t_{\mu_A}} (984 t_{2,2} t_{3,4}+864 t_{3,4} x_2+1640 t_{2,2} x_3+1440 x_2 x_3)+e^{10 t_{\mu_A}} (1400 t_{1,1}+6880 x_1)\\
&+e^{13 t_{\mu_A}} \left(-48136 t_{3,3}-\frac{263324 t_{3,4}^2}{5}-174608 t_{3,4} x_3-157260 x_3^2\right)+e^{15 t_{\mu_A}} (-111570 t_{2,2}-190080 x_2)\\
&+e^{19 t_{\mu_A}} (10377650 t_{3,4}+15134160 x_3)-584059200 e^{25 t_{\mu_A}},\\
&\phi^3_{29}=e^{t_{\mu_A}} \left(-t_{2,1}-\frac{t_{2,2}^2}{6}-2 t_{2,2} x_2-3 x_2^2\right)+e^{3 t_{\mu_A}} \left(-3 t_{3,2}-\frac{6}{5} t_{3,3} t_{3,4}-\frac{2 t_{3,4}^3}{25}\right)+e^{5 t_{\mu_A}} (11 t_{2,2} t_{3,4}+66 t_{3,4} x_2\\
&+20 t_{2,2} x_3)+e^{6 t_{\mu_A}} (48 t_{1,1}-192 x_1)+e^{9 t_{\mu_A}} \left(1431 t_{3,3}+\frac{933 t_{3,4}^2}{10}+1980 t_{3,4} x_3+3000 x_3^2\right)+e^{11 t_{\mu_A}} (-571 t_{2,2}+768 x_2)\\
&+e^{15 t_{\mu_A}} (-161184 t_{3,4}-377880 x_3)+15859680 e^{21 t_{\mu_A}}.
\end{align*}

\[
\phi^1_{33}=30 e^{5 t_{\mu_A}},\
\phi^2_{33}=0,\
\phi^3_{33}=0.
\]
\[
\phi^1_{34}=e^{3 t_{\mu_A}} \left(3 t_{3,3}+\frac{9 t_{3,4}^2}{10}+4 t_{3,4} x_3+5 x_3^2\right)+e^{5 t_{\mu_A}} (20 t_{2,2}+45 x_2)+e^{9 t_{\mu_A}} (-1287 t_{3,4}-865 x_3)+129510 e^{15 t_{\mu_A}},
\]
\[
\phi^2_{34}=\frac{x_2}{3}+e^{4 t_{\mu_A}} (8 t_{3,4}+10 x_3)-2090 e^{10 t_{\mu_A}},\
\phi^3_{34}=114 e^{6 t_{\mu_A}}.\quad
\phi^1_{35}=-e^{t_{\mu_A}},\
\phi^2_{35}=0,\
\phi^3_{35}=0.
\]
\[
\phi^1_{36}=e^{t_{\mu_A}} \left(-\frac{t_{3,4}}{5}-x_3\right)+35 e^{7 t_{\mu_A}},\
\phi^2_{36}=-2 e^{2 t_{\mu_A}},\
\phi^3_{36}=0.
\]
\[
\phi^1_{37}=e^{t_{\mu_A}} \left(-\frac{t_{3,3}}{5}-\frac{t_{3,4}^2}{50}-\frac{2}{5} t_{3,4} x_3-x_3^2\right)+e^{3 t_{\mu_A}} (-t_{2,2}-6 x_2)+e^{7 t_{\mu_A}} (77 t_{3,4}+120 x_3)-8001 e^{13 t_{\mu_A}},
\]
\[
\phi^2_{37}=e^{2 t_{\mu_A}} \left(-\frac{4 t_{3,4}}{5}-2 x_3\right)+120 e^{8 t_{\mu_A}},\
\phi^3_{37}=-4 e^{4 t_{\mu_A}}.
\]
\begin{align*}
&\phi^1_{38}=e^{t_{\mu_A}} \left(-\frac{t_{3,2}}{5}-\frac{1}{25} t_{3,3} t_{3,4}-\frac{t_{3,4}^3}{750}-\frac{2}{5} t_{3,3} x_3-\frac{2}{25} t_{3,4}^2 x_3-\frac{3}{5} t_{3,4} x_3^2-x_3^3\right)+e^{3 t_{\mu_A}} \left(-\frac{3}{5} t_{2,2} t_{3,4}-\frac{18}{5} t_{3,4} x_2\right.\\
&-6 x_2 x_3\Biggr)+2 e^{4 t_{\mu_A}} t_{1,1}+e^{7 t_{\mu_A}} \left(17 t_{3,3}+\frac{83 t_{3,4}^2}{2}+124 t_{3,4} x_3+165 x_3^2\right)+e^{9 t_{\mu_A}} (57 t_{2,2}+585 x_2)\\
&+e^{13 t_{\mu_A}} \left(-\frac{53553 t_{3,4}}{5}-16000 x_3\right)+687555 e^{19 t_{\mu_A}},
\end{align*}
\[
\phi^2_{38}=e^{2 t_{\mu_A}} \left(-\frac{4 t_{3,3}}{5}-\frac{4 t_{3,4}^2}{25}-\frac{6}{5} t_{3,4} x_3-2 x_3^2\right)-9 e^{4 t_{\mu_A}} x_2+e^{8 t_{\mu_A}} (176 t_{3,4}+330 x_3)-17000 e^{14 t_{\mu_A}},
\]
\[
\phi^3_{38}=e^{4 t_{\mu_A}} \left(-\frac{12 t_{3,4}}{5}-3 x_3\right)+170 e^{10 t_{\mu_A}}.
\]
\begin{align*}
&\phi^1_{39}=-x_1 x_2+e^{3 t_{\mu_A}} \left(-6 t_{2,2} t_{3,3}-\frac{9}{5} t_{2,2} t_{3,4}^2-18 t_{3,3} x_2-\frac{27}{5} t_{3,4}^2 x_2-8 t_{2,2} t_{3,4} x_3-24 t_{3,4} x_2 x_3-10 t_{2,2} x_3^2-30 x_2 x_3^2\right)\\
&+15 e^{4 t_{\mu_A}} x_1 x_3+e^{5 t_{\mu_A}} \left(-90 t_{2,2} x_2-180 x_2^2\right)+e^{7 t_{\mu_A}} \left(30 t_{3,2}+258 t_{3,3} t_{3,4}+\frac{373 t_{3,4}^3}{5}+150 t_{3,3} x_3+309 t_{3,4}^2 x_3\right.\\
&+450 t_{3,4} x_3^2+150 x_3^3\biggr)+e^{9 t_{\mu_A}} (1800 t_{2,2} t_{3,4}+7290 t_{3,4} x_2+2450 t_{2,2} x_3+6900 x_2 x_3)-300 e^{10 t_{\mu_A}} t_{1,1}\\
&+e^{13 t_{\mu_A}} \left(-27876 t_{3,3}-\frac{354414 t_{3,4}^2}{5}-136518 t_{3,4} x_3-62460 x_3^2\right)+e^{15 t_{\mu_A}} (-179100 t_{2,2}-730890 x_2)\\
&+e^{19 t_{\mu_A}} (12768810 t_{3,4}+10956150 x_3)-648482400 e^{25 t_{\mu_A}},\\
&\phi^2_{39}=e^{2 t_{\mu_A}} \left(-4 t_{3,2}-\frac{8}{5} t_{3,3} t_{3,4}-\frac{8 t_{3,4}^3}{75}-6 t_{3,3} x_3-\frac{9}{5} t_{3,4}^2 x_3-8 t_{3,4} x_3^2-10 x_3^3\right)+e^{4 t_{\mu_A}} (-36 t_{3,4} x_2+10 t_{2,2} x_3\\
&-60 x_2 x_3)+e^{5 t_{\mu_A}} (55 t_{1,1}+30 x_1)+e^{8 t_{\mu_A}} \left(300 t_{3,3}+506 t_{3,4}^2+1890 t_{3,4} x_3+2300 x_3^2\right)+e^{10 t_{\mu_A}} (-720 t_{2,2}+7920 x_2)\\
&+e^{14 t_{\mu_A}} (-149020 t_{3,4}-267390 x_3)+11540300 e^{20 t_{\mu_A}},\\
&\phi^3_{39}=e^{t_{\mu_A}} (-t_{1,1}-2 x_1)+e^{4 t_{\mu_A}} \left(-4 t_{3,3}-\frac{6 t_{3,4}^2}{5}\right)+e^{6 t_{\mu_A}} (24 t_{2,2}-270 x_2)+e^{10 t_{\mu_A}} (2960 t_{3,4}+4500 x_3)\\
&-368820 e^{16 t_{\mu_A}}.
\end{align*}

\begin{align*}
&\phi^1_{44}=-\frac{1}{3} x_1 x_2+e^{3 t_{\mu_A}} \left(-2 t_{2,2} t_{3,3}-\frac{3}{5} t_{2,2} t_{3,4}^2-\frac{8}{3} t_{2,2} t_{3,4} x_3-\frac{10}{3} t_{2,2} x_3^2\right)-5 e^{4 t_{\mu_A}} x_1 x_3+\frac{10}{3} e^{5 t_{\mu_A}} t_{2,2}^2\\
&+e^{7 t_{\mu_A}} \left(36 t_{3,2}+\frac{162}{5} t_{3,3} t_{3,4}+\frac{159 t_{3,4}^3}{25}+54 t_{3,3} x_3+\frac{201}{5} t_{3,4}^2 x_3+102 t_{3,4} x_3^2+90 x_3^3\right)+e^{9 t_{\mu_A}} \Biggl(162 t_{2,2} t_{3,4}\\
&\left.-18 t_{3,4} x_2+\frac{1100}{3} t_{2,2} x_3+90 x_2 x_3\right)+e^{10 t_{\mu_A}} (-270 t_{1,1}+288 x_1)+e^{13 t_{\mu_A}} \left(-4194 t_{3,3}-\frac{15291 t_{3,4}^2}{5}-14592 t_{3,4} x_3\right.\\
&\left.-16590 x_3^2\right)+e^{15 t_{\mu_A}} (-13320 t_{2,2}-5220 x_2)+e^{19 t_{\mu_A}} (578790 t_{3,4}+1180860 x_3)-32746050 e^{25 t_{\mu_A}},\\
&\phi^2_{44}=\frac{x_2^2}{3}+e^{2 t_{\mu_A}} \left(-\frac{4 t_{3,2}}{3}-\frac{8}{15} t_{3,3} t_{3,4}-\frac{8 t_{3,4}^3}{225}-2 t_{3,3} x_3-\frac{3}{5} t_{3,4}^2 x_3-\frac{8}{3} t_{3,4} x_3^2-\frac{10 x_3^3}{3}\right)+e^{4 t_{\mu_A}} \Biggl(-4 t_{3,4} x_2\\
&\left.+\frac{10}{3} t_{2,2} x_3\right)+e^{5 t_{\mu_A}} \left(\frac{25 t_{1,1}}{3}-10 x_1\right)+e^{8 t_{\mu_A}} \left(164 t_{3,3}+\frac{898 t_{3,4}^2}{15}+342 t_{3,4} x_3+\frac{1820 x_3^2}{3}\right)\\
&+e^{10 t_{\mu_A}} (-60 t_{2,2}+622 x_2)+e^{14 t_{\mu_A}} \left(-\frac{53368 t_{3,4}}{3}-39150 x_3\right)+\frac{3779600}{3} e^{20 t_{\mu_A}},
\end{align*}

\begin{align*}
&\phi^3_{44}=e^{t_{\mu_A}} \left(-\frac{t_{1,1}}{3}-\frac{2 x_1}{3}\right)+e^{4 t_{\mu_A}} \left(\frac{8 t_{3,3}}{3}+\frac{4 t_{3,4}^2}{5}+4 t_{3,4} x_3+5 x_3^2\right)\\
&+e^{6 t_{\mu_A}} (4 t_{2,2}+18 x_2)+e^{10 t_{\mu_A}} \left(-\frac{1648 t_{3,4}}{3}-730 x_3\right)+44640 e^{16 t_{\mu_A}}.\\
&\phi^1_{45}=e^{t_{\mu_A}} \left(-\frac{t_{2,2}}{3}-x_2\right)+15 e^{5 t_{\mu_A}} t_{3,4}-1725 e^{11 t_{\mu_A}},\
\phi^2_{45}=24 e^{6 t_{\mu_A}},\
\phi^3_{45}=-2 e^{2 t_{\mu_A}}\\
&\phi^1_{46}=e^{t_{\mu_A}} \left(-\frac{1}{15} t_{2,2} t_{3,4}-\frac{1}{5} t_{3,4} x_2-\frac{1}{3} t_{2,2} x_3-x_2 x_3\right)-e^{2 t_{\mu_A}} t_{1,1}+e^{5 t_{\mu_A}} \left(-t_{3,3}+\frac{27 t_{3,4}^2}{10}+3 t_{3,4} x_3-15 x_3^2\right)\\
&+e^{7 t_{\mu_A}} \left(\frac{161 t_{2,2}}{3}+62 x_2\right)+e^{11 t_{\mu_A}} (-655 t_{3,4}+870 x_3)+543 e^{17 t_{\mu_A}},\\
\end{align*}

\begin{align*}
&\phi^2_{46}=e^{2 t_{\mu_A}} \left(-\frac{2 t_{2,2}}{3}-x_2\right)+e^{6 t_{\mu_A}} \left(\frac{54 t_{3,4}}{5}-6 x_3\right)-180 e^{12 t_{\mu_A}},\
\phi^3_{46}=e^{2 t_{\mu_A}} \left(-\frac{2 t_{3,4}}{5}-x_3\right)+52 e^{8 t_{\mu_A}}.\\
&\phi^1_{47}=e^{t_{\mu_A}} \left(-\frac{1}{15} t_{2,2} t_{3,3}-\frac{1}{150} t_{2,2} t_{3,4}^2-\frac{1}{5} t_{3,3} x_2-\frac{1}{50} t_{3,4}^2 x_2-\frac{2}{15} t_{2,2} t_{3,4} x_3-\frac{2}{5} t_{3,4} x_2 x_3-\frac{1}{3} t_{2,2} x_3^2-x_2 x_3^2\right)\\
&+e^{2 t_{\mu_A}} \left(-\frac{2}{5} t_{1,1} t_{3,4}-t_{1,1} x_3-x_1 x_3\right)+e^{3 t_{\mu_A}} \left(-\frac{1}{3} t_{2,2}^2+t_{2,2} x_2\right)+e^{5 t_{\mu_A}} \left(4 t_{3,2}+\frac{21}{5} t_{3,3} t_{3,4}+\frac{43 t_{3,4}^3}{150}+14 t_{3,3} x_3\right.\\
&\left.+\frac{27}{5} t_{3,4}^2 x_3+17 t_{3,4} x_3^2+10 x_3^3\right)+e^{7 t_{\mu_A}} \left(\frac{142}{15} t_{2,2} t_{3,4}+\frac{109}{5} t_{3,4} x_2+42 t_{2,2} x_3+132 x_2 x_3\right)+e^{8 t_{\mu_A}} (92 t_{1,1}+96 x_1)\\
&+e^{11 t_{\mu_A}} \left(-1206 t_{3,3}-\frac{1929 t_{3,4}^2}{5}-2758 t_{3,4} x_3-2600 x_3^2\right)+e^{13 t_{\mu_A}} (-1770 t_{2,2}-5352 x_2)\\
&+e^{17 t_{\mu_A}} \left(\frac{601871 t_{3,4}}{5}+233253 x_3\right)-7840445 e^{23 t_{\mu_A}},\\
&\phi^2_{47}=e^{2 t_{\mu_A}} \left(-\frac{4}{15} t_{2,2} t_{3,4}-\frac{2}{5} t_{3,4} x_2-\frac{2}{3} t_{2,2} x_3\right)+e^{3 t_{\mu_A}} (-3 t_{1,1}-2 x_1)\\
&+e^{6 t_{\mu_A}} \left(\frac{84 t_{3,3}}{5}+\frac{162 t_{3,4}^2}{25}+\frac{58}{5} t_{3,4} x_3+4 x_3^2\right)+e^{8 t_{\mu_A}} (76 t_{2,2}+75 x_2)+e^{12 t_{\mu_A}} (-2028 t_{3,4}-1160 x_3)+121914 e^{18 t_{\mu_A}},
\end{align*}

\begin{align*}
&\phi^3_{47}=e^{2 t_{\mu_A}} \left(-\frac{2 t_{3,3}}{5}-\frac{t_{3,4}^2}{25}-\frac{2}{5} t_{3,4} x_3-x_3^2\right)+e^{4 t_{\mu_A}} \left(\frac{2 t_{2,2}}{3}+2 x_2\right)+e^{8 t_{\mu_A}} \left(\frac{14 t_{3,4}}{5}+93 x_3\right)-1008 e^{14 t_{\mu_A}}.\\
&\phi^1_{48}=e^{t_{\mu_A}} \left(-\frac{1}{15} t_{2,2} t_{3,2}-\frac{1}{75} t_{2,2} t_{3,3} t_{3,4}-\frac{t_{2,2} t_{3,4}^3}{2250}-\frac{1}{5} t_{3,2} x_2-\frac{1}{25} t_{3,3} t_{3,4} x_2-\frac{1}{750} t_{3,4}^3 x_2-\frac{2}{15} t_{2,2} t_{3,3} x_3\right.\\
&\left.-\frac{2}{75} t_{2,2} t_{3,4}^2 x_3-\frac{2}{5} t_{3,3} x_2 x_3-\frac{2}{25} t_{3,4}^2 x_2 x_3-\frac{1}{5} t_{2,2} t_{3,4} x_3^2-\frac{3}{5} t_{3,4} x_2 x_3^2-\frac{1}{3} t_{2,2} x_3^3-x_2 x_3^3\right)+e^{2 t_{\mu_A}} \Biggl(-\frac{2}{5} t_{1,1} t_{3,3}\\
&\left.-\frac{2}{25} t_{1,1} t_{3,4}^2-\frac{3}{5} t_{1,1} t_{3,4} x_3-\frac{3}{5} t_{3,4} x_1 x_3-t_{1,1} x_3^2-x_1 x_3^2\right)+e^{3 t_{\mu_A}} \left(-\frac{1}{5} t_{2,2}^2 t_{3,4}+\frac{3}{5} t_{2,2} t_{3,4} x_2+2 t_{2,2} x_2 x_3\right)\\
&+e^{4 t_{\mu_A}} \left(-\frac{4}{3} t_{1,1} t_{2,2}-t_{1,1} x_2+3 x_1 x_2\right)+e^{5 t_{\mu_A}} \left(-\frac{2}{5} t_{3,3}^2+\frac{27}{5} t_{3,2} t_{3,4}+\frac{34}{25} t_{3,3} t_{3,4}^2+\frac{3 t_{3,4}^4}{50}+\frac{48}{5} t_{3,3} t_{3,4} x_3\right.\\
&\left.+\frac{69}{25} t_{3,4}^3 x_3+2 t_{3,3} x_3^2+15 t_{3,4}^2 x_3^2+21 t_{3,4} x_3^3\right)+e^{7 t_{\mu_A}} \left(\frac{442}{15} t_{2,2} t_{3,3}-\frac{38}{75} t_{2,2} t_{3,4}^2+\frac{124}{5} t_{3,3} x_2+\frac{79}{25} t_{3,4}^2 x_2+\frac{458}{15} t_{2,2} t_{3,4} x_3\right.\\
&\left.+\frac{386}{5} t_{3,4} x_2 x_3+112 t_{2,2} x_3^2+192 x_2 x_3^2\right)+e^{8 t_{\mu_A}} \left(\frac{476}{5} t_{1,1} t_{3,4}+\frac{288}{5} t_{3,4} x_1+192 t_{1,1} x_3+246 x_1 x_3\right)\\
&+e^{9 t_{\mu_A}} \left(-23 t_{2,2}^2-201 t_{2,2} x_2\right)+e^{11 t_{\mu_A}} \left(-689 t_{3,2}-\frac{4921}{5} t_{3,3} t_{3,4}-\frac{32063 t_{3,4}^3}{150}-1704 t_{3,3} x_3-\frac{12828}{5} t_{3,4}^2 x_3\right.\\
&-6239 t_{3,4} x_3^2-2585 x_3^3\Biggr)+e^{13 t_{\mu_A}} \left(-124 t_{2,2} t_{3,4}-\frac{11796}{5} t_{3,4} x_2-\frac{23239}{3} t_{2,2} x_3-13531 x_2 x_3\right)\\
&+e^{14 t_{\mu_A}} (-10547 t_{1,1}-12000 x_1)+e^{17 t_{\mu_A}} \left(\frac{508921 t_{3,3}}{5}+\frac{5545087 t_{3,4}^2}{50}+\frac{2898999}{5} t_{3,4} x_3+489433 x_3^2\right)\\
&+e^{19 t_{\mu_A}} (128058 t_{2,2}+367230 x_2)+e^{23 t_{\mu_A}} (-18359274 t_{3,4}-34006980 x_3)+886406335 e^{29 t_{\mu_A}},
\end{align*}

\begin{align*}
&\phi^2_{48}=e^{2 t_{\mu_A}} \left(-\frac{4}{15} t_{2,2} t_{3,3}-\frac{4}{75} t_{2,2} t_{3,4}^2-\frac{2}{5} t_{3,3} x_2-\frac{2}{25} t_{3,4}^2 x_2-\frac{2}{5} t_{2,2} t_{3,4} x_3-\frac{2}{3} t_{2,2} x_3^2\right)+e^{3 t_{\mu_A}} \left(-\frac{9}{5} t_{1,1} t_{3,4}-\frac{6}{5} t_{3,4} x_1\right.\\
&-3 t_{1,1} x_3-2 x_1 x_3\Biggr)+e^{4 t_{\mu_A}} \left(-2 t_{2,1}-\frac{t_{2,2}^2}{3}-2 t_{2,2} x_2-9 x_2^2\right)+e^{6 t_{\mu_A}} \left(\frac{114 t_{3,2}}{5}+\frac{444}{25} t_{3,3} t_{3,4}+\frac{394 t_{3,4}^3}{125}+\frac{198}{5} t_{3,3} x_3\right.\\
&\left.+\frac{513}{25} t_{3,4}^2 x_3+\frac{282}{5} t_{3,4} x_3^2+54 x_3^3\right)+e^{8 t_{\mu_A}} \left(\frac{574}{15} t_{2,2} t_{3,4}+178 t_{3,4} x_2+62 t_{2,2} x_3\right)+e^{9 t_{\mu_A}} (366 t_{1,1}+492 x_1)
\end{align*}

\begin{align*}
&+e^{12 t_{\mu_A}} \left(-3954 t_{3,3}-\frac{14043 t_{3,4}^2}{5}-8680 t_{3,4} x_3-10130 x_3^2\right)+e^{14 t_{\mu_A}} \left(-\frac{20816 t_{2,2}}{3}-20149 x_2\right)\\
&+e^{18 t_{\mu_A}} \left(\frac{3055182 t_{3,4}}{5}+756580 x_3\right)-33881640 e^{24 t_{\mu_A}},\\
&\phi^3_{48}=e^{2 t_{\mu_A}} \left(-\frac{2 t_{3,2}}{5}-\frac{2}{25} t_{3,3} t_{3,4}-\frac{t_{3,4}^3}{375}-\frac{2}{5} t_{3,3} x_3-\frac{2}{25} t_{3,4}^2 x_3-\frac{3}{5} t_{3,4} x_3^2-x_3^3\right)+e^{4 t_{\mu_A}} \left(\frac{2}{5} t_{2,2} t_{3,4}+\frac{6}{5} t_{3,4} x_2\right.\\
&+2 t_{2,2} x_3\Biggr)+e^{5 t_{\mu_A}} (t_{1,1}+6 x_1)+e^{8 t_{\mu_A}} \left(-\frac{16 t_{3,3}}{5}-\frac{346 t_{3,4}^2}{25}+\frac{134}{5} t_{3,4} x_3+138 x_3^2\right)+e^{10 t_{\mu_A}} \left(-\frac{550 t_{2,2}}{3}-172 x_2\right)\\
&+e^{14 t_{\mu_A}} \left(\frac{20256 t_{3,4}}{5}-4417 x_3\right)-190116 e^{20 t_{\mu_A}}.\\
&\phi^1_{49}=-\frac{1}{3} t_{2,2} x_1 x_2-x_1 x_2^2+e^{2 t_{\mu_A}} \left(-2 t_{1,1} t_{3,2}-\frac{4}{5} t_{1,1} t_{3,3} t_{3,4}-\frac{4}{75} t_{1,1} t_{3,4}^3-3 t_{1,1} t_{3,3} x_3-\frac{9}{10} t_{1,1} t_{3,4}^2 x_3-3 t_{3,3} x_1 x_3\right.\\
&\left.-\frac{9}{10} t_{3,4}^2 x_1 x_3-4 t_{1,1} t_{3,4} x_3^2-4 t_{3,4} x_1 x_3^2-5 t_{1,1} x_3^3-5 x_1 x_3^3\right)+e^{3 t_{\mu_A}} \left(-2 t_{2,2}^2 t_{3,3}-\frac{3}{5} t_{2,2}^2 t_{3,4}^2-\frac{8}{3} t_{2,2}^2 t_{3,4} x_3-\frac{10}{3} t_{2,2}^2 x_3^2\right)\\
&+e^{4 t_{\mu_A}} (-8 t_{1,1} t_{2,2} t_{3,4}-12 t_{1,1} t_{3,4} x_2+27 t_{3,4} x_1 x_2-5 t_{1,1} t_{2,2} x_3-5 t_{2,2} x_1 x_3-15 t_{1,1} x_2 x_3-60 x_1 x_2 x_3)\\
&+e^{5 t_{\mu_A}} \left(5 t_{1,1}^2+10 t_{3,2} t_{3,3}+4 t_{3,3}^2 t_{3,4}+3 t_{3,2} t_{3,4}^2+\frac{22}{15} t_{3,3} t_{3,4}^3+\frac{2 t_{3,4}^5}{25}+30 t_{2,1} x_2-5 t_{2,2}^2 x_2+30 t_{3,1} x_3+45 t_{3,3}^2 x_3\right.\\
&-18 t_{3,2} t_{3,4} x_3+\frac{81}{5} t_{3,3} t_{3,4}^2 x_3+\frac{63}{20} t_{3,4}^4 x_3+30 t_{3,2} x_3^2+68 t_{3,3} t_{3,4} x_3^2+\frac{88}{5} t_{3,4}^3 x_3^2+160 t_{3,3} x_3^3+48 t_{3,4}^2 x_3^3+120 t_{3,4} x_3^4\\
&+150 x_3^5\Biggr)+e^{7 t_{\mu_A}} \left(214 t_{2,2} t_{3,2}+\frac{708}{5} t_{2,2} t_{3,3} t_{3,4}+\frac{1688}{75} t_{2,2} t_{3,4}^3+174 t_{3,2} x_2+\frac{528}{5} t_{3,3} t_{3,4} x_2\right.+\frac{386}{25} t_{3,4}^3 x_2\\
&+186 t_{2,2} t_{3,3} x_3+\frac{1079}{5} t_{2,2} t_{3,4}^2 x_3+576 t_{3,3} x_2 x_3+\frac{1104}{5} t_{3,4}^2 x_2 x_3+528 t_{2,2} t_{3,4} x_3^2+828 t_{3,4} x_2 x_3^2+260 t_{2,2} x_3^3\\
&+960 x_2 x_3^3\Biggr)+e^{8 t_{\mu_A}} \biggl(296 t_{1,1} t_{3,3}+\frac{1844}{5} t_{1,1} t_{3,4}^2+288 t_{3,3} x_1+\frac{432}{5} t_{3,4}^2 x_1+963 t_{1,1} t_{3,4} x_3+1674 t_{3,4} x_1 x_3+760 t_{1,1} x_3^2\\
&+1555 x_1 x_3^2\biggr)+e^{9 t_{\mu_A}} \left(-360 t_{2,1} t_{3,4}-108 t_{2,2}^2 t_{3,4}-882 t_{2,2} t_{3,4} x_2-450 t_{2,1} x_3-\frac{385}{3} t_{2,2}^2 x_3-1200 t_{2,2} x_2 x_3\right)\\
&+e^{10 t_{\mu_A}} (700 t_{1,1} t_{2,2}+1440 t_{2,2} x_1+1305 t_{1,1} x_2+1485 x_1 x_2)+e^{11 t_{\mu_A}} \Biggl(-1800 t_{3,1}-2748 t_{3,3}^2-4016 t_{3,2} t_{3,4}\\
&-\frac{20846}{5} t_{3,3} t_{3,4}^2-\frac{50116 t_{3,4}^4}{75}-810 t_{3,2} x_3-21180 t_{3,3} t_{3,4} x_3-\frac{38712}{5} t_{3,4}^3 x_3-25280 t_{3,3} x_3^2-28150 t_{3,4}^2 x_3^2\\
&-47570 t_{3,4} x_3^3-39000 x_3^4\Biggr)+e^{13 t_{\mu_A}} \left(-25496 t_{2,2} t_{3,3}-\frac{25184}{5} t_{2,2} t_{3,4}^2-40644 t_{3,3} x_2-\frac{51876}{5} t_{3,4}^2 x_2-58138 t_{2,2} t_{3,4} x_3\right.\\
&-109992 t_{3,4} x_2 x_3-71760 t_{2,2} x_3^2-175590 x_2 x_3^2\biggr)+e^{14 t_{\mu_A}} (-85818 t_{1,1} t_{3,4}-108000 t_{3,4} x_1-85755 t_{1,1} x_3\\
&-199035 x_1 x_3)+e^{15 t_{\mu_A}} \left(41400 t_{2,1}+8220 t_{2,2}^2+24660 t_{2,2} x_2-124740 x_2^2\right)+e^{17 t_{\mu_A}} \left(429246 t_{3,2}+\frac{7504862}{5} t_{3,3} t_{3,4}\right.\\
&\left.+\frac{12203119 t_{3,4}^3}{25}+2577924 t_{3,3} x_3+\frac{18217296}{5} t_{3,4}^2 x_3+7882462 t_{3,4} x_3^2+6053940 x_3^3\right)+e^{19 t_{\mu_A}} (1836892 t_{2,2} t_{3,4}\\
&+6134130 t_{3,4} x_2+4371240 t_{2,2} x_3+11865240 x_2 x_3)+e^{20 t_{\mu_A}} (5316640 t_{1,1}+10692000 x_1)+e^{23 t_{\mu_A}} \biggl(-116485460 t_{3,3}\\
&\left.-155871716 t_{3,4}^2-584406000 t_{3,4} x_3-576147300 x_3^2\right)+e^{25 t_{\mu_A}} (-157960260 t_{2,2}-520036830 x_2)\\
&+e^{29 t_{\mu_A}} (19266865810 t_{3,4}+29950525410 x_3)-786837196800 e^{35 t_{\mu_A}},\\
\end{align*}

\begin{align*}
&\phi^2_{49}=e^{2 t_{\mu_A}} \left(-\frac{4}{3} t_{2,2} t_{3,2}-\frac{8}{15} t_{2,2} t_{3,3} t_{3,4}-\frac{8}{225} t_{2,2} t_{3,4}^3-2 t_{3,2} x_2-\frac{4}{5} t_{3,3} t_{3,4} x_2-\frac{4}{75} t_{3,4}^3 x_2-2 t_{2,2} t_{3,3} x_3-\frac{3}{5} t_{2,2} t_{3,4}^2 x_3\right.\\
&\left.-\frac{8}{3} t_{2,2} t_{3,4} x_3^2-\frac{10}{3} t_{2,2} x_3^3\right)+e^{3 t_{\mu_A}} \Biggl(-9 t_{1,1} t_{3,3}-\frac{27}{10} t_{1,1} t_{3,4}^2-6 t_{3,3} x_1-\frac{9}{5} t_{3,4}^2 x_1-12 t_{1,1} t_{3,4} x_3-8 t_{3,4} x_1 x_3-15 t_{1,1} x_3^2\\
&-10 x_1 x_3^2\Biggr)+e^{4 t_{\mu_A}} \left(-8 t_{2,1} t_{3,4}-\frac{4}{3} t_{2,2}^2 t_{3,4}-8 t_{2,2} t_{3,4} x_2-36 t_{3,4} x_2^2-10 t_{2,1} x_3+\frac{5}{3} t_{2,2}^2 x_3\right)+e^{5 t_{\mu_A}} \left(-\frac{5}{3} t_{1,1} t_{2,2}\right.\\
&-10 t_{2,2} x_1-40 t_{1,1} x_2-36 x_1 x_2\Biggr)+e^{6 t_{\mu_A}} \biggl(36 t_{3,3}^2+84 t_{3,2} t_{3,4}+\frac{276}{5} t_{3,3} t_{3,4}^2+\frac{137 t_{3,4}^4}{25}-30 t_{3,2} x_3+198 t_{3,3} t_{3,4} x_3\\
&+\frac{311}{5} t_{3,4}^3 x_3+60 t_{3,3} x_3^2+234 t_{3,4}^2 x_3^2+270 t_{3,4} x_3^3\biggr)+e^{8 t_{\mu_A}} \left(440 t_{2,2} t_{3,3}+\frac{308}{3} t_{2,2} t_{3,4}^2+600 t_{3,3} x_2+712 t_{3,4}^2 x_2\right.\\
&\left.+630 t_{2,2} t_{3,4} x_3+\frac{4100}{3} t_{2,2} x_3^2\right)+e^{9 t_{\mu_A}} (2616 t_{1,1} t_{3,4}+3348 t_{3,4} x_1+3795 t_{1,1} x_3+3110 x_1 x_3)+e^{10 t_{\mu_A}} \Biggl(1670 t_{2,1}\\
&-\frac{245 t_{2,2}^2}{3}+3110 t_{2,2} x_2+7920 x_2^2\Biggr)+e^{12 t_{\mu_A}} \left(-11730 t_{3,2}-34248 t_{3,3} t_{3,4}-\frac{55018 t_{3,4}^3}{5}-38460 t_{3,3} x_3-49128 t_{3,4}^2 x_3\right.\\
&-91320 t_{3,4} x_3^2-39150 x_3^3\biggr)+e^{14 t_{\mu_A}} \left(-\frac{223078}{3} t_{2,2} t_{3,4}-245636 t_{3,4} x_2-124050 t_{2,2} x_3-41580 x_2 x_3\right)\\
&+e^{15 t_{\mu_A}} (-292425 t_{1,1}-398070 x_1)+e^{18 t_{\mu_A}} \left(3320334 t_{3,3}+\frac{22793661 t_{3,4}^2}{5}+10242792 t_{3,4} x_3+7544340 x_3^2\right)\\
&+e^{20 t_{\mu_A}} \left(\frac{23541890 t_{2,2}}{3}+17836720 x_2\right)+e^{24 t_{\mu_A}} (-615527100 t_{3,4}-590410650 x_3)+25783674900 e^{30 t_{\mu_A}},\\
&\phi^3_{49}=e^{t_{\mu_A}} \left(-\frac{1}{3} t_{1,1} t_{2,2}-\frac{2}{3} t_{2,2} x_1-t_{1,1} x_2-2 x_1 x_2\right)+e^{2 t_{\mu_A}} \left(-2 t_{3,1}-\frac{t_{3,3}^2}{5}-\frac{2}{5} t_{3,2} t_{3,4}-\frac{1}{25} t_{3,3} t_{3,4}^2-\frac{t_{3,4}^4}{1500}\right)\\
&+e^{4 t_{\mu_A}} \left(\frac{14}{3} t_{2,2} t_{3,3}+\frac{7}{5} t_{2,2} t_{3,4}^2+14 t_{3,3} x_2+\frac{21}{5} t_{3,4}^2 x_2+8 t_{2,2} t_{3,4} x_3+10 t_{2,2} x_3^2+30 x_2 x_3^2\right)+e^{5 t_{\mu_A}} (19 t_{1,1} t_{3,4}\\
&+54 t_{3,4} x_1-60 x_1 x_3)+e^{6 t_{\mu_A}} \left(54 t_{2,1}+7 t_{2,2}^2+42 t_{2,2} x_2-108 x_2^2\right)+e^{8 t_{\mu_A}} \left(472 t_{3,2}-\frac{406}{5} t_{3,3} t_{3,4}-\frac{5131 t_{3,4}^3}{75}\right.\\
&+960 t_{3,3} x_3+120 t_{3,4}^2 x_3+990 t_{3,4} x_3^2+1500 x_3^3\Biggr)+e^{10 t_{\mu_A}} \biggl(-\frac{6760}{3} t_{2,2} t_{3,4}-604 t_{3,4} x_2-2930 t_{2,2} x_3-1800 x_2 x_3\biggr)\\
&+e^{11 t_{\mu_A}} (-641 t_{1,1}+2970 x_1)+e^{14 t_{\mu_A}} \left(-48924 t_{3,3}+\frac{91584 t_{3,4}^2}{5}-121722 t_{3,4} x_3-244440 x_3^2\right)\\
&+e^{16 t_{\mu_A}} (157950 t_{2,2}-94122 x_2)+e^{20 t_{\mu_A}} (2898126 t_{3,4}+15130620 x_3)-377238060 e^{26 t_{\mu_A}}.
\end{align*}

\[
\phi^1_{55}=0,\
\phi^2_{55}=0,\
\phi^3_{55}=0.\quad
\phi^1_{56}=3 e^{3 t_{\mu_A}},\
\phi^2_{56}=0,\
\phi^3_{56}=0.
\]
\[
\phi^1_{57}=\frac{6}{5} e^{3 t_{\mu_A}} t_{3,4}-225 e^{9 t_{\mu_A}},\
\phi^2_{57}=8 e^{4 t_{\mu_A}},\
\phi^3_{57}=0.
\]
\[
\phi^1_{58}=e^{3 t_{\mu_A}} \left(\frac{6 t_{3,3}}{5}+\frac{6 t_{3,4}^2}{25}\right)+e^{5 t_{\mu_A}} (15 t_{2,2}+27 x_2)-615 e^{9 t_{\mu_A}} t_{3,4}+56640 e^{15 t_{\mu_A}},\\
\]
\[
\phi^2_{58}=\frac{24}{5} e^{4 t_{\mu_A}} t_{3,4}-1000 e^{10 t_{\mu_A}},\
\phi^3_{58}=\frac{x_3}{5}+54 e^{6 t_{\mu_A}}.\\
\]
\begin{align*}
&\phi^1_{59}=-x_1 x_3+e^{3 t_{\mu_A}} \left(-6 t_{3,2}-\frac{12}{5} t_{3,3} t_{3,4}-\frac{4 t_{3,4}^3}{25}-18 t_{3,3} x_3-\frac{27}{5} t_{3,4}^2 x_3-24 t_{3,4} x_3^2-30 x_3^3\right)+e^{5 t_{\mu_A}} (12 t_{2,2} t_{3,4}\\
&-36 t_{3,4} x_2-30 t_{2,2} x_3-180 x_2 x_3)+60 e^{6 t_{\mu_A}} t_{1,1}+e^{9 t_{\mu_A}} \left(1080 t_{3,3}+660 t_{3,4}^2+6210 t_{3,4} x_3+6000 x_3^2\right)\\
&+e^{11 t_{\mu_A}} (690 t_{2,2}+12240 x_2)+e^{15 t_{\mu_A}} (-341088 t_{3,4}-712530 x_3)+29208000 e^{21 t_{\mu_A}},\\
&\phi^2_{59}=e^{t_{\mu_A}} (-t_{1,1}-2 x_1)+e^{4 t_{\mu_A}} \left(-4 t_{3,3}-\frac{6 t_{3,4}^2}{5}-48 t_{3,4} x_3-60 x_3^2\right)+84 e^{6 t_{\mu_A}} t_{2,2}\\
&+e^{10 t_{\mu_A}} (2420 t_{3,4}+12000 x_3)-474660 e^{16 t_{\mu_A}},\\
&\phi^3_{59}=e^{2 t_{\mu_A}} (-4 t_{2,2}-6 x_2)+e^{6 t_{\mu_A}} (36 t_{3,4}-300 x_3)+14976 e^{12 t_{\mu_A}}.
\end{align*}

\begin{align*}
&\phi^1_{66}=e^{3 t_{\mu_A}} \left(\frac{6 t_{3,4}}{5}+6 x_3\right)-480 e^{9 t_{\mu_A}},\
\phi^2_{66}=12 e^{4 t_{\mu_A}},\
\phi^3_{66}=0.\\
&\phi^1_{67}=e^{3 t_{\mu_A}} \left(\frac{3 t_{3,3}}{5}+\frac{3 t_{3,4}^2}{10}+\frac{12}{5} t_{3,4} x_3+3 x_3^2\right)+e^{5 t_{\mu_A}} (19 t_{2,2}+42 x_2)+e^{9 t_{\mu_A}} (-858 t_{3,4}-585 x_3)+90144 e^{15 t_{\mu_A}},\\
&\phi^2_{67}=e^{4 t_{\mu_A}} \left(\frac{32 t_{3,4}}{5}+14 x_3\right)-1810 e^{10 t_{\mu_A}},\
\phi^3_{67}=\frac{x_3}{5}+60 e^{6 t_{\mu_A}}.\\
\end{align*}

\begin{align*}
&\phi^1_{68}=-\frac{1}{5} x_1 x_3+e^{3 t_{\mu_A}} \left(-\frac{9 t_{3,2}}{5}-\frac{3}{5} t_{3,3} t_{3,4}-\frac{3 t_{3,4}^3}{250}-\frac{6}{5} t_{3,3} x_3-\frac{3}{5} t_{3,4}^2 x_3-3 t_{3,4} x_3^2-3 x_3^3\right)\\
&+e^{5 t_{\mu_A}} \left(\frac{24}{5} t_{2,2} t_{3,4}+\frac{9}{5} t_{3,4} x_2+9 t_{2,2} x_3+9 x_2 x_3\right)+9 e^{6 t_{\mu_A}} t_{1,1}+e^{9 t_{\mu_A}} \left(159 t_{3,3}+\frac{45 t_{3,4}^2}{2}+255 t_{3,4} x_3+705 x_3^2\right)\\
&+e^{11 t_{\mu_A}} (-1106 t_{2,2}-930 x_2)+e^{15 t_{\mu_A}} \left(-\frac{13326 t_{3,4}}{5}-37866 x_3\right)+375135 e^{21 t_{\mu_A}},\\
&\phi^2_{68}=e^{t_{\mu_A}} \left(-\frac{t_{1,1}}{5}-\frac{2 x_1}{5}\right)+e^{4 t_{\mu_A}} \left(-\frac{4 t_{3,3}}{5}+\frac{6 t_{3,4}^2}{25}-\frac{6}{5} t_{3,4} x_3-6 x_3^2\right)\\
&+e^{6 t_{\mu_A}} \left(\frac{114 t_{2,2}}{5}+27 x_2\right)+e^{10 t_{\mu_A}} (-274 t_{3,4}+410 x_3)+16148 e^{16 t_{\mu_A}},\\
&\phi^3_{68}=\frac{x_3^2}{5}+e^{2 t_{\mu_A}} \left(-\frac{4 t_{2,2}}{5}-\frac{6 x_2}{5}\right)+e^{6 t_{\mu_A}} \left(\frac{54 t_{3,4}}{5}-9 x_3\right)-\frac{3804}{5} e^{12 t_{\mu_A}}.\\
&\phi^1_{69}=-\frac{1}{5} t_{3,4} x_1 x_3-x_1 x_3^2+e^{2 t_{\mu_A}} (-2 t_{1,1} t_{2,2}-3 t_{1,1} x_2-3 x_1 x_2)+e^{3 t_{\mu_A}} \left(-\frac{6}{5} t_{3,2} t_{3,4}-\frac{12}{25} t_{3,3} t_{3,4}^2-\frac{4 t_{3,4}^4}{125}+6 t_{3,2} x_3\right.\\
&\left.-\frac{6}{5} t_{3,3} t_{3,4} x_3-\frac{23}{25} t_{3,4}^3 x_3-\frac{24}{5} t_{3,4}^2 x_3^2-6 t_{3,4} x_3^3\right)+e^{5 t_{\mu_A}} \left(10 t_{2,2} t_{3,3}+\frac{27}{5} t_{2,2} t_{3,4}^2+42 t_{3,3} x_2+\frac{27}{5} t_{3,4}^2 x_2-18 t_{2,2} t_{3,4} x_3\right.\\
&-90 t_{2,2} x_3^2-90 x_2 x_3^2\biggr)+e^{6 t_{\mu_A}} (30 t_{1,1} t_{3,4}-60 t_{1,1} x_3+20 x_1 x_3)+e^{7 t_{\mu_A}} \left(168 t_{2,2}^2+522 t_{2,2} x_2+540 x_2^2\right)\\
&+e^{9 t_{\mu_A}} \left(-60 t_{3,2}-42 t_{3,3} t_{3,4}+\frac{301 t_{3,4}^3}{5}+360 t_{3,3} x_3+1230 t_{3,4}^2 x_3+4650 t_{3,4} x_3^2+1500 x_3^3\right)+e^{11 t_{\mu_A}} (-8474 t_{2,2} t_{3,4}\\
&-16302 t_{3,4} x_2+5670 t_{2,2} x_3+3240 x_2 x_3)-2712 e^{12 t_{\mu_A}} t_{1,1}+e^{15 t_{\mu_A}} \left(4500 t_{3,3}+\frac{236862 t_{3,4}^2}{5}-399738 t_{3,4} x_3\right.\\
&-487740 x_3^2\biggr)+e^{17 t_{\mu_A}} (683538 t_{2,2}+1398690 x_2)+e^{21 t_{\mu_A}} (-6361830 t_{3,4}+31002000 x_3)+73886880 e^{27 t_{\mu_A}},
\end{align*}

\begin{align*}
&\phi^2_{69}=e^{t_{\mu_A}} \left(-\frac{1}{5} t_{1,1} t_{3,4}-\frac{2}{5} t_{3,4} x_1-t_{1,1} x_3-2 x_1 x_3\right)+e^{2 t_{\mu_A}} \left(-2 t_{2,1}-\frac{t_{2,2}^2}{3}\right)+e^{4 t_{\mu_A}} \left(6 t_{3,2}+\frac{8}{5} t_{3,3} t_{3,4}\right.\\
&\left.-\frac{2 t_{3,4}^3}{25}+14 t_{3,3} x_3-\frac{27}{5} t_{3,4}^2 x_3-36 t_{3,4} x_3^2-30 x_3^3\right)+e^{6 t_{\mu_A}} \left(\frac{354}{5} t_{2,2} t_{3,4}+108 t_{3,4} x_2+54 t_{2,2} x_3+180 x_2 x_3\right)\\
&+e^{7 t_{\mu_A}} (68 t_{1,1}+40 x_1)+e^{10 t_{\mu_A}} \left(-370 t_{3,3}-683 t_{3,4}^2+2270 t_{3,4} x_3+9000 x_3^2\right)+e^{12 t_{\mu_A}} (-12570 t_{2,2}\\
&-23760 x_2)+e^{16 t_{\mu_A}} (165932 t_{3,4}-452490 x_3)-3828000 e^{22 t_{\mu_A}},\\
&\phi^3_{69}=e^{2 t_{\mu_A}} \left(-\frac{4}{5} t_{2,2} t_{3,4}-\frac{6}{5} t_{3,4} x_2-2 t_{2,2} x_3\right)+e^{3 t_{\mu_A}} (-9 t_{1,1}-6 x_1)+e^{6 t_{\mu_A}} \left(36 t_{3,3}+18 t_{3,4}^2-18 t_{3,4} x_3-300 x_3^2\right)\\
&+e^{8 t_{\mu_A}} (836 t_{2,2}+1200 x_2)+e^{12 t_{\mu_A}} \left(-\frac{63324 t_{3,4}}{5}+20988 x_3\right)+492660 e^{18 t_{\mu_A}}.\\
\end{align*}

\begin{align*}
&\phi^1_{77}=-\frac{1}{5} x_1 x_3+e^{3 t_{\mu_A}} \left(-\frac{12 t_{3,2}}{5}-\frac{12}{25} t_{3,3} t_{3,4}-\frac{2 t_{3,4}^3}{125}-\frac{18}{5} t_{3,3} x_3-\frac{3}{25} t_{3,4}^2 x_3-\frac{12}{5} t_{3,4} x_3^2-6 x_3^3\right)+e^{5 t_{\mu_A}} \left(\frac{28}{5} t_{2,2} t_{3,4}\right.\\
&\left.+\frac{24}{5} t_{3,4} x_2+26 t_{2,2} x_3+12 x_2 x_3\right)+e^{6 t_{\mu_A}} (36 t_{1,1}+32 x_1)+e^{9 t_{\mu_A}} \left(102 t_{3,3}-\frac{369 t_{3,4}^2}{5}-90 t_{3,4} x_3+750 x_3^2\right)\\
&+e^{11 t_{\mu_A}} (-1854 t_{2,2}-1260 x_2)+e^{15 t_{\mu_A}} \left(\frac{126678 t_{3,4}}{5}-12474 x_3\right)-1220640 e^{21 t_{\mu_A}},\\
&\phi^2_{77}=e^{t_{\mu_A}} \left(-\frac{t_{1,1}}{5}-\frac{2 x_1}{5}\right)+e^{4 t_{\mu_A}} \left(-\frac{12 t_{3,3}}{5}+\frac{14 t_{3,4}^2}{25}+\frac{8}{5} t_{3,4} x_3+4 x_3^2\right)+\frac{44}{5} e^{6 t_{\mu_A}} t_{2,2}\\
&+e^{10 t_{\mu_A}} (-172 t_{3,4}-420 x_3)+16884 e^{16 t_{\mu_A}},\\
&\phi^3_{77}=\frac{x_3^2}{5}+e^{2 t_{\mu_A}} \left(-\frac{4 t_{2,2}}{5}-\frac{6 x_2}{5}\right)+e^{6 t_{\mu_A}} (12 t_{3,4}-8 x_3)-\frac{4104}{5} e^{12 t_{\mu_A}}.
\end{align*}

\begin{align*}
&\phi^1_{78}=-\frac{1}{25} t_{3,4} x_1 x_3-\frac{1}{5} x_1 x_3^2+e^{2 t_{\mu_A}} \left(-\frac{2}{5} t_{1,1} t_{2,2}-\frac{3}{5} t_{1,1} x_2-\frac{6}{5} x_1 x_2\right)+e^{3 t_{\mu_A}} \left(\frac{6 t_{3,3}^2}{25}-\frac{6}{25} t_{3,2} t_{3,4}-\frac{9}{125} t_{3,3} t_{3,4}^2\right.\\
&\left.-\frac{4 t_{3,4}^4}{625}+\frac{6}{25} t_{3,3} t_{3,4} x_3-\frac{3}{125} t_{3,4}^3 x_3+\frac{6}{5} t_{3,3} x_3^2\right)+e^{5 t_{\mu_A}} \left(-\frac{3}{5} t_{2,2} t_{3,3}+\frac{11}{10} t_{2,2} t_{3,4}^2+\frac{27}{5} t_{3,3} x_2+\frac{159}{50} t_{3,4}^2 x_2\right.\\
&+\frac{24}{5} t_{2,2} t_{3,4} x_3+\frac{54}{5} t_{3,4} x_2 x_3+7 t_{2,2} x_3^2+15 x_2 x_3^2\biggr)+e^{6 t_{\mu_A}} \left(\frac{84}{5} t_{1,1} t_{3,4}+\frac{96}{5} t_{3,4} x_1+15 t_{1,1} x_3+47 x_1 x_3\right)\\
&+e^{7 t_{\mu_A}} \left(\frac{3 t_{2,2}^2}{5}-\frac{63}{5} t_{2,2} x_2\right)+e^{9 t_{\mu_A}} \left(87 t_{3,2}-\frac{384}{5} t_{3,3} t_{3,4}-\frac{866 t_{3,4}^3}{25}-156 t_{3,3} x_3-\frac{849}{5} t_{3,4}^2 x_3-228 t_{3,4} x_3^2+105 x_3^3\right)\\
&+e^{11 t_{\mu_A}} \left(-\frac{2531}{5} t_{2,2} t_{3,4}-1365 t_{3,4} x_2-1670 t_{2,2} x_3-2370 x_2 x_3\right)+e^{12 t_{\mu_A}} \left(-\frac{14202 t_{1,1}}{5}-4000 x_1\right)\\
&+e^{15 t_{\mu_A}} \left(\frac{68847 t_{3,3}}{5}+\frac{217941 t_{3,4}^2}{10}+\frac{287124}{5} t_{3,4} x_3+8757 x_3^2\right)+e^{17 t_{\mu_A}} \left(\frac{404558 t_{2,2}}{5}+141738 x_2\right)\\
&+e^{21 t_{\mu_A}} (-3840648 t_{3,4}-3069810 x_3)+175921251 e^{27 t_{\mu_A}},\\
\end{align*}

\begin{align*}
&\phi^2_{78}=e^{t_{\mu_A}} \left(-\frac{1}{25} t_{1,1} t_{3,4}-\frac{2}{25} t_{3,4} x_1-\frac{1}{5} t_{1,1} x_3-\frac{2}{5} x_1 x_3\right)+e^{2 t_{\mu_A}} \left(-\frac{2 t_{2,1}}{5}-\frac{t_{2,2}^2}{15}\right)+e^{4 t_{\mu_A}} \left(-2 t_{3,2}+\frac{16}{25} t_{3,3} t_{3,4}\right.\\
&\left.+\frac{22 t_{3,4}^3}{375}+\frac{22}{25} t_{3,4}^2 x_3-4 x_3^3\right)+e^{6 t_{\mu_A}} \left(\frac{174}{25} t_{2,2} t_{3,4}+\frac{54}{5} t_{3,4} x_2+\frac{114}{5} t_{2,2} x_3\right)+e^{7 t_{\mu_A}} \left(\frac{353 t_{1,1}}{5}+94 x_1\right)\\
&+e^{10 t_{\mu_A}} \left(-298 t_{3,3}-\frac{1309 t_{3,4}^2}{5}-572 t_{3,4} x_3+290 x_3^2\right)+e^{12 t_{\mu_A}} (-2134 t_{2,2}-2025 x_2)+e^{16 t_{\mu_A}} (74992 t_{3,4}+30398 x_3)\\
&-4046580 e^{22 t_{\mu_A}},\\
&\phi^3_{78}=\frac{1}{25} t_{3,4} x_3^2+\frac{x_3^3}{5}+e^{2 t_{\mu_A}} \left(-\frac{4}{25} t_{2,2} t_{3,4}-\frac{6}{25} t_{3,4} x_2-\frac{4}{5} t_{2,2} x_3\right)+e^{3 t_{\mu_A}} \left(-\frac{12 t_{1,1}}{5}-\frac{12 x_1}{5}\right)+e^{6 t_{\mu_A}} \left(\frac{54 t_{3,3}}{5}\right.\\
&\left.+\frac{147 t_{3,4}^2}{25}+\frac{66}{5} t_{3,4} x_3-17 x_3^2\right)+e^{8 t_{\mu_A}} \left(\frac{358 t_{2,2}}{5}+66 x_2\right)+e^{12 t_{\mu_A}} \left(-\frac{54654 t_{3,4}}{25}-\frac{3079 x_3}{5}\right)+123324 e^{18 t_{\mu_A}}.
\end{align*}

\begin{align*}
&\phi^1_{79}=-\frac{1}{5} t_{3,3} x_1 x_3-\frac{1}{50} t_{3,4}^2 x_1 x_3-\frac{2}{5} t_{3,4} x_1 x_3^2-x_1 x_3^3+e^{2 t_{\mu_A}} \left(-\frac{4}{5} t_{1,1} t_{2,2} t_{3,4}-\frac{6}{5} t_{1,1} t_{3,4} x_2-\frac{6}{5} t_{3,4} x_1 x_2-2 t_{1,1} t_{2,2} x_3\right.\\
&-3 t_{2,2} x_1 x_3-3 t_{1,1} x_2 x_3-12 x_1 x_2 x_3\Biggr)+e^{3 t_{\mu_A}} \biggl(-3 t_{1,1}^2-\frac{6}{5} t_{3,2} t_{3,3}-\frac{12}{25} t_{3,3}^2 t_{3,4}-\frac{3}{25} t_{3,2} t_{3,4}^2-\frac{2}{25} t_{3,3} t_{3,4}^3-\frac{2 t_{3,4}^5}{625}\\
&+8 t_{2,2}^2 x_2+24 t_{2,2} x_2^2+18 x_2^3-\frac{18}{5} t_{3,3}^2 x_3+\frac{12}{5} t_{3,2} t_{3,4} x_3-\frac{12}{25} t_{3,3} t_{3,4}^2 x_3-\frac{11}{250} t_{3,4}^4 x_3+6 t_{3,2} x_3^2-\frac{12}{5} t_{3,3} t_{3,4} x_3^2\\
&\left.-\frac{8}{25} t_{3,4}^3 x_3^2-6 t_{3,3} x_3^3-\frac{3}{5} t_{3,4}^2 x_3^3\right)+e^{5 t_{\mu_A}} \left(2 t_{2,2} t_{3,2}+\frac{36}{5} t_{2,2} t_{3,3} t_{3,4}+\frac{112}{75} t_{2,2} t_{3,4}^3+12 t_{3,2} x_2+\frac{72}{5} t_{3,3} t_{3,4} x_2\right.\\
&\left.+\frac{116}{25} t_{3,4}^3 x_2+52 t_{2,2} t_{3,3} x_3+12 t_{2,2} t_{3,4}^2 x_3+24 t_{3,3} x_2 x_3+\frac{144}{5} t_{3,4}^2 x_2 x_3+16 t_{2,2} t_{3,4} x_3^2+24 t_{3,4} x_2 x_3^2-10 t_{2,2} x_3^3-60 x_2 x_3^3\right)\\
&+e^{6 t_{\mu_A}} \biggl(96 t_{1,1} t_{3,3}+\frac{168}{5} t_{1,1} t_{3,4}^2+96 t_{3,3} x_1+\frac{144}{5} t_{3,4}^2 x_1+90 t_{1,1} t_{3,4} x_3+257 t_{3,4} x_1 x_3+60 t_{1,1} x_3^2+300 x_1 x_3^2\biggr)\\
&+e^{7 t_{\mu_A}} \left(-\frac{84}{5} t_{2,2}^2 t_{3,4}-\frac{1656}{5} t_{2,2} t_{3,4} x_2-432 t_{3,4} x_2^2+118 t_{2,2}^2 x_3+312 t_{2,2} x_2 x_3\right)+e^{8 t_{\mu_A}} (424 t_{1,1} t_{2,2}+480 t_{2,2} x_1\\
&+780 t_{1,1} x_2+1575 x_1 x_2)+e^{9 t_{\mu_A}} \left(-264 t_{3,3}^2-6 t_{3,2} t_{3,4}-\frac{1152}{5} t_{3,3} t_{3,4}^2-\frac{1456 t_{3,4}^4}{25}-720 t_{3,2} x_3-516 t_{3,3} t_{3,4} x_3\right.\\
&-384 t_{3,4}^3 x_3-120 t_{3,3} x_3^2-24 t_{3,4}^2 x_3^2+1590 t_{3,4} x_3^3\biggr)+e^{11 t_{\mu_A}} \biggl(-5028 t_{2,2} t_{3,3}-\frac{11004}{5} t_{2,2} t_{3,4}^2-7362 t_{3,3} x_2\\
&-\frac{14571}{5} t_{3,4}^2 x_2-15608 t_{2,2} t_{3,4} x_3-15744 t_{3,4} x_2 x_3-4060 t_{2,2} x_3^2-510 x_2 x_3^2\biggr)+e^{12 t_{\mu_A}} \left(-\frac{120024}{5} t_{1,1} t_{3,4}-36000 t_{3,4} x_1\right.\\
&-12912 t_{1,1} x_3-41388 x_1 x_3\biggr)+e^{13 t_{\mu_A}} \left(-6918 t_{2,2}^2-3294 t_{2,2} x_2+16524 x_2^2\right)+e^{15 t_{\mu_A}} \biggl(14292 t_{3,2}+\frac{873846}{5} t_{3,3} t_{3,4}\\
&+\frac{2109309 t_{3,4}^3}{25}+87942 t_{3,3} x_3+286941 t_{3,4}^2 x_3-44814 t_{3,4} x_3^2-168660 x_3^3\biggr)+e^{17 t_{\mu_A}} \left(\frac{5938786}{5} t_{2,2} t_{3,4}+1473906 t_{3,4} x_2\right.\\
&+1595538 t_{2,2} x_3+1224180 x_2 x_3\Biggr)+e^{18 t_{\mu_A}} (2077020 t_{1,1}+3564000 x_1)+e^{21 t_{\mu_A}} \left(-14782620 t_{3,3}-31873638 t_{3,4}^2\right.\\
&\left.-42125310 t_{3,4} x_3+8344500 x_3^2\right)+e^{23 t_{\mu_A}} (-97157950 t_{2,2}-110617962 x_2)+e^{27 t_{\mu_A}} (3797143338 t_{3,4}+1547617680 x_3)\\
&-136685952000 e^{33 t_{\mu_A}},\\
\end{align*}

\begin{align*}
&\phi^2_{79}=e^{t_{\mu_A}} \left(-\frac{1}{5} t_{1,1} t_{3,3}-\frac{1}{50} t_{1,1} t_{3,4}^2-\frac{2}{5} t_{3,3} x_1-\frac{1}{25} t_{3,4}^2 x_1-\frac{2}{5} t_{1,1} t_{3,4} x_3-\frac{4}{5} t_{3,4} x_1 x_3-t_{1,1} x_3^2-2 x_1 x_3^2\right)\\
&+e^{2 t_{\mu_A}} \left(-\frac{4}{5} t_{2,1} t_{3,4}-\frac{2}{15} t_{2,2}^2 t_{3,4}-2 t_{2,1} x_3-\frac{1}{3} t_{2,2}^2 x_3+4 t_{2,2} x_2 x_3+6 x_2^2 x_3\right)+e^{3 t_{\mu_A}} (-7 t_{1,1} t_{2,2}-6 t_{2,2} x_1-6 t_{1,1} x_2\\
&-16 x_1 x_2)+e^{4 t_{\mu_A}} \left(-\frac{4}{5} t_{3,3}^2+\frac{12}{5} t_{3,2} t_{3,4}+\frac{16}{25} t_{3,3} t_{3,4}^2+\frac{t_{3,4}^4}{25}+18 t_{3,2} x_3+\frac{16}{5} t_{3,3} t_{3,4} x_3+\frac{6}{5} t_{3,4}^3 x_3+20 t_{3,3} x_3^2\right.\\
&\left.-\frac{6}{5} t_{3,4}^2 x_3^2-12 t_{3,4} x_3^3\right)+e^{6 t_{\mu_A}} \left(\frac{184}{5} t_{2,2} t_{3,3}+\frac{732}{25} t_{2,2} t_{3,4}^2+48 t_{3,3} x_2+\frac{288}{5} t_{3,4}^2 x_2+\frac{458}{5} t_{2,2} t_{3,4} x_3-36 t_{3,4} x_2 x_3\right.\\
&+44 t_{2,2} x_3^2\biggr)+e^{7 t_{\mu_A}} \left(\frac{1561}{5} t_{1,1} t_{3,4}+514 t_{3,4} x_1+228 t_{1,1} x_3+600 x_1 x_3\right)+e^{8 t_{\mu_A}} \left(120 t_{2,1}+128 t_{2,2}^2+144 t_{2,2} x_2\right)\\
&+e^{10 t_{\mu_A}} \left(-780 t_{3,2}-1956 t_{3,3} t_{3,4}-\frac{5166 t_{3,4}^3}{5}-2280 t_{3,3} x_3-2608 t_{3,4}^2 x_3+260 t_{3,4} x_3^2+3750 x_3^3\right)\\
&+e^{12 t_{\mu_A}} (-18148 t_{2,2} t_{3,4}-23604 t_{3,4} x_2-25060 t_{2,2} x_3-3996 x_2 x_3)+e^{13 t_{\mu_A}} (-46692 t_{1,1}-82776 x_1)\\
&+e^{16 t_{\mu_A}} \left(300708 t_{3,3}+\frac{2668822 t_{3,4}^2}{5}+677214 t_{3,4} x_3-257400 x_3^2\right)+e^{18 t_{\mu_A}} (2113848 t_{2,2}+2138400 x_2)\\
&+e^{22 t_{\mu_A}} (-80787900 t_{3,4}-28500750 x_3)+3363095520 e^{28 t_{\mu_A}},\\
&\phi^3_{79}=e^{2 t_{\mu_A}} \left(-\frac{4}{5} t_{2,2} t_{3,3}-\frac{2}{25} t_{2,2} t_{3,4}^2-\frac{6}{5} t_{3,3} x_2-\frac{3}{25} t_{3,4}^2 x_2-\frac{4}{5} t_{2,2} t_{3,4} x_3-2 t_{2,2} x_3^2\right)+e^{3 t_{\mu_A}} \left(-\frac{18}{5} t_{1,1} t_{3,4}-\frac{12}{5} t_{3,4} x_1\right.\\
&-6 t_{1,1} x_3\Biggr)+e^{4 t_{\mu_A}} \left(-16 t_{2,1}+\frac{4 t_{2,2}^2}{3}-2 t_{2,2} x_2-12 x_2^2\right)+e^{6 t_{\mu_A}} \left(24 t_{3,2}+\frac{156}{5} t_{3,3} t_{3,4}+\frac{142 t_{3,4}^3}{25}-24 t_{3,3} x_3+\frac{108}{5} t_{3,4}^2 x_3\right.\\
&-78 t_{3,4} x_3^2\-300 x_3^3\biggr)+e^{8 t_{\mu_A}} \left(\frac{1672}{5} t_{2,2} t_{3,4}+678 t_{3,4} x_2+1224 t_{2,2} x_3\right)+e^{9 t_{\mu_A}} (2559 t_{1,1}+3150 x_1)+e^{12 t_{\mu_A}} \left(-\frac{65124 t_{3,3}}{5}\right.\\
&\left.-\frac{289362 t_{3,4}^2}{25}-\frac{67224}{5} t_{3,4} x_3+32988 x_3^2\right)+e^{14 t_{\mu_A}} (-101052 t_{2,2}-101232 x_2)+e^{18 t_{\mu_A}} (2974884 t_{3,4}-24060 x_3)\\
&-146515476 e^{24 t_{\mu_A}}.
\end{align*}

\begin{align*}
&\phi^1_{88}=-\frac{1}{25} t_{3,3} x_1 x_3-\frac{1}{250} t_{3,4}^2 x_1 x_3-\frac{2}{25} t_{3,4} x_1 x_3^2-\frac{1}{5} x_1 x_3^3+e^{2 t_{\mu_A}} \left(-\frac{4}{25} t_{1,1} t_{2,2} t_{3,4}-\frac{6}{25} t_{1,1} t_{3,4} x_2-\frac{12}{25} t_{3,4} x_1 x_2\right.\\
&\left.-\frac{2}{5} t_{1,1} t_{2,2} x_3-\frac{3}{5} t_{2,2} x_1 x_3-\frac{3}{5} t_{1,1} x_2 x_3-\frac{12}{5} x_1 x_2 x_3\right)+e^{3 t_{\mu_A}} \left(-\frac{3}{5} t_{1,1}^2-\frac{12}{125} t_{3,3}^2 t_{3,4}+\frac{6}{125} t_{3,2} t_{3,4}^2-\frac{6}{625} t_{3,3} t_{3,4}^3\right.\\
&-\frac{2 t_{3,4}^5}{3125}+\frac{8}{5} t_{2,2}^2 x_2+\frac{24}{5} t_{2,2} x_2^2+\frac{18 x_2^3}{5}+\frac{6}{25} t_{3,3}^2 x_3+\frac{12}{125} t_{3,3} t_{3,4}^2 x_3+\frac{21 t_{3,4}^4 x_3}{1250}+\frac{12}{25} t_{3,3} t_{3,4} x_3^2+\frac{24}{125} t_{3,4}^3 x_3^2\\
&\left.+\frac{6}{5} t_{3,3} x_3^3+\frac{9}{25} t_{3,4}^2 x_3^3\right)+e^{5 t_{\mu_A}} \left(-\frac{6}{5} t_{2,2} t_{3,2}+\frac{26}{25} t_{2,2} t_{3,3} t_{3,4}+\frac{17}{125} t_{2,2} t_{3,4}^3-\frac{18}{5} t_{3,2} x_2+\frac{102}{25} t_{3,3} t_{3,4} x_2\right.\\
&\left.+\frac{141}{125} t_{3,4}^3 x_2+\frac{16}{5} t_{2,2} t_{3,3} x_3+\frac{96}{25} t_{2,2} t_{3,4}^2 x_3+12 t_{3,3} x_2 x_3+\frac{36}{5} t_{3,4}^2 x_2 x_3+\frac{66}{5} t_{2,2} t_{3,4} x_3^2+18 t_{3,4} x_2 x_3^2+12 t_{2,2} x_3^3+18 x_2 x_3^3\right)\\
&+e^{6 t_{\mu_A}} \left(\frac{12}{5} t_{1,1} t_{3,3}+\frac{48}{5} t_{1,1} t_{3,4}^2+\frac{288}{25} t_{3,4}^2 x_1+\frac{144}{5} t_{1,1} t_{3,4} x_3+\frac{343}{5} t_{3,4} x_1 x_3+18 t_{1,1} x_3^2+62 x_1 x_3^2\right)\\
&+e^{7 t_{\mu_A}} \left(-\frac{414}{25} t_{2,2}^2 t_{3,4}-\frac{2826}{25} t_{2,2} t_{3,4} x_2-\frac{648}{5} t_{3,4} x_2^2-\frac{42}{5} t_{2,2}^2 x_3-\frac{108}{5} t_{2,2} x_2 x_3\right)+e^{8 t_{\mu_A}} \left(\frac{104}{5} t_{1,1} t_{2,2}+6 t_{1,1} x_2\right.\\
&+129 x_1 x_2\Biggr)+e^{9 t_{\mu_A}} \left(-\frac{36}{5} t_{3,3}^2+\frac{294}{5} t_{3,2} t_{3,4}-\frac{1356}{25} t_{3,3} t_{3,4}^2-\frac{499 t_{3,4}^4}{25}-\frac{1572}{5} t_{3,3} t_{3,4} x_3-\frac{4218}{25} t_{3,4}^3 x_3-300 t_{3,3} x_3^2\right.\\
&-486 t_{3,4}^2 x_3^2-468 t_{3,4} x_3^3\biggr)+e^{11 t_{\mu_A}} \left(-\frac{1338}{5} t_{2,2} t_{3,3}+\frac{1581}{5} t_{2,2} t_{3,4}^2-936 t_{3,3} x_2+702 t_{3,4}^2 x_2-\frac{9536}{5} t_{2,2} t_{3,4} x_3\right.\\
&-3084 t_{3,4} x_2 x_3-3202 t_{2,2} x_3^2-3840 x_2 x_3^2\Biggr)+e^{12 t_{\mu_A}} \left(-\frac{64104}{25} t_{1,1} t_{3,4}-4800 t_{3,4} x_1-\frac{23952}{5} t_{1,1} x_3-\frac{57378}{5} x_1 x_3\right)
\end{align*}

\begin{align*}
&+e^{13 t_{\mu_A}} \left(\frac{12132 t_{2,2}^2}{5}+\frac{60246}{5} t_{2,2} x_2+\frac{57024 x_2^2}{5}\right)+e^{15 t_{\mu_A}} \left(-\frac{35256 t_{3,2}}{5}+\frac{584526}{25} t_{3,3} t_{3,4}+\frac{1315567 t_{3,4}^3}{125}\right.\\
&+\frac{244866}{5} t_{3,3} x_3+\frac{2187471}{25} t_{3,4}^2 x_3+\frac{691266}{5} t_{3,4} x_3^2+39012 x_3^3\biggr)+e^{17 t_{\mu_A}} \left(-\frac{967924}{25} t_{2,2} t_{3,4}-\frac{459714}{5} t_{3,4} x_2\right.\\
&+\frac{1304718}{5} t_{2,2} x_3+323604 x_2 x_3\biggr)+e^{18 t_{\mu_A}} (215772 t_{1,1}+500000 x_1)+e^{21 t_{\mu_A}} \left(-2200020 t_{3,3}-\frac{15303492 t_{3,4}^2}{5}\right.\\
&-14407686 t_{3,4} x_3-8280720 x_3^2\Biggr)+e^{23 t_{\mu_A}} \left(-2026918 t_{2,2}-\frac{3189792 x_2}{5}\right)+e^{27 t_{\mu_A}} \left(\frac{2026944318 t_{3,4}}{5}+675234492 x_3\right)\\
&-17094510720 e^{33 t_{\mu_A}},\\
&\phi^2_{88}=e^{t_{\mu_A}} \left(-\frac{1}{25} t_{1,1} t_{3,3}-\frac{1}{250} t_{1,1} t_{3,4}^2-\frac{2}{25} t_{3,3} x_1-\frac{1}{125} t_{3,4}^2 x_1-\frac{2}{25} t_{1,1} t_{3,4} x_3-\frac{4}{25} t_{3,4} x_1 x_3-\frac{1}{5} t_{1,1} x_3^2-\frac{2}{5} x_1 x_3^2\right)\\
&+e^{2 t_{\mu_A}} \left(-\frac{4}{25} t_{2,1} t_{3,4}-\frac{2}{75} t_{2,2}^2 t_{3,4}-\frac{2}{5} t_{2,1} x_3-\frac{1}{15} t_{2,2}^2 x_3+\frac{4}{5} t_{2,2} x_2 x_3+\frac{6}{5} x_2^2 x_3\right)+e^{3 t_{\mu_A}} \left(-\frac{7}{5} t_{1,1} t_{2,2}-\frac{6}{5} t_{2,2} x_1\right.\\
&\left.-\frac{6}{5} t_{1,1} x_2-\frac{24}{5} x_1 x_2\right)+e^{4 t_{\mu_A}} \left(\frac{4 t_{3,3}^2}{5}+\frac{12}{25} t_{3,2} t_{3,4}+\frac{16}{125} t_{3,3} t_{3,4}^2+\frac{11 t_{3,4}^4}{625}+\frac{6}{5} t_{3,2} x_3+\frac{76}{25} t_{3,3} t_{3,4} x_3\right.\\
&\left.+\frac{64}{125} t_{3,4}^3 x_3+4 t_{3,3} x_3^2+\frac{78}{25} t_{3,4}^2 x_3^2+\frac{24}{5} t_{3,4} x_3^3\right)+e^{6 t_{\mu_A}} \left(\frac{304}{25} t_{2,2} t_{3,3}+\frac{252}{125} t_{2,2} t_{3,4}^2+\frac{156}{5} t_{3,3} x_2+\frac{36}{5} t_{3,4}^2 x_2\right.\\
&\left.+\frac{258}{25} t_{2,2} t_{3,4} x_3-\frac{108}{5} t_{3,4} x_2 x_3+\frac{144}{5} t_{2,2} x_3^2\right)+e^{7 t_{\mu_A}} \left(\frac{1811}{25} t_{1,1} t_{3,4}+\frac{686}{5} t_{3,4} x_1+\frac{538}{5} t_{1,1} x_3+124 x_1 x_3\right)\\
&+e^{8 t_{\mu_A}} \left(124 t_{2,1}+\frac{274 t_{2,2}^2}{15}+\frac{684}{5} t_{2,2} x_2+243 x_2^2\right)+e^{10 t_{\mu_A}} \left(-256 t_{3,2}-\frac{4604}{5} t_{3,3} t_{3,4}-\frac{15854 t_{3,4}^3}{75}-1212 t_{3,3} x_3\right.\\
&-930 t_{3,4}^2 x_3-1928 t_{3,4} x_3^2-710 x_3^3\biggr)+e^{12 t_{\mu_A}} \left(-\frac{16552}{5} t_{2,2} t_{3,4}-7896 t_{3,4} x_2-4452 t_{2,2} x_3+\frac{9504}{5} x_2 x_3\right)\\
&+e^{13 t_{\mu_A}} \left(-\frac{65962 t_{1,1}}{5}-\frac{114756 x_1}{5}\right)+e^{16 t_{\mu_A}} \left(\frac{660536 t_{3,3}}{5}+\frac{3543388 t_{3,4}^2}{25}+\frac{1320902}{5} t_{3,4} x_3+161904 x_3^2\right)\\
&+e^{18 t_{\mu_A}} \left(\frac{2245908 t_{2,2}}{5}+705888 x_2\right)+e^{22 t_{\mu_A}} (-23681508 t_{3,4}-15860810 x_3)+1067265624 e^{28 t_{\mu_A}},
\end{align*}

\begin{align*}
&\phi^3_{88}=\frac{1}{25} t_{3,3} x_3^2+\frac{1}{250} t_{3,4}^2 x_3^2+\frac{2}{25} t_{3,4} x_3^3+\frac{x_3^4}{5}+e^{2 t_{\mu_A}} \left(-\frac{4}{25} t_{2,2} t_{3,3}-\frac{2}{125} t_{2,2} t_{3,4}^2-\frac{6}{25} t_{3,3} x_2-\frac{3}{125} t_{3,4}^2 x_2\right.\\
&\left.-\frac{8}{25} t_{2,2} t_{3,4} x_3-t_{2,2} x_3^2\right)+e^{3 t_{\mu_A}} \left(-\frac{24}{25} t_{1,1} t_{3,4}-\frac{24}{25} t_{3,4} x_1-\frac{12}{5} t_{1,1} x_3\right)+e^{4 t_{\mu_A}} \left(-\frac{24 t_{2,1}}{5}-\frac{18}{5} t_{2,2} x_2-\frac{36 x_2^2}{5}\right)\\
&+e^{6 t_{\mu_A}} \left(\frac{36 t_{3,2}}{5}+\frac{264}{25} t_{3,3} t_{3,4}+\frac{48 t_{3,4}^3}{25}+12 t_{3,3} x_3+\frac{252}{25} t_{3,4}^2 x_3+\frac{53}{5} t_{3,4} x_3^2-26 x_3^3\right)+e^{8 t_{\mu_A}} \left(\frac{1316}{25} t_{2,2} t_{3,4}\right.\\
&\left.+\frac{762}{5} t_{3,4} x_2+\frac{628}{5} t_{2,2} x_3\right)+e^{9 t_{\mu_A}} \left(\frac{1281 t_{1,1}}{5}+258 x_1\right)+e^{12 t_{\mu_A}} \left(-\frac{49104 t_{3,3}}{25}-\frac{240552 t_{3,4}^2}{125}-\frac{83208}{25} t_{3,4} x_3-6 x_3^2\right)\\
&+e^{14 t_{\mu_A}} \left(-\frac{44412 t_{2,2}}{5}-\frac{74268 x_2}{5}\right)+e^{18 t_{\mu_A}} \left(\frac{1815468 t_{3,4}}{5}+173128 x_3\right)-\frac{83805096}{5} e^{24 t_{\mu_A}}.
\end{align*}

\begin{align*}
&\phi^1_{89}=-\frac{1}{5} t_{3,2} x_1 x_3-\frac{1}{25} t_{3,3} t_{3,4} x_1 x_3-\frac{1}{750} t_{3,4}^3 x_1 x_3-\frac{2}{5} t_{3,3} x_1 x_3^2-\frac{2}{25} t_{3,4}^2 x_1 x_3^2-\frac{3}{5} t_{3,4} x_1 x_3^3-x_1 x_3^4\\
&+e^{2 t_{\mu_A}} \left(-\frac{4}{5} t_{1,1} t_{2,2} t_{3,3}-\frac{4}{25} t_{1,1} t_{2,2} t_{3,4}^2-\frac{6}{5} t_{1,1} t_{3,3} x_2-\frac{6}{25} t_{1,1} t_{3,4}^2 x_2-\frac{6}{5} t_{3,3} x_1 x_2-\frac{6}{25} t_{3,4}^2 x_1 x_2-\frac{6}{5} t_{1,1} t_{2,2} t_{3,4} x_3\right.\\
&\left.-\frac{9}{5} t_{2,2} t_{3,4} x_1 x_3-\frac{9}{5} t_{1,1} t_{3,4} x_2 x_3-\frac{36}{5} t_{3,4} x_1 x_2 x_3-2 t_{1,1} t_{2,2} x_3^2-2 t_{2,2} x_1 x_3^2-3 t_{1,1} x_2 x_3^2-12 x_1 x_2 x_3^2\right)+e^{3 t_{\mu_A}} \left(-\frac{6}{5} t_{3,2}^2\right.\\
&-\frac{9}{5} t_{1,1}^2 t_{3,4}-\frac{18}{25} t_{3,2} t_{3,3} t_{3,4}-\frac{12}{125} t_{3,3}^2 t_{3,4}^2-\frac{1}{25} t_{3,2} t_{3,4}^3-\frac{6}{625} t_{3,3} t_{3,4}^4-\frac{2 t_{3,4}^6}{9375}+\frac{24}{5} t_{2,2}^2 t_{3,4} x_2+\frac{72}{5} t_{2,2} t_{3,4} x_2^2+\frac{54}{5} t_{3,4} x_2^3\\
&-3 t_{1,1}^2 x_3-\frac{6}{5} t_{3,2} t_{3,3} x_3+\frac{6}{25} t_{3,3}^2 t_{3,4} x_3-\frac{3}{5} t_{3,2} t_{3,4}^2 x_3+\frac{2}{125} t_{3,3} t_{3,4}^3 x_3+\frac{7 t_{3,4}^5 x_3}{1250}-4 t_{1,1} x_1 x_3+8 t_{2,2}^2 x_2 x_3+24 t_{2,2} x_2^2 x_3\\
&\left.+18 x_2^3 x_3-\frac{6}{5} t_{3,2} t_{3,4} x_3^2+\frac{12}{25} t_{3,3} t_{3,4}^2 x_3^2+\frac{8}{125} t_{3,4}^4 x_3^2+\frac{6}{5} t_{3,3} t_{3,4} x_3^3+\frac{3}{25} t_{3,4}^3 x_3^3\right)+e^{4 t_{\mu_A}} \left(-4 t_{1,1} t_{2,2}^2+6 t_{1,1} t_{2,2} x_2\right.\\
&\left.+3 t_{2,2} x_1 x_2+18 t_{1,1} x_2^2-9 x_1 x_2^2\right)+e^{5 t_{\mu_A}} \biggl(4 t_{2,2} t_{3,3}^2+\frac{18}{5} t_{2,2} t_{3,2} t_{3,4}+\frac{74}{25} t_{2,2} t_{3,3} t_{3,4}^2+\frac{36}{125} t_{2,2} t_{3,4}^4+\frac{84}{5} t_{3,3}^2 x_2\\
&+\frac{246}{25} t_{3,3} t_{3,4}^2 x_2+\frac{144}{125} t_{3,4}^4 x_2-6 t_{2,2} t_{3,2} x_3+\frac{144}{5} t_{2,2} t_{3,3} t_{3,4} x_3+\frac{236}{25} t_{2,2} t_{3,4}^3 x_3+\frac{288}{5} t_{3,3} t_{3,4} x_2 x_3+\frac{72}{5} t_{3,4}^3 x_2 x_3\\
&+28 t_{2,2} t_{3,3} x_3^2+\frac{216}{5} t_{2,2} t_{3,4}^2 x_3^2+60 t_{3,3} x_2 x_3^2+\frac{288}{5} t_{3,4}^2 x_2 x_3^2+54 t_{2,2} t_{3,4} x_3^3+72 t_{3,4} x_2 x_3^3\biggr)+e^{6 t_{\mu_A}} \Biggl(6 t_{1,1} t_{3,2}\\
&+\frac{288}{5} t_{1,1} t_{3,3} t_{3,4}+\frac{412}{25} t_{1,1} t_{3,4}^3+\frac{288}{5} t_{3,3} t_{3,4} x_1+\frac{432}{25} t_{3,4}^3 x_1+105 t_{1,1} t_{3,3} x_3+\frac{1023}{10} t_{1,1} t_{3,4}^2 x_3+161 t_{3,3} x_1 x_3\\
&\left.+\frac{1969}{10} t_{3,4}^2 x_1 x_3+210 t_{1,1} t_{3,4} x_3^2+442 t_{3,4} x_1 x_3^2+135 t_{1,1} x_3^3+375 x_1 x_3^3\right)+e^{7 t_{\mu_A}} \left(\frac{226}{5} t_{2,2}^2 t_{3,3}-\frac{1089}{25} t_{2,2}^2 t_{3,4}^2\right.\\
&+\frac{264}{5} t_{2,2} t_{3,3} x_2-\frac{8226}{25} t_{2,2} t_{3,4}^2 x_2+72 t_{3,3} x_2^2-\frac{1944}{5} t_{3,4}^2 x_2^2-\frac{6}{5} t_{2,2}^2 t_{3,4} x_3-\frac{2304}{5} t_{2,2} t_{3,4} x_2 x_3-648 t_{3,4} x_2^2 x_3+198 t_{2,2}^2 x_3^2\\
&+252 t_{2,2} x_2 x_3^2\Biggr)+e^{8 t_{\mu_A}} \left(\frac{1432}{5} t_{1,1} t_{2,2} t_{3,4}+288 t_{2,2} t_{3,4} x_1+138 t_{1,1} t_{3,4} x_2+1065 t_{3,4} x_1 x_2+639 t_{1,1} t_{2,2} x_3\right.\\
&+1026 t_{2,2} x_1 x_3+1125 t_{1,1} x_2 x_3+4500 x_1 x_2 x_3\Biggr)+e^{9 t_{\mu_A}} \Biggl(105 t_{1,1}^2-144 t_{2,2}^3+282 t_{3,2} t_{3,3}-312 t_{3,3}^2 t_{3,4}+165 t_{3,2} t_{3,4}^2\\
&-\frac{884}{5} t_{3,3} t_{3,4}^3-\frac{3722 t_{3,4}^5}{125}-810 t_{2,1} x_2-1749 t_{2,2}^2 x_2-5112 t_{2,2} x_2^2-4320 x_2^3-810 t_{3,1} x_3-543 t_{3,3}^2 x_3+336 t_{3,2} t_{3,4} x_3\\
&-\frac{6867}{5} t_{3,3} t_{3,4}^2 x_3-\frac{7157}{20} t_{3,4}^4 x_3-1410 t_{3,2} x_3^2-2592 t_{3,3} t_{3,4} x_3^2-\frac{7064}{5} t_{3,4}^3 x_3^2-3360 t_{3,3} x_3^3-2394 t_{3,4}^2 x_3^3-3240 t_{3,4} x_3^4\\
&-4050 x_3^5\biggr)+e^{11 t_{\mu_A}} \left(-118 t_{2,2} t_{3,2}-\frac{26846}{5} t_{2,2} t_{3,3} t_{3,4}+\frac{31369}{75} t_{2,2} t_{3,4}^3+2118 t_{3,2} x_2-\frac{50232}{5} t_{3,3} t_{3,4} x_2+\frac{50168}{25} t_{3,4}^3 x_2\right.\\
&-7068 t_{2,2} t_{3,3} x_3-\frac{46044}{5} t_{2,2} t_{3,4}^2 x_3-16704 t_{3,3} x_2 x_3-\frac{37944}{5} t_{3,4}^2 x_2 x_3-22822 t_{2,2} t_{3,4} x_3^2-29256 t_{3,4} x_2 x_3^2+740 t_{2,2} x_3^3\\
&-11760 x_2 x_3^3\biggr)+e^{12 t_{\mu_A}} \left(-\frac{50424}{5} t_{1,1} t_{3,3}-\frac{369624}{25} t_{1,1} t_{3,4}^2-12000 t_{3,3} x_1-25200 t_{3,4}^2 x_1-\frac{242061}{5} t_{1,1} t_{3,4} x_3\right.\\
&\left.-\frac{550014}{5} t_{3,4} x_1 x_3-23112 t_{1,1} x_3^2-77512 x_1 x_3^2\right)+e^{13 t_{\mu_A}} \left(9720 t_{2,1} t_{3,4}+\frac{88006}{5} t_{2,2}^2 t_{3,4}+\frac{660918}{5} t_{2,2} t_{3,4} x_2\right.\\
&\left.+\frac{849672}{5} t_{3,4} x_2^2+12150 t_{2,1} x_3-15533 t_{2,2}^2 x_3+19536 t_{2,2} x_2 x_3+57024 x_2^2 x_3\right)+e^{14 t_{\mu_A}} (-33778 t_{1,1} t_{2,2}-60000 t_{2,2} x_1\\
&-30264 t_{1,1} x_2-243192 x_1 x_2)+e^{15 t_{\mu_A}} \biggl(48600 t_{3,1}+56700 t_{3,3}^2-\frac{353112}{5} t_{3,2} t_{3,4}+\frac{4588314}{25} t_{3,3} t_{3,4}^2+\frac{4335086 t_{3,4}^4}{125}\\
&\left.+115494 t_{3,2} x_3+\frac{2899374}{5} t_{3,3} t_{3,4} x_3+\frac{8402551}{25} t_{3,4}^3 x_3+634248 t_{3,3} x_3^2+\frac{4278966}{5} t_{3,4}^2 x_3^2+1039614 t_{3,4} x_3^3+1053000 x_3^4\right)
\end{align*}

\begin{align*}
&+e^{17 t_{\mu_A}} \left(\frac{2736506}{5} t_{2,2} t_{3,3}-\frac{2065049}{25} t_{2,2} t_{3,4}^2+1223070 t_{3,3} x_2-\frac{5167743}{5} t_{3,4}^2 x_2+\frac{11899794}{5} t_{2,2} t_{3,4} x_3\right.\\
&+2418840 t_{3,4} x_2 x_3+873698 t_{2,2} x_3^2+2627220 x_2 x_3^2\Biggr)+e^{18 t_{\mu_A}} (2880264 t_{1,1} t_{3,4}+6638400 t_{3,4} x_1+3090885 t_{1,1} x_3\\
&+9691725 x_1 x_3)+e^{19 t_{\mu_A}} \left(-1117800 t_{2,1}-688278 t_{2,2}^2-7162794 t_{2,2} x_2-9575280 x_2^2\right)+e^{21 t_{\mu_A}} \Biggl(584970 t_{3,2}\\
&-38455356 t_{3,3} t_{3,4}-\frac{81542494 t_{3,4}^3}{5}-57636180 t_{3,3} x_3-102794646 t_{3,4}^2 x_3-158851140 t_{3,4} x_3^2-131244600 x_3^3\biggr)\\
&+e^{23 t_{\mu_A}} \left(-28927128 t_{2,2} t_{3,4}+\frac{253873764}{5} t_{3,4} x_2-113054640 t_{2,2} x_3-177923952 x_2 x_3\right)+e^{24 t_{\mu_A}} (-155768424 t_{1,1}\\
&-445500000 x_1)+e^{27 t_{\mu_A}} \left(2424423000 t_{3,3}+\frac{18908795082 t_{3,4}^2}{5}+12444036816 t_{3,4} x_3+10511125560 x_3^2\right)\\
&+e^{29 t_{\mu_A}} (2803944506 t_{2,2}+2255008080 x_2)+e^{33 t_{\mu_A}} (-388000760040 t_{3,4}-545906954550 x_3)+14162630484960 e^{39 t_{\mu_A}},\\
&\phi^2_{89}=e^{t_{\mu_A}} \left(-\frac{1}{5} t_{1,1} t_{3,2}-\frac{1}{25} t_{1,1} t_{3,3} t_{3,4}-\frac{1}{750} t_{1,1} t_{3,4}^3-\frac{2}{5} t_{3,2} x_1-\frac{2}{25} t_{3,3} t_{3,4} x_1-\frac{1}{375} t_{3,4}^3 x_1-\frac{2}{5} t_{1,1} t_{3,3} x_3\right.\\
&\left.-\frac{2}{25} t_{1,1} t_{3,4}^2 x_3-\frac{4}{5} t_{3,3} x_1 x_3-\frac{4}{25} t_{3,4}^2 x_1 x_3-\frac{3}{5} t_{1,1} t_{3,4} x_3^2-\frac{6}{5} t_{3,4} x_1 x_3^2-t_{1,1} x_3^3-2 x_1 x_3^3\right)+e^{2 t_{\mu_A}} \left(-\frac{4}{5} t_{2,1} t_{3,3}\right.\\
&-\frac{2}{15} t_{2,2}^2 t_{3,3}-\frac{4}{25} t_{2,1} t_{3,4}^2-\frac{2}{75} t_{2,2}^2 t_{3,4}^2-\frac{6}{5} t_{2,1} t_{3,4} x_3-\frac{1}{5} t_{2,2}^2 t_{3,4} x_3+\frac{12}{5} t_{2,2} t_{3,4} x_2 x_3+\frac{18}{5} t_{3,4} x_2^2 x_3-2 t_{2,1} x_3^2-\frac{1}{3} t_{2,2}^2 x_3^2\\
&+4 t_{2,2} x_2 x_3^2+6 x_2^2 x_3^2\biggr)+e^{3 t_{\mu_A}} \left(-\frac{21}{5} t_{1,1} t_{2,2} t_{3,4}-\frac{18}{5} t_{2,2} t_{3,4} x_1-\frac{18}{5} t_{1,1} t_{3,4} x_2-\frac{48}{5} t_{3,4} x_1 x_2-6 t_{1,1} t_{2,2} x_3\right.\\
&-4 t_{2,2} x_1 x_3-48 x_1 x_2 x_3\biggr)+e^{4 t_{\mu_A}} \biggl(-10 t_{1,1}^2-4 t_{2,1} t_{2,2}-\frac{2 t_{2,2}^3}{3}+\frac{8}{5} t_{3,2} t_{3,3}+\frac{4}{5} t_{3,3}^2 t_{3,4}+\frac{6}{25} t_{3,2} t_{3,4}^2+\frac{76}{375} t_{3,3} t_{3,4}^3\\
&+\frac{7 t_{3,4}^5}{625}-8 t_{1,1} x_1+24 t_{3,1} x_3+8 t_{3,3}^2 x_3+6 t_{3,2} t_{3,4} x_3+\frac{142}{25} t_{3,3} t_{3,4}^2 x_3+\frac{71}{125} t_{3,4}^4 x_3+54 t_{3,2} x_3^2+\frac{168}{5} t_{3,3} t_{3,4} x_3^2\\
&\left.+\frac{138}{25} t_{3,4}^3 x_3^2+92 t_{3,3} x_3^3+\frac{144}{5} t_{3,4}^2 x_3^3+96 t_{3,4} x_3^4+120 x_3^5\right)+e^{6 t_{\mu_A}} \biggl(\frac{204}{5} t_{2,2} t_{3,2}+\frac{1164}{25} t_{2,2} t_{3,3} t_{3,4}+\frac{1084}{125} t_{2,2} t_{3,4}^3\\
&+54 t_{3,2} x_2+\frac{468}{5} t_{3,3} t_{3,4} x_2+\frac{468}{25} t_{3,4}^3 x_2+\frac{78}{5} t_{2,2} t_{3,3} x_3+\frac{1113}{25} t_{2,2} t_{3,4}^2 x_3+72 t_{3,3} x_2 x_3-\frac{252}{5} t_{3,4}^2 x_2 x_3+\frac{42}{5} t_{2,2} t_{3,4} x_3^2\\
&-108 t_{3,4} x_2 x_3^2-186 t_{2,2} x_3^3\biggr)+e^{7 t_{\mu_A}} \left(\frac{1391}{5} t_{1,1} t_{3,3}+\frac{12587}{50} t_{1,1} t_{3,4}^2+322 t_{3,3} x_1+\frac{1969}{5} t_{3,4}^2 x_1+\frac{2774}{5} t_{1,1} t_{3,4} x_3\right.\\
&+884 t_{3,4} x_1 x_3+213 t_{1,1} x_3^2+750 x_1 x_3^2\biggr)+e^{8 t_{\mu_A}} \left(596 t_{2,1} t_{3,4}+\frac{482}{15} t_{2,2}^2 t_{3,4}+\frac{2592}{5} t_{2,2} t_{3,4} x_2+972 t_{3,4} x_2^2+390 t_{2,1} x_3\right.\\
&+59 t_{2,2}^2 x_3-600 t_{2,2} x_2 x_3-1440 x_2^2 x_3\biggr)+e^{9 t_{\mu_A}} (1416 t_{1,1} t_{2,2}+2052 t_{2,2} x_1+2520 t_{1,1} x_2+7820 x_1 x_2)\\
&+e^{10 t_{\mu_A}} \biggl(-1440 t_{3,1}-1648 t_{3,3}^2-1988 t_{3,2} t_{3,4}-\frac{16518}{5} t_{3,3} t_{3,4}^2-\frac{47047 t_{3,4}^4}{75}-6630 t_{3,2} x_3-9696 t_{3,3} t_{3,4} x_3\\
&\left.-3024 t_{3,4}^3 x_3-16950 t_{3,3} x_3^2-9465 t_{3,4}^2 x_3^2-24190 t_{3,4} x_3^3-31200 x_3^4\right)+e^{12 t_{\mu_A}} \left(-14532 t_{2,2} t_{3,3}-\frac{72024}{5} t_{2,2} t_{3,4}^2\right.\\
&\left.-37224 t_{3,3} x_2-\frac{178704}{5} t_{3,4}^2 x_2-20108 t_{2,2} t_{3,4} x_3+\frac{95112}{5} t_{3,4} x_2 x_3+12050 t_{2,2} x_3^2+9504 x_2 x_3^2\right)\\
&+e^{13 t_{\mu_A}} \left(-\frac{606656}{5} t_{1,1} t_{3,4}-\frac{1100028}{5} t_{3,4} x_1-82696 t_{1,1} x_3-155024 x_1 x_3\right)+e^{14 t_{\mu_A}} \left(-73424 t_{2,1}-\frac{74872 t_{2,2}^2}{3}\right.\\
&-166050 t_{2,2} x_2-213840 x_2^2\Biggr)+e^{16 t_{\mu_A}} \left.(378972 t_{3,2}+\frac{7113124}{5} t_{3,3} t_{3,4}+\frac{16080922 t_{3,4}^3}{25}+1685018 t_{3,3} x_3\right.\\
\end{align*}

\begin{align*}
&+1682971 t_{3,4}^2 x_3+3127632 t_{3,4} x_3^2+3654120 x_3^3\Biggr)+e^{18 t_{\mu_A}} \left(\frac{24792324}{5} t_{2,2} t_{3,4}+10311552 t_{3,4} x_2+2003564 t_{2,2} x_3\right.\\
&-554400 x_2 x_3\biggr)+e^{19 t_{\mu_A}} (10436151 t_{1,1}+19383450 x_1)+e^{22 t_{\mu_A}} \left(-119109720 t_{3,3}-202718352 t_{3,4}^2-297398170 t_{3,4} x_3\right.\\
&\left.-269136900 x_3^2\right)+e^{24 t_{\mu_A}} (-363099900 t_{2,2}-646842240 x_2)+e^{28 t_{\mu_A}} (23176371060 t_{3,4}+15870520800 x_3)\\
&-860105264700 e^{34 t_{\mu_A}},\\
&\phi^3_{89}=e^{2 t_{\mu_A}} \left(-\frac{4}{5} t_{2,2} t_{3,2}-\frac{4}{25} t_{2,2} t_{3,3} t_{3,4}-\frac{2}{375} t_{2,2} t_{3,4}^3-\frac{6}{5} t_{3,2} x_2-\frac{6}{25} t_{3,3} t_{3,4} x_2-\frac{1}{125} t_{3,4}^3 x_2-\frac{4}{5} t_{2,2} t_{3,3} x_3\right.\\
&\left.-\frac{4}{25} t_{2,2} t_{3,4}^2 x_3-\frac{6}{5} t_{2,2} t_{3,4} x_3^2-2 t_{2,2} x_3^3\right)+e^{3 t_{\mu_A}} \left(-\frac{18}{5} t_{1,1} t_{3,3}-\frac{18}{25} t_{1,1} t_{3,4}^2-\frac{12}{5} t_{3,3} x_1-\frac{12}{25} t_{3,4}^2 x_1-\frac{18}{5} t_{1,1} t_{3,4} x_3\right.\\
&-6 t_{1,1} x_3^2+24 x_1 x_3^2\biggr)+e^{4 t_{\mu_A}} \left(-\frac{48}{5} t_{2,1} t_{3,4}+\frac{4}{5} t_{2,2}^2 t_{3,4}-\frac{6}{5} t_{2,2} t_{3,4} x_2-\frac{36}{5} t_{3,4} x_2^2-36 t_{2,1} x_3-48 t_{2,2} x_2 x_3\right.\\
&-72 x_2^2 x_3\biggr)+e^{5 t_{\mu_A}} (-t_{1,1} t_{2,2}+6 t_{2,2} x_1+3 t_{1,1} x_2-18 x_1 x_2)+e^{6 t_{\mu_A}} \left(-18 t_{3,1}+\frac{63 t_{3,3}^2}{5}+18 t_{3,2} t_{3,4}+\frac{351}{25} t_{3,3} t_{3,4}^2\right.\\
&\left.+\frac{129 t_{3,4}^4}{100}+36 t_{3,2} x_3+\frac{324}{5} t_{3,3} t_{3,4} x_3+12 t_{3,4}^3 x_3-84 t_{3,3} x_3^2+12 t_{3,4}^2 x_3^2-138 t_{3,4} x_3^3-300 x_3^4\right)+e^{8 t_{\mu_A}} \left(\frac{1522}{5} t_{2,2} t_{3,3}\right.\\
&\left.+\frac{1447}{25} t_{2,2} t_{3,4}^2+534 t_{3,3} x_2+231 t_{3,4}^2 x_2+\frac{5152}{5} t_{2,2} t_{3,4} x_3+1296 t_{3,4} x_2 x_3+1614 t_{2,2} x_3^2-810 x_2 x_3^2\right)\\
&+e^{9 t_{\mu_A}} \left(\frac{7857}{5} t_{1,1} t_{3,4}+2130 t_{3,4} x_1+4824 t_{1,1} x_3+180 x_1 x_3\right)+e^{10 t_{\mu_A}} \left(6470 t_{2,1}-\frac{23 t_{2,2}^2}{3}+8548 t_{2,2} x_2+14820 x_2^2\right)\\
&+e^{12 t_{\mu_A}} \left(-\frac{68424 t_{3,2}}{5}-\frac{568974}{25} t_{3,3} t_{3,4}-\frac{660829 t_{3,4}^3}{125}-\frac{119124}{5} t_{3,3} x_3-\frac{679824}{25} t_{3,4}^2 x_3-\frac{36186}{5} t_{3,4} x_3^2+46488 x_3^3\right)\\
&+e^{14 t_{\mu_A}} \left(-\frac{588876}{5} t_{2,2} t_{3,4}-\frac{1659996}{5} t_{3,4} x_2-249308 t_{2,2} x_3-65448 x_2 x_3\right)+e^{15 t_{\mu_A}} (-433896 t_{1,1}-486384 x_1)\\
&+e^{18 t_{\mu_A}} \left(3490344 t_{3,3}+\frac{19796394 t_{3,4}^2}{5}+6493320 t_{3,4} x_3-1275720 x_3^2\right)+e^{20 t_{\mu_A}} (16938594 t_{2,2}+28878030 x_2)\\
&+e^{24 t_{\mu_A}} \left(-\frac{3397991478 t_{3,4}}{5}-292688964 x_3\right)+29981411520 e^{30 t_{\mu_A}}.
\end{align*}

\begin{align*}
&\phi^1_{99}=-e^{-t_{\mu_A}} x_1^2 x_2 x_3+e^{2 t_{\mu_A}} \left(-4 t_{1,1} t_{2,2} t_{3,2}-\frac{8}{5} t_{1,1} t_{2,2} t_{3,3} t_{3,4}-\frac{8}{75} t_{1,1} t_{2,2} t_{3,4}^3-6 t_{1,1} t_{3,2} x_2-\frac{12}{5} t_{1,1} t_{3,3} t_{3,4} x_2\right.\\
&-\frac{4}{25} t_{1,1} t_{3,4}^3 x_2-12 t_{3,2} x_1 x_2-\frac{24}{5} t_{3,3} t_{3,4} x_1 x_2-\frac{8}{25} t_{3,4}^3 x_1 x_2-6 t_{1,1} t_{2,2} t_{3,3} x_3-\frac{9}{5} t_{1,1} t_{2,2} t_{3,4}^2 x_3-12 t_{2,2} t_{3,3} x_1 x_3\\
&-\frac{18}{5} t_{2,2} t_{3,4}^2 x_1 x_3-9 t_{1,1} t_{3,3} x_2 x_3-\frac{27}{10} t_{1,1} t_{3,4}^2 x_2 x_3-54 t_{3,3} x_1 x_2 x_3-\frac{81}{5} t_{3,4}^2 x_1 x_2 x_3-8 t_{1,1} t_{2,2} t_{3,4} x_3^2-16 t_{2,2} t_{3,4} x_1 x_3^2\\
&-12 t_{1,1} t_{3,4} x_2 x_3^2-72 t_{3,4} x_1 x_2 x_3^2-10 t_{1,1} t_{2,2} x_3^3-20 t_{2,2} x_1 x_3^3-15 t_{1,1} x_2 x_3^3-90 x_1 x_2 x_3^3\biggr)+e^{3 t_{\mu_A}} \Biggl(-9 t_{1,1}^2 t_{3,3}\\
&-\frac{27}{10} t_{1,1}^2 t_{3,4}^2+24 t_{2,2}^2 t_{3,3} x_2+\frac{36}{5} t_{2,2}^2 t_{3,4}^2 x_2+72 t_{2,2} t_{3,3} x_2^2+\frac{108}{5} t_{2,2} t_{3,4}^2 x_2^2+54 t_{3,3} x_2^3+\frac{81}{5} t_{3,4}^2 x_2^3+24 t_{3,2}^2 x_3\\
&-12 t_{1,1}^2 t_{3,4} x_3+\frac{96}{5} t_{3,2} t_{3,3} t_{3,4} x_3+\frac{96}{25} t_{3,3}^2 t_{3,4}^2 x_3+\frac{32}{25} t_{3,2} t_{3,4}^3 x_3+\frac{64}{125} t_{3,3} t_{3,4}^4 x_3+\frac{32 t_{3,4}^6 x_3}{1875}-24 t_{1,1} t_{3,4} x_1 x_3\\
&+32 t_{2,2}^2 t_{3,4} x_2 x_3+96 t_{2,2} t_{3,4} x_2^2 x_3+72 t_{3,4} x_2^3 x_3-15 t_{1,1}^2 x_3^2+72 t_{3,2} t_{3,3} x_3^2+\frac{144}{5} t_{3,3}^2 t_{3,4} x_3^2+\frac{108}{5} t_{3,2} t_{3,4}^2 x_3^2\\
&+\frac{264}{25} t_{3,3} t_{3,4}^3 x_3^2+\frac{72}{125} t_{3,4}^5 x_3^2-30 t_{1,1} x_1 x_3^2+40 t_{2,2}^2 x_2 x_3^2+120 t_{2,2} x_2^2 x_3^2+90 x_2^3 x_3^2+54 t_{3,3}^2 x_3^3+96 t_{3,2} t_{3,4} x_3^3
\end{align*}

\begin{align*}
&+\frac{354}{5} t_{3,3} t_{3,4}^2 x_3^3+\frac{371}{50} t_{3,4}^4 x_3^3+120 t_{3,2} x_3^4+192 t_{3,3} t_{3,4} x_3^4+\frac{232}{5} t_{3,4}^3 x_3^4+180 t_{3,3} x_3^5+150 t_{3,4}^2 x_3^5+240 t_{3,4} x_3^6+150 x_3^7\Biggr)\\
&+e^{4 t_{\mu_A}} \left(-16 t_{1,1} t_{2,2}^2 t_{3,4}+24 t_{1,1} t_{2,2} t_{3,4} x_2+24 t_{2,2} t_{3,4} x_1 x_2+72 t_{1,1} t_{3,4} x_2^2-72 t_{3,4} x_1 x_2^2-10 t_{1,1} t_{2,2}^2 x_3-20 t_{2,2}^2 x_1 x_3\right.\\
&\left.+45 t_{1,1} t_{2,2} x_2 x_3-90 t_{2,2} x_1 x_2 x_3+90 t_{1,1} x_2^2 x_3-360 x_1 x_2^2 x_3\right)+e^{5 t_{\mu_A}} \left(-35 t_{1,1}^2 t_{2,2}+8 t_{2,2} t_{3,2} t_{3,3}+\frac{16}{5} t_{2,2} t_{3,3}^2 t_{3,4}\right.\\
&+\frac{12}{5} t_{2,2} t_{3,2} t_{3,4}^2+\frac{88}{75} t_{2,2} t_{3,3} t_{3,4}^3+\frac{8}{125} t_{2,2} t_{3,4}^5-30 t_{1,1}^2 x_2+40 t_{2,2}^3 x_2+120 t_{3,2} t_{3,3} x_2+48 t_{3,3}^2 t_{3,4} x_2+36 t_{3,2} t_{3,4}^2 x_2\\
&+\frac{88}{5} t_{3,3} t_{3,4}^3 x_2+\frac{24}{25} t_{3,4}^5 x_2+120 t_{1,1} x_1 x_2+120 t_{2,2}^2 x_2^2+90 t_{2,2} x_2^3+120 t_{2,2} t_{3,3}^2 x_3-96 t_{2,2} t_{3,2} t_{3,4} x_3+\frac{168}{5} t_{2,2} t_{3,3} t_{3,4}^2 x_3\\
&+\frac{206}{25} t_{2,2} t_{3,4}^4 x_3+180 t_{3,3}^2 x_2 x_3+288 t_{3,2} t_{3,4} x_2 x_3+\frac{1116}{5} t_{3,3} t_{3,4}^2 x_2 x_3+\frac{597}{25} t_{3,4}^4 x_2 x_3-240 t_{2,2} t_{3,2} x_3^2+64 t_{2,2} t_{3,3} t_{3,4} x_3^2\\
&+\frac{208}{5} t_{2,2} t_{3,4}^3 x_3^2+360 t_{3,2} x_2 x_3^2+816 t_{3,3} t_{3,4} x_2 x_3^2+\frac{1056}{5} t_{3,4}^3 x_2 x_3^2+20 t_{2,2} t_{3,3} x_3^3+6 t_{2,2} t_{3,4}^2 x_3^3+840 t_{3,3} x_2 x_3^3\\
&+828 t_{3,4}^2 x_2 x_3^3-240 t_{2,2} t_{3,4} x_3^4+1440 t_{3,4} x_2 x_3^4-300 t_{2,2} x_3^5+900 x_2 x_3^5\biggr)+e^{6 t_{\mu_A}} \biggl(252 t_{1,1} t_{3,3}^2-36 t_{1,1} t_{3,2} t_{3,4}\\
&+\frac{684}{5} t_{1,1} t_{3,3} t_{3,4}^2+\frac{543}{25} t_{1,1} t_{3,4}^4+288 t_{3,3}^2 x_1+\frac{864}{5} t_{3,3} t_{3,4}^2 x_1+\frac{648}{25} t_{3,4}^4 x_1-480 t_{1,1} t_{3,2} x_3+522 t_{1,1} t_{3,3} t_{3,4} x_3\\
&+\frac{1007}{5} t_{1,1} t_{3,4}^3 x_3-120 t_{3,2} x_1 x_3+1356 t_{3,3} t_{3,4} x_1 x_3+418 t_{3,4}^3 x_1 x_3+240 t_{1,1} t_{3,3} x_3^2+576 t_{1,1} t_{3,4}^2 x_3^2+1440 t_{3,3} x_1 x_3^2\\
&+1536 t_{3,4}^2 x_1 x_3^2+390 t_{1,1} t_{3,4} x_3^3+2340 t_{3,4} x_1 x_3^3-300 t_{1,1} x_3^4+1200 x_1 x_3^4\biggr)+e^{7 t_{\mu_A}} \left(576 t_{2,2}^2 t_{3,2}-\frac{368}{5} t_{2,2}^2 t_{3,3} t_{3,4}\right.\\
&-\frac{1896}{25} t_{2,2}^2 t_{3,4}^3+1944 t_{2,2} t_{3,2} x_2-\frac{7872}{5} t_{2,2} t_{3,3} t_{3,4} x_2-\frac{16344}{25} t_{2,2} t_{3,4}^3 x_2+2160 t_{3,2} x_2^2-1656 t_{3,3} t_{3,4} x_2^2-\frac{3492}{5} t_{3,4}^3 x_2^2\\
&+464 t_{2,2}^2 t_{3,3} x_3-\frac{264}{5} t_{2,2}^2 t_{3,4}^2 x_3+216 t_{2,2} t_{3,3} x_2 x_3-\frac{14076}{5} t_{2,2} t_{3,4}^2 x_2 x_3-720 t_{3,3} x_2^2 x_3-2808 t_{3,4}^2 x_2^2 x_3+912 t_{2,2}^2 t_{3,4} x_3^2\\
&-3312 t_{2,2} t_{3,4} x_2 x_3^2-3240 t_{3,4} x_2^2 x_3^2+1590 t_{2,2}^2 x_3^3+360 t_{2,2} x_2 x_3^3\biggr)+e^{8 t_{\mu_A}} \left(2632 t_{1,1} t_{2,2} t_{3,3}+\frac{5068}{5} t_{1,1} t_{2,2} t_{3,4}^2\right.\\
&+2880 t_{2,2} t_{3,3} x_1+864 t_{2,2} t_{3,4}^2 x_1+4560 t_{1,1} t_{3,3} x_2-240 t_{1,1} t_{3,4}^2 x_2+9360 t_{3,3} x_1 x_2+3624 t_{3,4}^2 x_1 x_2+4956 t_{1,1} t_{2,2} t_{3,4} x_3\\
&+8496 t_{2,2} t_{3,4} x_1 x_3+4005 t_{1,1} t_{3,4} x_2 x_3+27270 t_{3,4} x_1 x_2 x_3+7520 t_{1,1} t_{2,2} x_3^2+10720 t_{2,2} x_1 x_3^2+8400 t_{1,1} x_2 x_3^2\\
&+32400 x_1 x_2 x_3^2\biggr)+e^{9 t_{\mu_A}} \biggl(-720 t_{3,2}^2-1440 t_{3,3}^3+1845 t_{1,1}^2 t_{3,4}-1536 t_{2,2}^3 t_{3,4}-360 t_{3,2} t_{3,3} t_{3,4}-\frac{6624}{5} t_{3,3}^2 t_{3,4}^2\\
&+\frac{132}{5} t_{3,2} t_{3,4}^3-\frac{9312}{25} t_{3,3} t_{3,4}^4-\frac{4708 t_{3,4}^6}{125}-13176 t_{2,2}^2 t_{3,4} x_2-25488 t_{2,2} t_{3,4} x_2^2-12150 t_{3,4} x_2^3+3600 t_{1,1}^2 x_3-2760 t_{2,2}^3 x_3\\
&-6480 t_{3,2} t_{3,3} x_3-9048 t_{3,3}^2 t_{3,4} x_3-3768 t_{3,2} t_{3,4}^2 x_3-\frac{27768}{5} t_{3,3} t_{3,4}^3 x_3-\frac{17038}{25} t_{3,4}^5 x_3+3600 t_{1,1} x_1 x_3-23180 t_{2,2}^2 x_2 x_3\\
&-54960 t_{2,2} x_2^2 x_3-39600 x_2^3 x_3-15360 t_{3,3}^2 x_3^2-15840 t_{3,2} t_{3,4} x_3^2-27360 t_{3,3} t_{3,4}^2 x_3^2-\frac{26736}{5} t_{3,4}^4 x_3^2-24000 t_{3,2} x_3^3\\
&-61980 t_{3,3} t_{3,4} x_3^3-24034 t_{3,4}^3 x_3^3-57600 t_{3,3} x_3^4-63600 t_{3,4}^2 x_3^4-93900 t_{3,4} x_3^5-60000 x_3^6\biggr)+e^{10 t_{\mu_A}} \left(2680 t_{1,1} t_{2,2}^2\right.\\
&\left.+7200 t_{2,2}^2 x_1+4890 t_{1,1} t_{2,2} x_2+37380 t_{2,2} x_1 x_2+3600 t_{1,1} x_2^2+67680 x_1 x_2^2\right)+e^{11 t_{\mu_A}} \biggl(-18048 t_{2,2} t_{3,3}^2\\
&-24016 t_{2,2} t_{3,2} t_{3,4}-\frac{33136}{5} t_{2,2} t_{3,3} t_{3,4}^2+\frac{140824}{75} t_{2,2} t_{3,4}^4-44640 t_{3,3}^2 x_2-76200 t_{3,2} t_{3,4} x_2-21864 t_{3,3} t_{3,4}^2 x_2\\
&+\frac{22852}{5} t_{3,4}^4 x_2+13920 t_{2,2} t_{3,2} x_3-76608 t_{2,2} t_{3,3} t_{3,4} x_3-\frac{34048}{5} t_{2,2} t_{3,4}^3 x_3-72000 t_{3,2} x_2 x_3-176040 t_{3,3} t_{3,4} x_2 x_3\\
&-12972 t_{3,4}^3 x_2 x_3-55960 t_{2,2} t_{3,3} x_3^2-60164 t_{2,2} t_{3,4}^2 x_3^2-261600 t_{3,3} x_2 x_3^2-187320 t_{3,4}^2 x_2 x_3^2-77860 t_{2,2} t_{3,4} x_3^3\\
&-475800 t_{3,4} x_2 x_3^3+22200 t_{2,2} x_3^4-360000 x_2 x_3^4\biggr)+e^{12 t_{\mu_A}} \left(10656 t_{1,1} t_{3,2}-\frac{689688}{5} t_{1,1} t_{3,3} t_{3,4}-\frac{803796}{25} t_{1,1} t_{3,4}^3\right.\\
&-216000 t_{3,3} t_{3,4} x_1-64800 t_{3,4}^3 x_1-60816 t_{1,1} t_{3,3} x_3-\frac{1042074}{5} t_{1,1} t_{3,4}^2 x_3-230472 t_{3,3} x_1 x_3-\frac{2470608}{5} t_{3,4}^2 x_1 x_3\\
\end{align*}

\begin{align*}
&-321888 t_{1,1} t_{3,4} x_3^2-910896 t_{3,4} x_1 x_3^2-21360 t_{1,1} x_3^3-438120 x_1 x_3^3\biggr)+e^{13 t_{\mu_A}} \left(-68024 t_{2,2}^2 t_{3,3}+\frac{556684}{5} t_{2,2}^2 t_{3,4}^2\right.\\
&-150792 t_{2,2} t_{3,3} x_2+\frac{3054972}{5} t_{2,2} t_{3,4}^2 x_2-160128 t_{3,3} x_2^2+\frac{2350008}{5} t_{3,4}^2 x_2^2+137048 t_{2,2}^2 t_{3,4} x_3+1611144 t_{2,2} t_{3,4} x_2 x_3\\
&\left.+1743696 t_{3,4} x_2^2 x_3-204240 t_{2,2}^2 x_3^2+321480 t_{2,2} x_2 x_3^2+285120 x_2^2 x_3^2\right)+e^{14 t_{\mu_A}} (-493452 t_{1,1} t_{2,2} t_{3,4}-1080000 t_{2,2} t_{3,4} x_1\\
&-646416 t_{1,1} t_{3,4} x_2-3380736 t_{3,4} x_1 x_2-749670 t_{1,1} t_{2,2} x_3-1260540 t_{2,2} x_1 x_3-665505 t_{1,1} x_2 x_3-3973590 x_1 x_2 x_3)\\
&+e^{15 t_{\mu_A}} \left(-193185 t_{1,1}^2+111840 t_{2,2}^3+159408 t_{3,2} t_{3,3}+\frac{5920416}{5} t_{3,3}^2 t_{3,4}+\frac{2103912}{5} t_{3,2} t_{3,4}^2+\frac{17101296}{25} t_{3,3} t_{3,4}^3\right.\\
&+\frac{7979808 t_{3,4}^5}{125}+861840 t_{2,2}^2 x_2+1547640 t_{2,2} x_2^2+614790 x_2^3+1490160 t_{3,3}^2 x_3+1577424 t_{3,2} t_{3,4} x_3+\frac{23950128}{5} t_{3,3} t_{3,4}^2 x_3\\
&+\frac{25157276}{25} t_{3,4}^4 x_3+2099160 t_{3,2} x_3^2+10896144 t_{3,3} t_{3,4} x_3^2+\frac{27179808}{5} t_{3,4}^3 x_3^2+8679420 t_{3,3} x_3^3+14474376 t_{3,4}^2 x_3^3\\
&+19799760 t_{3,4} x_3^4+11699700 x_3^5\biggr)+e^{17 t_{\mu_A}} \left(2543256 t_{2,2} t_{3,2}+\frac{42484072}{5} t_{2,2} t_{3,3} t_{3,4}-\frac{35982556}{25} t_{2,2} t_{3,4}^3+8075160 t_{3,2} x_2\right.\\
&+19675632 t_{3,3} t_{3,4} x_2-1700232 t_{3,4}^3 x_2+10344384 t_{2,2} t_{3,3} x_3+\frac{12508056}{5} t_{2,2} t_{3,4}^2 x_3+24713640 t_{3,3} x_2 x_3\\
&+1013472 t_{3,4}^2 x_2 x_3+18909632 t_{2,2} t_{3,4} x_3^2+54765120 t_{3,4} x_2 x_3^2+9743940 t_{2,2} x_3^3+56039400 x_2 x_3^3\biggr)\\
&+e^{18 t_{\mu_A}} \left(12174840 t_{1,1} t_{3,3}+18742176 t_{1,1} t_{3,4}^2+21384000 t_{3,3} x_1+46915200 t_{3,4}^2 x_1+43516650 t_{1,1} t_{3,4} x_3\right.\\
&\left.+117255780 t_{3,4} x_1 x_3+25045200 t_{1,1} x_3^2+86731200 x_1 x_3^2\right)+e^{19 t_{\mu_A}} \left(-9348032 t_{2,2}^2 t_{3,4}-60201816 t_{2,2} t_{3,4} x_2\right.\\
&\left.-34585560 t_{3,4} x_2^2-7344120 t_{2,2}^2 x_3-118365480 t_{2,2} x_2 x_3-130798800 x_2^2 x_3\right)+e^{20 t_{\mu_A}} (44237480 t_{1,1} t_{2,2}+106920000 t_{2,2} x_1\\
&+73369200 t_{1,1} x_2+322908000 x_1 x_2)+e^{21 t_{\mu_A}} \Biggl(-111676320 t_{3,3}^2-105057720 t_{3,2} t_{3,4}-372047472 t_{3,3} t_{3,4}^2\\
&-\frac{312471528 t_{3,4}^4}{5}-135330000 t_{3,2} x_3-1028764080 t_{3,3} t_{3,4} x_3
-520356832 t_{3,4}^3 x_3-938318400 t_{3,3} x_3^2-1821978240 t_{3,4}^2 x_3^2\\
&-2728752300 t_{3,4} x_3^3-1595040000 x_3^4\Biggr)+e^{23 t_{\mu_A}} \biggl(-818066800 t_{2,2} t_{3,3}-251287096 t_{2,2} t_{3,4}^2-1837633056 t_{3,3} x_2\\
&-\frac{4231559784}{5} t_{3,4}^2 x_2-790690320 t_{2,2} t_{3,4} x_3-989596008 t_{3,4} x_2 x_3-1433740200 t_{2,2} x_3^2-4342613760 x_2 x_3^2\biggr)\\
&+e^{24 t_{\mu_A}} (-3018152376 t_{1,1} t_{3,4}-8019000000 t_{3,4} x_1-2451801360 t_{1,1} x_3-7485858720 x_1 x_3)+e^{25 t_{\mu_A}} \left(-147207960 t_{2,2}^2\right.\\
&\left.+698942520 t_{2,2} x_2-208802880 x_2^2\right)+e^{27 t_{\mu_A}} \left(6205115520 t_{3,2}+57329810784 t_{3,3} t_{3,4}+29106573360 t_{3,4}^3\right.\\
&\left.+61873919280 t_{3,3} x_3+118547035032 t_{3,4}^2 x_3+236876184000 t_{3,4} x_3^2+156672016950 x_3^3\right)+e^{29 t_{\mu_A}} (86399501864 t_{2,2} t_{3,4}\\
&+193415433600 t_{3,4} x_2+60645688440 t_{2,2} x_3+100959556500 x_2 x_3)+e^{30 t_{\mu_A}} (144023219100 t_{1,1}+396940500000 x_1)\\
&+e^{33 t_{\mu_A}} \left(-2655547480800 t_{3,3}-5537714448240 t_{3,4}^2-11844161039400 t_{3,4} x_3-10402175520000 x_3^2\right)\\
&+e^{35 t_{\mu_A}} (-4835694794400 t_{2,2}-9566526202560 x_2)+e^{39 t_{\mu_A}} (452466384691680 t_{3,4}+424218905102400 x_3)\\
&-13317074584740000 e^{45 t_{\mu_A}},
\end{align*}

\begin{align*}
&\phi^2_{99}=-t_{1,1} x_1 x_2-2 x_1^2 x_2+e^{2 t_{\mu_A}} \left(-4 t_{2,1} t_{3,2}-\frac{2}{3} t_{2,2}^2 t_{3,2}-\frac{8}{5} t_{2,1} t_{3,3} t_{3,4}-\frac{4}{15} t_{2,2}^2 t_{3,3} t_{3,4}-\frac{8}{75} t_{2,1} t_{3,4}^3-\frac{4}{225} t_{2,2}^2 t_{3,4}^3\right.\\
&-6 t_{2,1} t_{3,3} x_3-t_{2,2}^2 t_{3,3} x_3-\frac{9}{5} t_{2,1} t_{3,4}^2 x_3-\frac{3}{10} t_{2,2}^2 t_{3,4}^2 x_3+12 t_{2,2} t_{3,3} x_2 x_3+\frac{18}{5} t_{2,2} t_{3,4}^2 x_2 x_3+18 t_{3,3} x_2^2 x_3+\frac{27}{5} t_{3,4}^2 x_2^2 x_3\\
&-8 t_{2,1} t_{3,4} x_3^2-\frac{4}{3} t_{2,2}^2 t_{3,4} x_3^2+16 t_{2,2} t_{3,4} x_2 x_3^2+24 t_{3,4} x_2^2 x_3^2-10 t_{2,1} x_3^3-\frac{5}{3} t_{2,2}^2 x_3^3+20 t_{2,2} x_2 x_3^3+30 x_2^2 x_3^3\biggr)\\
&+e^{3 t_{\mu_A}} \left(-24 t_{1,1} t_{2,2} t_{3,3}-\frac{36}{5} t_{1,1} t_{2,2} t_{3,4}^2-24 t_{2,2} t_{3,3} x_1-\frac{36}{5} t_{2,2} t_{3,4}^2 x_1-18 t_{1,1} t_{3,3} x_2-\frac{27}{5} t_{1,1} t_{3,4}^2 x_2-52 t_{3,3} x_1 x_2\right.
\end{align*}

\begin{align*}
&-\frac{78}{5} t_{3,4}^2 x_1 x_2-32 t_{1,1} t_{2,2} t_{3,4} x_3-32 t_{2,2} t_{3,4} x_1 x_3-144 t_{3,4} x_1 x_2 x_3-40 t_{1,1} t_{2,2} x_3^2-40 t_{2,2} x_1 x_3^2-180 x_1 x_2 x_3^2\biggr)\\
&+e^{4 t_{\mu_A}} \left(12 t_{3,2}^2-48 t_{1,1}^2 t_{3,4}-16 t_{2,1} t_{2,2} t_{3,4}-\frac{8}{3} t_{2,2}^3 t_{3,4}+\frac{48}{5} t_{3,2} t_{3,3} t_{3,4}+\frac{48}{25} t_{3,3}^2 t_{3,4}^2+\frac{16}{25} t_{3,2} t_{3,4}^3+\frac{32}{125} t_{3,3} t_{3,4}^4\right.\\
&+\frac{16 t_{3,4}^6}{1875}-48 t_{1,1} t_{3,4} x_1-60 t_{1,1}^2 x_3-10 t_{2,1} t_{2,2} x_3-\frac{5}{3} t_{2,2}^3 x_3+90 t_{3,2} t_{3,3} x_3+36 t_{3,3}^2 t_{3,4} x_3+27 t_{3,2} t_{3,4}^2 x_3+\frac{66}{5} t_{3,3} t_{3,4}^3 x_3\\
&+\frac{18}{25} t_{3,4}^5 x_3-60 t_{1,1} x_1 x_3+20 t_{2,2}^2 x_2 x_3+30 t_{2,2} x_2^2 x_3+108 t_{3,3}^2 x_3^2+120 t_{3,2} t_{3,4} x_3^2+\frac{564}{5} t_{3,3} t_{3,4}^2 x_3^2+\frac{323}{25} t_{3,4}^4 x_3^2\\
&\left.+150 t_{3,2} x_3^3+348 t_{3,3} t_{3,4} x_3^3+\frac{452}{5} t_{3,4}^3 x_3^3+360 t_{3,3} x_3^4+300 t_{3,4}^2 x_3^4+480 t_{3,4} x_3^5+300 x_3^6\right)+e^{5 t_{\mu_A}} \biggl(-20 t_{1,1} t_{2,1}\\
&\left.-\frac{130}{3} t_{1,1} t_{2,2}^2-40 t_{2,2}^2 x_1-30 t_{1,1} t_{2,2} x_2+84 t_{2,2} x_1 x_2\right)+e^{6 t_{\mu_A}} \left(48 t_{2,2} t_{3,3}^2+204 t_{2,2} t_{3,2} t_{3,4}+\frac{552}{5} t_{2,2} t_{3,3} t_{3,4}^2\right.\\
&+\frac{244}{25} t_{2,2} t_{3,4}^4+144 t_{3,3}^2 x_2+432 t_{3,2} t_{3,4} x_2+\frac{1296}{5} t_{3,3} t_{3,4}^2 x_2+\frac{612}{25} t_{3,4}^4 x_2-90 t_{2,2} t_{3,2} x_3+262 t_{2,2} t_{3,3} t_{3,4} x_3\\
&+87 t_{2,2} t_{3,4}^3 x_3+720 t_{3,2} x_2 x_3-36 t_{3,3} t_{3,4} x_2 x_3-78 t_{3,4}^3 x_2 x_3-280 t_{2,2} t_{3,3} x_3^2+228 t_{2,2} t_{3,4}^2 x_3^2-432 t_{3,4}^2 x_2 x_3^2\\
&-90 t_{2,2} t_{3,4} x_3^3-540 t_{3,4} x_2 x_3^3-600 t_{2,2} x_3^4\biggr)+e^{7 t_{\mu_A}} \left(36 t_{1,1} t_{3,2}+\frac{9432}{5} t_{1,1} t_{3,3} t_{3,4}+\frac{14064}{25} t_{1,1} t_{3,4}^3-240 t_{3,2} x_1\right.\\
&+2712 t_{3,3} t_{3,4} x_1+836 t_{3,4}^3 x_1+1584 t_{1,1} t_{3,3} x_3+\frac{12936}{5} t_{1,1} t_{3,4}^2 x_3+2880 t_{3,3} x_1 x_3+3072 t_{3,4}^2 x_1 x_3+4272 t_{1,1} t_{3,4} x_3^2\\
&+4680 t_{3,4} x_1 x_3^2+2040 t_{1,1} x_3^3+2400 x_1 x_3^3\biggr)+e^{8 t_{\mu_A}} \biggl(480 t_{2,1} t_{3,3}+976 t_{2,2}^2 t_{3,3}+2240 t_{2,1} t_{3,4}^2+\frac{1712}{15} t_{2,2}^2 t_{3,4}^2\\
&+672 t_{2,2} t_{3,3} x_2+\frac{9648}{5} t_{2,2} t_{3,4}^2 x_2+3888 t_{3,4}^2 x_2^2+5670 t_{2,1} t_{3,4} x_3+663 t_{2,2}^2 t_{3,4} x_3-4680 t_{2,2} t_{3,4} x_2 x_3-4050 t_{3,4} x_2^2 x_3\\
&+4400 t_{2,1} x_3^2+\frac{5980}{3} t_{2,2}^2 x_3^2-8800 t_{2,2} x_2 x_3^2-13200 x_2^2 x_3^2\biggr)+e^{9 t_{\mu_A}} (11028 t_{1,1} t_{2,2} t_{3,4}+16992 t_{2,2} t_{3,4} x_1+12690 t_{1,1} t_{3,4} x_2\\
&+35540 t_{3,4} x_1 x_2+15920 t_{1,1} t_{2,2} x_3+21440 t_{2,2} x_1 x_3+14400 t_{1,1} x_2 x_3+64800 x_1 x_2 x_3)+e^{10 t_{\mu_A}} \Biggl(6480 t_{1,1}^2+3580 t_{2,1} t_{2,2}\\
&+\frac{1070 t_{2,2}^3}{3}-2820 t_{3,2} t_{3,3}-12312 t_{3,3}^2 t_{3,4}-9822 t_{3,2} t_{3,4}^2-\frac{53572}{5} t_{3,3} t_{3,4}^3-\frac{31712 t_{3,4}^5}{25}+7200 t_{1,1} x_1+720 t_{2,2}^2 x_2\\
&-21990 t_{3,3}^2 x_3-36810 t_{3,2} t_{3,4} x_3-48846 t_{3,3} t_{3,4}^2 x_3-\frac{92391}{10} t_{3,4}^4 x_3-37200 t_{3,2} x_3^2-112800 t_{3,3} t_{3,4} x_3^2-41120 t_{3,4}^3 x_3^2\\
&-105450 t_{3,3} x_3^3-118515 t_{3,4}^2 x_3^3-187800 t_{3,4} x_3^4-120000 x_3^5\Biggr)+e^{12 t_{\mu_A}} \biggl(-37320 t_{2,2} t_{3,2}-134976 t_{2,2} t_{3,3} t_{3,4}\\
&-\frac{238568}{5} t_{2,2} t_{3,4}^3-95040 t_{3,2} x_2-218016 t_{3,3} t_{3,4} x_2-\frac{623952}{5} t_{3,4}^3 x_2-142000 t_{2,2} t_{3,3} x_3-123218 t_{2,2} t_{3,4}^2 x_3\\
&-187488 t_{3,3} x_2 x_3+\frac{97668}{5} t_{3,4}^2 x_2 x_3-118240 t_{2,2} t_{3,4} x_3^2+275616 t_{3,4} x_2 x_3^2+12150 t_{2,2} x_3^3+47520 x_2 x_3^3\biggr)\\
&+e^{13 t_{\mu_A}} \left(-310524 t_{1,1} t_{3,3}-\frac{2823906}{5} t_{1,1} t_{3,4}^2-460944 t_{3,3} x_1-\frac{4941216}{5} t_{3,4}^2 x_1-1299552 t_{1,1} t_{3,4} x_3-1821792 t_{3,4} x_1 x_3\right.\\
&-826440 t_{1,1} x_3^2-876240 x_1 x_3^2\biggr)+e^{14 t_{\mu_A}} \left(-767776 t_{2,1} t_{3,4}-\frac{743612}{3} t_{2,2}^2 t_{3,4}-1179120 t_{2,2} t_{3,4} x_2-1710720 t_{3,4} x_2^2\right.\\
&-1018710 t_{2,1} x_3-313395 t_{2,2}^2 x_3+14940 t_{2,2} x_2 x_3+204930 x_2^2 x_3\biggr)+e^{15 t_{\mu_A}} (-1454700 t_{1,1} t_{2,2}-2521080 t_{2,2} x_1\\
&-2105730 t_{1,1} x_2-6105780 x_1 x_2)+e^{16 t_{\mu_A}} \left(2165568 t_{3,3}^2+3609984 t_{3,2} t_{3,4}+\frac{40336272}{5} t_{3,3} t_{3,4}^2+\frac{54600584 t_{3,4}^4}{25}\right.\\
&+5544810 t_{3,2} x_3+19824522 t_{3,3} t_{3,4} x_3+8788881 t_{3,4}^3 x_3+16810920 t_{3,3} x_3^2+21775608 t_{3,4}^2 x_3^2+34558560 t_{3,4} x_3^3\\
&+23399400 x_3^4\biggr)+e^{18 t_{\mu_A}} \left(18291996 t_{2,2} t_{3,3}+\frac{153341514}{5} t_{2,2} t_{3,4}^2+31220640 t_{3,3} x_2+63970560 t_{3,4}^2 x_2\right.\\
\end{align*}

\begin{align*}
&+44706608 t_{2,2} t_{3,4} x_3+13201020 t_{3,4} x_2 x_3+14586360 t_{2,2} x_3^2-20908800 x_2 x_3^2\biggr)+e^{19 t_{\mu_A}} (130053624 t_{1,1} t_{3,4}\\
&+234511560 t_{3,4} x_1+126689040 t_{1,1} x_3+173462400 x_1 x_3)+e^{20 t_{\mu_A}} \left(70426400 t_{2,1}+\frac{107803060 t_{2,2}^2}{3}+154487520 t_{2,2} x_2\right.\\
&+188179200 x_2^2\Biggr)+e^{22 t_{\mu_A}} \left(-322723200 t_{3,2}-1658859840 t_{3,3} t_{3,4}-1104177104 t_{3,4}^3-1751296950 t_{3,3} x_3\right.\\
&\left.-2789043165 t_{3,4}^2 x_3-3852264600 t_{3,4} x_3^2-2730360000 x_3^3\right)+e^{24 t_{\mu_A}} (-5796453600 t_{2,2} t_{3,4}-10520500320 t_{3,4} x_2\\
&-3593143350 t_{2,2} x_3-1178290080 x_2 x_3)+e^{25 t_{\mu_A}} (-8480495040 t_{1,1}-14971717440 x_1)+e^{28 t_{\mu_A}} \bigl(98859098880 t_{3,3}\\
&\left.+229725744576 t_{3,4}^2+338692970880 t_{3,4} x_3+212350641900 x_3^2\right)+e^{30 t_{\mu_A}} (321584010960 t_{2,2}+508631389200 x_2)\\
&+e^{34 t_{\mu_A}} (-20567085718800 t_{3,4}-13652202240000 x_3)+651652721212800 e^{40 t_{\mu_A}},
\end{align*}

\begin{align*}
&\phi^3_{99}=e^{t_{\mu_A}} \left(-t_{1,1} t_{2,1}-\frac{1}{6} t_{1,1} t_{2,2}^2-2 t_{2,1} x_1-\frac{1}{3} t_{2,2}^2 x_1-2 t_{1,1} t_{2,2} x_2-4 t_{2,2} x_1 x_2-3 t_{1,1} x_2^2-6 x_1 x_2^2\right)+e^{2 t_{\mu_A}} \biggl(-4 t_{2,2} t_{3,1}\\
&-\frac{2}{5} t_{2,2} t_{3,3}^2-\frac{4}{5} t_{2,2} t_{3,2} t_{3,4}-\frac{2}{25} t_{2,2} t_{3,3} t_{3,4}^2-\frac{1}{750} t_{2,2} t_{3,4}^4-6 t_{3,1} x_2-\frac{3}{5} t_{3,3}^2 x_2-\frac{6}{5} t_{3,2} t_{3,4} x_2-\frac{3}{25} t_{3,3} t_{3,4}^2 x_2-\frac{1}{500} t_{3,4}^4 x_2\\
&\left.+12 t_{3,2} x_2 x_3+\frac{24}{5} t_{3,3} t_{3,4} x_2 x_3+\frac{8}{25} t_{3,4}^3 x_2 x_3+18 t_{3,3} x_2 x_3^2+\frac{27}{5} t_{3,4}^2 x_2 x_3^2+24 t_{3,4} x_2 x_3^3+30 x_2 x_3^4\right)+e^{3 t_{\mu_A}} \biggl(-21 t_{1,1} t_{3,2}\\
&-\frac{42}{5} t_{1,1} t_{3,3} t_{3,4}-\frac{14}{25} t_{1,1} t_{3,4}^3-18 t_{3,2} x_1-\frac{36}{5} t_{3,3} t_{3,4} x_1-\frac{12}{25} t_{3,4}^3 x_1-18 t_{1,1} t_{3,3} x_3-\frac{27}{5} t_{1,1} t_{3,4}^2 x_3-24 t_{1,1} t_{3,4} x_3^2\\
&-30 t_{1,1} x_3^3\biggr)+e^{4 t_{\mu_A}} \left(-52 t_{2,1} t_{3,3}+\frac{10}{3} t_{2,2}^2 t_{3,3}-\frac{78}{5} t_{2,1} t_{3,4}^2+t_{2,2}^2 t_{3,4}^2-14 t_{2,2} t_{3,3} x_2-\frac{21}{5} t_{2,2} t_{3,4}^2 x_2-48 t_{3,3} x_2^2\right.\\
&\left.-\frac{72}{5} t_{3,4}^2 x_2^2-48 t_{2,1} t_{3,4} x_3+8 t_{2,2}^2 t_{3,4} x_3-48 t_{2,2} t_{3,4} x_2 x_3-60 t_{2,1} x_3^2+10 t_{2,2}^2 x_3^2-90 t_{2,2} x_2 x_3^2\right)+e^{5 t_{\mu_A}} (23 t_{1,1} t_{2,2} t_{3,4}\\
&+78 t_{2,2} t_{3,4} x_1+30 t_{1,1} t_{3,4} x_2-36 t_{3,4} x_1 x_2+90 t_{1,1} t_{2,2} x_3-120 t_{1,1} x_2 x_3-720 x_1 x_2 x_3)+e^{6 t_{\mu_A}} \biggl(84 t_{1,1}^2+144 t_{2,1} t_{2,2}\\
&+4 t_{2,2}^3+204 t_{3,2} t_{3,3}-36 t_{3,1} t_{3,4}+78 t_{3,3}^2 t_{3,4}+54 t_{3,2} t_{3,4}^2+\frac{146}{5} t_{3,3} t_{3,4}^3+\frac{81 t_{3,4}^5}{50}+288 t_{1,1} x_1+90 t_{2,1} x_2+273 t_{2,2}^2 x_2\\
&+432 t_{2,2} x_2^2+144 t_{3,3}^2 x_3+144 t_{3,2} t_{3,4} x_3+144 t_{3,3} t_{3,4}^2 x_3+\frac{84}{5} t_{3,4}^4 x_3+180 t_{3,2} x_3^2+264 t_{3,3} t_{3,4} x_3^2+\frac{312}{5} t_{3,4}^3 x_3^2+240 t_{3,3} x_3^3\\
&+72 t_{3,4}^2 x_3^3\biggr)+e^{8 t_{\mu_A}} \left(2324 t_{2,2} t_{3,2}+\frac{3868}{5} t_{2,2} t_{3,3} t_{3,4}+\frac{1138}{75} t_{2,2} t_{3,4}^3+4080 t_{3,2} x_2+2946 t_{3,3} t_{3,4} x_2+503 t_{3,4}^3 x_2\right.\\
&\left.+3552 t_{2,2} t_{3,3} x_3+\frac{5568}{5} t_{2,2} t_{3,4}^2 x_3-1440 t_{3,3} x_2 x_3-96 t_{3,4}^2 x_2 x_3+5436 t_{2,2} t_{3,4} x_3^2-2970 t_{3,4} x_2 x_3^2+6720 t_{2,2} x_3^3-6000 x_2 x_3^3\right)\\
&+e^{9 t_{\mu_A}} \left(12783 t_{1,1} t_{3,3}+\frac{40749}{10} t_{1,1} t_{3,4}^2+18270 t_{3,3} x_1+6633 t_{3,4}^2 x_1+19242 t_{1,1} t_{3,4} x_3+25440 t_{1,1} x_3^2\right)\\
&+e^{10 t_{\mu_A}} \left(36080 t_{2,1} t_{3,4}-\frac{8096}{3} t_{2,2}^2 t_{3,4}+39220 t_{2,2} t_{3,4} x_2+82320 t_{3,4} x_2^2+48000 t_{2,1} x_3-8960 t_{2,2}^2 x_3+10380 t_{2,2} x_2 x_3\right)\\
&+e^{11 t_{\mu_A}} (20699 t_{1,1} t_{2,2}+71310 t_{2,2} x_1+83568 t_{1,1} x_2+125856 x_1 x_2)+e^{12 t_{\mu_A}} \left(125856 t_{3,1}-\frac{384372 t_{3,3}^2}{5}-\frac{610044}{5} t_{3,2} t_{3,4}\right.\\
&-\frac{2983272}{25} t_{3,3} t_{3,4}^2-\frac{1840531 t_{3,4}^4}{125}+72000 t_{3,2} x_3-196560 t_{3,3} t_{3,4} x_3-70152 t_{3,4}^3 x_3+100500 t_{3,3} x_3^2-100230 t_{3,4}^2 x_3^2\\
&+204000 t_{3,4} x_3^3+450000 x_3^4\Biggr)+e^{14 t_{\mu_A}} \left(-1057056 t_{2,2} t_{3,3}-\frac{2936184}{5} t_{2,2} t_{3,4}^2-1362384 t_{3,3} x_2-\frac{9687726}{5} t_{3,4}^2 x_2\right.
\end{align*}

\begin{align*}
&-1909968 t_{2,2} t_{3,4} x_3+77952 t_{3,4} x_2 x_3-2989260 t_{2,2} x_3^2+427410 x_2 x_3^2\biggr)+e^{15 t_{\mu_A}} (-4862592 t_{1,1} t_{3,4}-6648192 t_{3,4} x_1\\
&-6363810 t_{1,1} x_3)+e^{16 t_{\mu_A}} \left(-7993620 t_{2,1}-599310 t_{2,2}^2-11865582 t_{2,2} x_2-17087760 x_2^2\right)+e^{18 t_{\mu_A}} \bigl(13279500 t_{3,2}\\
&\left.+58919244 t_{3,3} t_{3,4}+\frac{118070466 t_{3,4}^3}{5}+16506720 t_{3,3} x_3+46114704 t_{3,4}^2 x_3+7337760 t_{3,4} x_3^2-79182000 x_3^3\right)\\
&+e^{20 t_{\mu_A}} (289168248 t_{2,2} t_{3,4}+582458130 t_{3,4} x_2+337148640 t_{2,2} x_3-11652000 x_2 x_3)+e^{21 t_{\mu_A}} (598087680 t_{1,1}\\
&+640800000 x_1)+e^{24 t_{\mu_A}} \left(-4658254068 t_{3,3}-\frac{42954595602 t_{3,4}^2}{5}-6772067424 t_{3,4} x_3+3614130720 x_3^2\right)\\
&+e^{26 t_{\mu_A}} (-25117901460 t_{2,2}-36050857920 x_2)+e^{30 t_{\mu_A}} (1001073115680 t_{3,4}+176276934000 x_3)-34893343710720 e^{36 t_{\mu_A}}.
\end{align*}
}
%%%%%%%%%%%%%%%%%%%%%%%%%%%%%%%%%%%%%%%%%%%%%%%%%%%%%%%%%%%%%%%%%%%%%%%%%%%%%%%
\subsection{$\psi$ for $A=(2,3,5)$}
\
{\tiny
\begin{align*}
\psi ^1=&x_1+e^{3 t_{\mu_A}} \left(18 t_{3,3}+\frac{27 t_{3,4}^2}{5}+48 t_{3,4} x_3+90 x_3^2\right)+e^{5 t_{\mu_A}} (-210 t_{2,2}-360 x_2)\\
&+e^{9 t_{\mu_A}} (13230 t_{3,4}-12000 x_3)-2708910 e^{15 t_{\mu_A}},\\
\psi ^2=&-48 e^{4 t_{\mu_A}} t_{3,4}+40800 e^{10 t_{\mu_A}},\\
\psi ^3=&-1440 e^{6 t_{\mu_A}}.
\end{align*}
}
%%%%%%%%%%%%%%%%%%%%%%%%%%%%%%%%%%%%%%%%%%%%%%%%%%%%%%%%%%%%%%%%%%%%%%%%%%%%%%%%%%%%%%%%%%%%
%%%%%%%%%%%%%%%%%%%%%%%%%%%%%%%%%%%%%%%%%%%%%%%%%%%%%%%%%%%%%%%%%%%%%%%%%%%%%%%%%%%%%%%%%%%%

\end{document}